\documentclass{svjour2}                    
\smartqed  
\usepackage{graphicx}
\usepackage{amsmath}  \usepackage{amsfonts}
\usepackage{amssymb} \usepackage{latexsym}
\usepackage{color}
\usepackage{graphicx}\usepackage{subfigure}
\usepackage{epstopdf}
\usepackage[title]{appendix}
\allowdisplaybreaks
\usepackage{psfrag}
\allowdisplaybreaks \allowdisplaybreaks[0]
\newtheorem{thm}{Theorem}[section]
\numberwithin{equation}{section}
\numberwithin{figure}{section}
\numberwithin{table}{section}
\numberwithin{thm}{section}
\journalname{Journal of Statistical Physics}
\begin{document}

\title{Globally hyperbolic moment model of arbitrary order  for
	one-dimensional special relativistic Boltzmann equation\thanks{
	This work was partially supported by
	the National Natural Science Foundation
	of China (Nos.  91330205 \& 11421101) and the National Key Research and Development Program of China
(No. 2016YFB0200603).
}
}


\author{Yangyu Kuang         \and
        Huazhong Tang 
}


\institute{Y.Y. Kuang \at LMAM, School of Mathematical Sciences\\
Peking University\\ Beijing 100871, P.R. China\\
              \email{kyy@pku.edu.cn}           
           \and
           H.Z. Tang \at
            HEDPS, CAPT \& LMAM\\ School of Mathematical Sciences\\ Peking University\\
            Beijing 100871, P.R. China\\
           \at School of Mathematics and Computational Science\\
           Xiangtan University\\  Xiangtan 411105,  Hunan Province, P.R. China\\
               \email{hztang@math.pku.edu.cn}
}

\date{Received: date / Accepted: date}

\maketitle

\begin{abstract}
This paper extends the model reduction method by the operator projection
to the one-dimensional special relativistic Boltzmann equation.
 The derivation of arbitrary order globally hyperbolic moment system is built on
our careful study of  two families of the complicate
'Grad type orthogonal polynomials  depending on a parameter.
We derive their   recurrence relations, calculate
their derivatives with respect to the independent variable and parameter  respectively,
and study   their zeros and  coefficient matrices in the  recurrence formulas.
%
%
%
Some properties of the moment system are also proved. They include
the eigenvalues and their bound as well as eigenvectors, hyperbolicity, characteristic fields, linear stability,  and Lorentz covariance.
A semi-implicit numerical scheme is presented to solve a Cauchy problem of our
hyperbolic moment system in order to verify the convergence behavior of the moment method.
The results show that the solutions of our hyperbolic moment system {converge} to the solution of the special relativistic Boltzmann equation as the order of the hyperbolic moment system increases.

\end{abstract}

\section{Introduction}
\label{sec:intro}

The beginning of the relativistic kinetic theory goes back to 1911 when
an equilibrium distribution function was derived for a relativistic gas \cite{Juttner1911}.
Thirty years later, the covariant formulation of the {relativistic} Boltzmann equation was proposed in \cite{LM1940} to describe the statistical behavior of a thermodynamic system not in thermodynamic equilibrium.
The   transport coefficients {were} determined from the  Boltzmann equation by using the Chapman-Enskog methodology  in \cite{Israel1963}.
Different from a non-relativistic monatomic gas,  a relativistic gas has
a bulk viscosity. It has called the attention of many researchers to a number of applications of  this theory:  the
effect of neutrino viscosity on the evolution of the universe, the study of galaxy
formation, neutron stars, and controlled thermonuclear fusion, etc.
%
%
The readers are referred to the monographs \cite{RB:2002,RK:1980} for  more detailed
descriptions. 

The relativistic kinetic  theory  is attracting increasing attention  in recent years,
but  it has been used relatively sparsely to model phenomenological matter in comparison to fluid models.
In the non-relativistic case,
the kinetic theory has been studied intensively as a mathematical subject during several decades, and   also played an important role from an engineering point of view,
see e.g. \cite{RB:1988,Cha:1970}.
From the Boltzmann equation one could determine the  distribution function
hence the transport coefficients of  gases, however this task was not so easy.
 Hilbert showed that an  approximate solution of the integro-differential
equation  could be {obtained} from a power series expansion of a parameter
(being proportional to the mean free path).
Chapman and Enskog calculated   independently
 the transport coefficients for gases whose molecules interacted according to any kind of spherically symmetric potential function.
Another method  proposed  by Grad \cite{GRAD:1949,GRAD2:1949} is to
expand the distribution
function in terms of tensorial Hermite polynomials and introduce
the balance equations corresponding to higher order moments of the distribution function.
%
%
The crucial ingredient of the  Chapman-Enskog  method is
the assumption that in the hydrodynamic regime
the distribution function can be expressed as a function of the hydrodynamic variables and their gradients.
%
The Chapman-Enskog  method has been
extended to the relativistic cases, see e.g. \cite{CHA2:2008,CHA:2012,CHA:2008,NHH:1983,NHH:1985}.
Unfortunately, it is difficult to  derive the equations of relativistic fluid dynamics from the kinetic theory \cite{DM:2012}.
The moment method can avoid  such difficulty
and is  also generalized to the relativistic cases,
see {e.g. \cite{Ander:1970,IS:1976,IS:1979,IS2:1979,Kranys:1972,Stewart:1977}}.
However, the moment method cannot reflect the influence of the Knudsen number.
Combining the Chapman-Enskog method with
the moment method has been attempted  \cite{DM:2012,RDH:2013}.
%

It is difficult to derive  the relativistic moment system  of higher order
since the family of orthogonal polynomials can not be found easily.
Several authors  have tried to construct the family of orthogonal polynomials analogous to the Hermite polynomials, see e.g \cite{GRADP:1974,RK:1980}.
%
Their application can be found in \cite{DM:2012,RDH:2013,PRO:1998}.
%
%
Unfortunately,
there is no explicit expression of the moment systems
 if the order of the moment system is larger than $3$.
%
Moreover, { the hyperbolicity of  existing general moment systems is not proved,
 even for the  second order moment system (e.g. the general Israel and Stewart system).
 For a special  case with heat conduction and no viscosity,
 Hiscock and Lindblom  proved that the Israel and Stewart moment system in  the
  Landau frame was globally hyperbolic and linearly stable,
  but they also showed that the Israel and Stewart moment system in
     the Eckart frame was not globally hyperbolic and linearly stable.
     The readers are referred to \cite{NHH:1987,NHH:1988,NNH:1989}.
 Following the approach used in \cite{NHH:1987,NHH:1988},
 it is easy to show that  the above conclusion is not true if the viscosity exists,   that is,  the Israel and Stewart moment system in the Landau frame
 is not globally hyperbolic too if  the viscosity exists. There does not exist any result on the hyperbolicity or  loss of hyperbolicity  of  (existing) general higher-order moment systems for the relativistic kinetic equation. Such
proof is very difficult and challenging.
The loss of hyperbolicity will   cause the solution blow-up when the distribution is far away from the equilibrium state. Even for the non-relativistic case,  increasing the number of moments could not  avoid  such blow-up \cite{FAR:2012}.}

  {Up to now, there has been some latest progress on the Grad moment method in the non-relativistic case.}
  A regularization  was presented in  \cite{1DB:2013} for the 1D Grad moment system to achieve global hyperbolicity.
It was based on the observation that the characteristic polynomial of the Jacobian of the flux
in Grad's moment system is independent of the intermediate moments,
 and  further extended to the multi-dimensional case \cite{HY:2014,HG:2014}.
 The quadrature based projection methods were used
to derive hyperbolic PDE systems for the solution of the Boltzmann equation
  \cite{QI:2014,HG2:2014}  
 by  using some quadrature rule instead of the exact integration.
In the 1D case,  it is similar to the  regularization in  \cite{1DB:2013}. 
Those contributions  led to  well understanding the hyperbolicity of  the Grad moment systems.
Based on the operator projection,  a general framework of model reduction technique was recently presented  in \cite{MR:2014}.
It  projected the time and space derivatives in the kinetic equation   into a finite-dimensional weighted polynomial space synchronously,
and might give most of the existing moment systems mentioned above.
 %
%
The aim of this paper is to extend the model reduction method by the operator projection  \cite{MR:2014}
to the one-dimensional special relativistic Boltzmann equation and
 derive corresponding globally hyperbolic moment system of arbitrary order.
 The key  is  to   choose the weight function and define  the  polynomial spaces and their basis as well as the projection operator.
  The theoretical {foundations} of our moment method {are} the properties
  of two families of the complicate Grad type orthogonal polynomials depending on a parameter.

The paper is organized as follows.
Section \ref{sec:RB} introduces the special relativistic Boltzmann equation
and some macroscopic quantities defined via the kinetic theory.
 Section \ref{sec:orth} gives  two families of orthogonal polynomials dependent on a parameter,  
and studies their properties: recurrence relations, derivative relations with respect to the variable and the parameter, zeros, and the eigenvalues and eigenvectors of the recurrence matrices. Section \ref{sec:moment} derives the moment system
of the special relativistic Boltzmann equation and Section \ref{sec:prop} studies its properties:  the eigenvalues and its
bound as well as eigenvectors, hyperbolicity,  characteristic fields,   linear stability, and Lorentz covariance.
 Section \ref{sec:NE}
presents a semi-implicit numerical scheme and conducts a numerical experiment to check the convergence of the proposed hyperbolic moment system.
Section \ref{sec:conclud} concludes the paper. 
{
To make the main message of the paper less dilute, all proofs of theorems, lemmas and corollaries in Sections \ref{sec:RB}-\ref{sec:NE}
are given in the Appendices \ref{sec:App2}-\ref{sec:App6} respectively.}

\section{Preliminaries and notations}  
\label{sec:RB}

In the special relativistic kinetic theory of gases \cite{RB:2002},
a microscopic gas particle of rest  mass $m$ is characterized by the $(D+1)$ space-time
coordinates $(x^{\alpha})=(x^0,\vec{x})$ and momentum $(D+1)$-vectors $(p^{\alpha})=(p^{0},\vec{p})$, where
$x^0=ct$, $c$ denotes the speed of light in vacuum, and $t$ and $\vec{x}$ are the time and $D$-dimensional spatial coordinates,
 respectively. 
Besides  the contravariant notation (e.g. $p^{\alpha}$),  the covariant notation such as $p_{\alpha}$
 will be also used in the following and the covariant $p_{\alpha}$ is related to  the contravariant  $p^{\alpha}$ by
 \[
 p_{\alpha}=g_{\alpha\beta}p^{\beta}, \quad p^{\alpha}=g^{\alpha\beta}p_{\beta},
 \]
%
%
where $(g^{\alpha\beta})$ denotes the  Minkowski  space-time metric tensor
and is chosen as \\$(g^{\alpha\beta})={\rm diag}\{1,-\vec{I}_{D}\}$,
$\vec{I}_{D}$ is the $D\times D$ identity matrix, $(g_{\alpha\beta})$ denotes
the inverse of $(g^{\alpha\beta})$, and the Einstein summation convention over repeated indices is used.
%
For a free relativistic particle,  one has the relativistic energy-momentum
relation (aka ``on-shell'' or ``mass-shell'' condition)
$E^2-\vec p^2 c^2=m^2 c^4$. If putting $p^0= c^{-1}E=\sqrt{\vec{p}^2+m^2c^2}$,
then the ``mass-shell'' condition {can} be rewritten as
$p^{\alpha}p_{\alpha}=m^2c^2$.

As in the non-relativistic case, the relativistic Boltzmann equation
describes the evolution of the one-particle distribution function of an ideal gas in the phase space
spanned by the space-time coordinates  $(x^{\alpha})$ and momentum { (D+1)}-vectors of particles
 $(p^{\alpha})$.
The one-particle distribution function depends only
on $(\vec x,\vec p, t)$ and is defined in such a way that
$f(\vec x, \vec p, t) d^D\vec x d^D\vec p$
gives the number of particles at time $t$   in the volume element $d^D \vec x d^D\vec p $.
%
For a single gas the Boltzmann equation reads  \cite{RB:2002}
\begin{equation}
\label{eq:Boltz}
p^{\alpha}\frac{\partial f}{\partial x^{\alpha}}=Q(f,f),
\end{equation}
where the collision term $Q(f,f)$
  depends  on the product of the distribution functions of two particles at collision, e.g.
$$
Q(f,f)=\int_{\mathbb{R}^{{D}}}\int_{\mathbb{S}^{{D-1}}_{+}}\left(f_{*}'f'-f_{*}f\right)B d{\Omega} \frac{d^{{D}}\vec{p}_{*}}{p_{*}^{0}},
$$
where $f$ and $f_*$ are the {distributions} depending on  the momenta before a collision, while
 $f'$ and $f_*'$  depend on  the momenta   after the collision,
 $d\Omega$ denotes the element of the solid angle,
the collision kernel $B=\sigma \sqrt{(p_{*}^{\alpha}p_{\alpha})^{2}-m^2c^2}$ for a single non degenerate gas (e.g. electron gas),
and $\sigma$ denotes the differential cross section of collision,.
The collision term  satisfies
\begin{equation}
\label{eq:colcon}
    \int_{\mathbb{R}^{{D}}}Q(f,f)\frac{d^{D}\vec{p}}{p^{0}}=0, \quad \int_{\mathbb{R}^{{D}}}p^{\alpha}Q(f,f)\frac{d^{D}\vec{p}}{p^{0}}=0,
\end{equation}
so that  $1$ and $p^{\alpha}$ are called {\em collision invariants}.
%
Moreover,  the Boltzmann equation \eqref{eq:Boltz} should satisfy
the entropy dissipation relation (in the sense of classical statistics)  
\[
\int_{\mathbb{R}^{D}}Q(f,f)\ln(f)\frac{d^D\vec p}{p^{0}}\leq0,
\]
where the  equal sign corresponds to the local thermodynamic equilibrium.

In kinetic theory the macroscopic description of gas  can
be represented by { the first and second} moments of the distribution function $f$, namely,
the partial particle {(D+1)}-flow $N^\alpha$
and the partial energy-momentum tensor
$T^{\alpha\beta}$,  which are defined by
%
\begin{equation}
\label{eq:NTab}
  N^{\alpha}=c\int_{{\mathbb{R}^{D}}} p^{\alpha}f\frac{d^{D}\vec{p}}{p^{0}},\quad T^{\alpha\beta}=c\int_{{\mathbb{R}^{D}}} p^{\alpha}p^{\beta}f\frac{d^{D}\vec{p}}{p^{0}}.
\end{equation}
They  {can} be decomposed into the following forms (i.e. the Landau-Lifshitz decomposition)
\begin{align}
N^{\alpha}&={m^{-1}}\rho U^{\alpha}+n^{\alpha},\label{eq:Ndiv} \\
T^{\alpha\beta}&=c^{-2}\varepsilon U^{\alpha}U^{\beta}-\Delta^{\alpha\beta}(P_{0}+\Pi)+ \pi^{\alpha\beta},\label{eq:Tabdiv}
\end{align}
where $(U^{\alpha})=\left(\gamma(\vec{u}) c, \gamma(\vec{u})\vec u\right) $
denotes the macroscopic velocity $(D+1)$-vector of gases,
$\gamma(\vec{u})=(1-c^{-2}|\vec{u}|^{2})^{-\frac{1}{2}}$
is the Lorentz factor, $\Delta^{\alpha\beta}$ is   defined by
\begin{equation}
\label{eq:sym-tensor}
\Delta^{\alpha\beta}:=g^{\alpha\beta}-c^{-2}U^{\alpha}U^{\beta},
\end{equation}
which is symmetric and the projector onto the $D$-dimensional subspace orthogonal to $U^\alpha$, that is,
satisfies $\Delta^{\alpha\beta}U_{\beta}=0$.
Here,  the mass density $\rho$, the particle-diffusion current $n^{\alpha}$, the energy density $\varepsilon$, the shear-stress tensor $\pi^{\alpha\beta}$, and the sum of thermodynamic pressure $P_{0}$ and bulk viscous pressure $\Pi$ are defined and related to the distribution $f$ by
\begin{equation}
\label{eq:variable}
\begin{aligned}
&\rho:=c^{-2}{m}U_{\alpha}N^{\alpha}=c^{-1}{m}\int_{\mathbb{R}^{D}} Ef\frac{ d^{D}\vec{p} }{p^{0}},\\
&n^{\alpha}:=\Delta^{\alpha}_{\beta}N^{\beta}=c\int_{\mathbb{R}^{D}} p^{<\alpha>}f\frac{d^{D}\vec{p}}{p^{0}},\\
&\varepsilon:=c^{-2}U_{\alpha}U_{\beta}T^{\alpha\beta}=c^{-1}\int_{\mathbb{R}^{D}} E^2f\frac{ d^{D}\vec{p}}{p^{0}},\\
&\pi^{\alpha\beta}:= \Delta_{\mu\nu}^{\alpha\beta} T^{\mu\nu}=c\int_{\mathbb{R}^{D}} p^{<\alpha \beta>}f\frac{ d^{D}\vec{p} }{p^{0}},\\
&P_{0}+\Pi:=-D^{-1}\Delta_{\alpha\beta}T^{\alpha\beta}=D^{-1}c^{-1}\int_{\mathbb{R}^{D}} (E^2-{m^2}c^4)f\frac{ d^{D}\vec{p}}{p^{0}},
\end{aligned}
\end{equation}
where  $E:=U_{\alpha}p^{\alpha}$ here and hereafter, $p^{<\alpha>}:=\Delta_{\beta}^{\alpha}p^{\beta}$,
$p^{<\alpha \beta>}:=\Delta^{\alpha\beta}_{\mu\nu}{p^{\mu}p^{\nu}}$,
 and
 \[
 \Delta^{\alpha\beta}_{\mu\nu}:= \frac{1}{2}\left( \Delta_{\mu}^{\alpha}\Delta_{\nu}^{\beta}+ \Delta_{\mu}^{\beta}\Delta_{\nu}^{\alpha}-2D^{-1}\Delta_{\mu\nu}\Delta^{\alpha\beta} \right).
 \]
{ It is obvious to obtain
\begin{equation}
\label{eq:divpp}
U_{\alpha}p^{<\alpha>}=0.
\end{equation}
}
 It is not {difficult} to verify the following identity
\begin{equation}
\label{eq:divp}
p^{\alpha}=c^{-2}EU^{\alpha}+p^{<\alpha>}.
\end{equation}

Multiplying the special relativistic Boltzmann equation \eqref{eq:Boltz}
by $1$ and $p^{\alpha}$ respectively, integrating both sides over $\mathbb{R}^{D}$ in terms of
$\vec{p}$, and using \eqref{eq:colcon} gives  the following conservation laws
\begin{equation}
\label{eq:Tab}
    \partial_{\alpha} N^{\alpha}=0,\quad
    \partial_{\alpha} T^{\alpha\beta}=0.
\end{equation}

 \begin{remark}
  It is common to choose $U^\alpha$ as the velocity of either energy transport
   (the Landau-Lifshitz frame) \cite{Lau:1949}), i.e.
  \begin{equation}
  \label{eq:landau}
  U_{\beta}T^{\alpha\beta}=\varepsilon U^{\alpha},
  \end{equation}
  i.e.
   	\begin{equation}
   	\label{eq:condition-222}
   	\Delta^{\alpha}_{\beta}T^{\beta\gamma}U_{\gamma}=c\int_{{\mathbb{R}^{D}}} Ep^{<\alpha>}f\frac{d^D\vec p}{p^{0}}=0,
   	\end{equation}
  or particle transport (the Eckart frame) \cite{EC:1940}), i.e.
  in which the velocity is specified by the flow of particles
  	 \[
  	 N^{\alpha}={m^{-1}}\rho U^{\alpha},
  	 \]
  	   i.e.
  	   \begin{equation*}
  	   \Delta^{\alpha}_{\beta}{N^{\beta}}=c \int_{{\mathbb{R}^{D}}}   \Delta^{\alpha}_{\beta}{p^{\beta}}    f \frac{d^D\vec p}{p^{0}}=0.
  	   \end{equation*}
The former {{can}} be applied to multicomponent gas while the latter is only used for single component gas.
 This work will be done in the Landau-Lifshitz frame \eqref{eq:landau}.
 \end{remark}

 \begin{remark}\label{remark2.2}
 	At the local thermodynamic equilibrium, $n^\alpha$, $\Pi$, and $\pi^{\alpha\beta}$
 	will be zero.
 \end{remark}

 \begin{remark}
 	In order to simplify  the collision term, several simple collision models have been proposed, see \cite{RB:2002}.
 	Similar to  the  BGK (Bhatnagar-Gross-Krook) model  
 	in the non-relativistic theory,
 	two simple relativistic  collision models  are the Marle model \cite{MM:1965}
 	\begin{equation}
 	Q(f,f)=-\frac{m}{\tau}(f-f^{(0)}),\label{eq:colmarle}
 	\end{equation}
 	and
 	the Anderson-Witting model \cite{AW:1974}
 	\begin{equation}
 	Q(f,f)=-\frac{U_{\alpha}p^{\alpha}}{\tau { c^2}}(f-f^{(0)}),\label{eq:colAW}
 	\end{equation}
 	where  $f^{(0)}=f^{(0)}(\vec x, \vec p, t)$  denotes the distribution function at  the local thermodynamic equilibrium,
 	and $\tau$ is the relaxation time and  may rely on $\rho$, $\theta$.
 In the non-relativistic limit,  both   models \eqref{eq:colmarle} and  \eqref{eq:colAW} tend to the BGK model.
 However, 	the Marle model \eqref{eq:colmarle} does not satisfy the constraints of the collision terms in \eqref{eq:colcon}. 
{The relaxation time $\tau$  can be defined by
\[
\tau=\frac{1}{n\pi d^2\bar{g}},
\]
where $n$ denotes the particle number density, $d$ denotes the diameter of gas particles,
and $\bar{g}$ is proportional to the mean relative speed $\bar{\xi}$  between two particles, e.g.
$\bar{g}=\sqrt{2}\bar{\xi}$ or $\bar{\xi}$ \cite{RB:2002}.
In the non-relativistic case, $\bar{\xi}=4\sqrt{\frac{kT}{\pi m}}$,
but
the expression of $\bar{\xi}$ in relativistic case is very complicate,
see Section 8.2 of book \cite{RB:2002}.
Usually,  $\bar{\xi}$ or $\bar g$ is suitably approximated,
 for example, $\bar{g}\approx c$ (that is, $\bar{g}$ is
  approximated by using the ultra-relativistic limit).
 Under such simple approximation, one has
 \[
\tau\approx \frac{1}{n\pi d^2 c}=\frac{m}{\rho \pi d^2 c}.
\]}
 \end{remark}


This paper will only consider the one-dimensional  form of
 relativistic Boltzmann equation \eqref{eq:Boltz}.
In this case, the vector notations $\vec{x}$ and $\vec{p}$ will be replaced  with  $x$ or $x^1$ and $p$
or $p^{1}$, respectively,
the  Greek indices $\alpha$ and $\beta$ run from 0 to 1,
and \eqref{eq:Boltz}  reduces to the following form
\begin{equation}
\label{eq:1DBoltz}
p^{0}\frac{\partial f}{\partial ct}+p^{1}\frac{\partial f}{\partial x}=Q(f,f),\quad t\in\mathbb{R}^{+},\ x\in\mathbb{R}.
\end{equation}
In the $1$D case, the shear-stress tensor $\pi^{\alpha\beta}$ disappears even though
the local-equilibrium is departed from, and the
 local-equilibrium distribution $f^{(0)}$ {can} be explicitly given by
\begin{equation}
  \label{eq:equm1}
   f^{(0)}=\rho g^{(0)}, \quad g^{(0)}=\frac{1}{2{m^{2}}cK_{1}(\zeta)}\exp\left(-\zeta E\right),
\end{equation}
which is {like} the Maxwell-J\"{u}ttner distribution  {\cite{RB:2002} for the case of $D=3$ and Maxwell gas}
\begin{equation*}{
   f^{(0)}=\rho g^{(0)}, \quad g^{(0)}=\frac{\zeta}{4\pi m^{4}c^{3}K_{2}(\zeta)}\exp\left(-\zeta E\right),}
\end{equation*}
and
obeys the common prescription that the mass density $\rho$ and energy  density $\varepsilon$
are completely determined by the local-equilibrium distribution $f^{(0)}$ alone,
that is,
\begin{equation}
\label{eq:condition}
\rho =\rho_{0},\quad
\varepsilon=\varepsilon_{0}.
\end{equation}
In \eqref{eq:equm1},  $\zeta={(k_{B}T)^{-1}(mc^2)}$ is  the ratio between the particle rest energy $mc^2$
and the thermal energy of the gas $k_BT$, $k_B$ denotes the Boltzmann constant,  $T$ is the {thermodynamic} temperature,
and $K_{n}(\zeta)$ denotes the modified Bessel function of the second kind, defined by
\begin{equation}
\label{eq:bessel}
   K_{n}(\zeta)=\int_{0}^{\infty}\cosh(n\vartheta)\exp(-\zeta\cosh\vartheta)d\vartheta,
\end{equation}
satisfying the recurrence relation
\begin{equation}
\label{eq:besselrec}
   K_{n+1}(\zeta)=K_{n-1}(\zeta)+2n\zeta^{-1}K_{n}(\zeta).
\end{equation}
For $\zeta\gg 1$ the particles behave as non-relativistic,
and for $\zeta\ll 1$ they behave as ultra-relativistic.

Similar to  \eqref{eq:variable}, from the knowledge of
the equilibrium distribution function $f^{(0)}$ it is also possible to determine the values of
some macroscopic variables by
\begin{equation}
 \label{eq:variable0}
    \begin{aligned}
    &\rho_{0}:=c^{-1}{m}\int_{\mathbb{R}} Ef^{(0)}\frac{dp}{p^{0}},\\
    &n^{\alpha}_{0}:=c\int_{\mathbb{R}} p^{<\alpha>}f^{(0)}\frac{dp}{p^{0}}=0,\\
    &\varepsilon_{0}:=c^{-1}\int_{\mathbb{R}} E^2f^{(0)}\frac{dp}{p^{0}}={\rho }c^2\left(G(\zeta)-\zeta^{-1}\right),\\
    &P_{0}:=c^{-1}\int_{\mathbb{R}} (E^2-{m^2}c^4)f^{(0)}\frac{dp}{p^{0}}={m^{-1}}\rho k_BT={\rho c^2\zeta^{-1}},
    \end{aligned}
\end{equation}
where $G(\zeta):=K_{1}^{-1}(\zeta)K_{2}(\zeta)$.
%
Now, the conservation laws \eqref{eq:Tab} become
\begin{equation}
\label{eq:conser}
\begin{aligned}
&\frac{\partial \left(\rho U^{0}\right)}{\partial ct}+\frac{\partial \left(\rho U^{1}\right)}{\partial x}=0,\\
&\frac{\partial \left(c^{-2}\rho h U^{0}U^{1}\right)}{\partial ct}+\frac{\partial  \left(c^{-2}\rho h U^{1}U^{1}+P_{0}\right)}{\partial x}=0,\\
&\frac{\partial \left(c^{-2}\rho h U^{0}U^{0}-P_{0}\right)}{\partial ct}+\frac{\partial \left(c^{-2}\rho h U^{0}U^{1}\right)}{\partial x}=0,
\end{aligned}
\end{equation}
where $h:=\rho^{-1}(\varepsilon+P_{0})={c^2}G(\zeta)$ denotes the specific enthalpy.
{They are} just the macroscopic equations of  special relativistic hydrodynamics (RHD). In other words, when $f=f^{(0)}$, the special relativistic Boltzmann equation \eqref{eq:1DBoltz} {can} lead to the RHD equations \eqref{eq:conser}.
{We aim at finding reduced model equations to describe states with $f\neq f^{(0)}$.}
This paper will extend the moment method by operator projection  \cite{MR:2014}
to \eqref{eq:1DBoltz} and derive its arbitrary order moment model in Section \ref{sec:moment}.

Before ending this section,  we discuss the macroscopic variables calculated by
a given distribution $f$, in other words,  for the nonnegative distribution  $f(x,p,t)$,
which is not identically zero, can the physically admissible macroscopic states $\{\rho,u,\theta  ={\zeta}^{-1}\}$
satisfying $\rho>0,|u|<{c}$ and $\theta>0$
be obtained?

\begin{thm}
	\label{thm:admissible}
	For the nonnegative distribution  $f(x,p,t)$,  which is not identically zero,
the   density current $N^{\alpha}$ and  energy-momentum tensor $T^{\alpha\beta}$
calculated by \eqref{eq:NTab}
satisfy
	\begin{equation}
		\label{eq:admissible}
		(T^{00}+T^{11})^2>4(T^{01})^2, \quad N^{0}-{c^{-1}u}N^{1}>0, \quad {c^{-2}\rho^{-1}(T^{00}-c^{-1}uT^{01})}>1,
	\end{equation}
	where the macroscopic velocity 	$u$ is the unique solution satisfying $|u|<{c}$ of the {quadratic} equation
	\begin{equation}
		\label{eq:admu1}
		{T^{01}c^{-2}u^2-(T^{00}+T^{11})c^{-1}u+T^{01}=0},
	\end{equation}
{ which has  a solution satisfying ${|u|<c}$ and
\begin{equation}
\label{eq:admu}
u=\left\{
    \begin{array}{ll}
      \frac{T^{00}+T^{11}-\sqrt{(T^{00}+T^{11})^2-4(T^{01})^2}}{2T^{01}c^{-1}}, & \hbox{$T^{01}\neq0$,} \\
      0, & \hbox{$T^{01}=0$.}
    \end{array}
  \right.
\end{equation}
}
And the positive mass density $\rho$ is calculated by
	\begin{equation}
		\label{eq:admrho}
		\rho={c^{-1}m\frac{N^{0}-c^{-1}u N^{1}}{\sqrt{1-c^{-2}u^2}}}.
	\end{equation}
	Furthermore, the equation
	\begin{equation}
		\label{eq:admT}
		{ G(\theta^{-1})-\theta=c^{-2}\rho^{-1}(T^{00}-c^{-1}uT^{01})},
	\end{equation}
	has a unique  positive solution $\theta$ in the interval $(0,+\infty)$.
\end{thm}

Furthermore, the following conclusion holds.
	\begin{thm}
		\label{lem:admissible1}
		Under the {assumptions} of Theorem \ref{thm:admissible}, the bulk viscous pressure $\Pi$
		satisfies
		\[
		\Pi >  -\rho {c^{2}}\theta.
		\]
	\end{thm}


\begin{remark}
{ The proofs of those theorems are given
 in the Appendix \ref{sec:App2}.}	 Theorem \ref{thm:admissible} provides a recovery procedure
	 of the admissible primitive variables $\rho, u$, and $\theta$
	 from the nonnegative distribution  $f(x,p,t)$ or the given density current $N^{\alpha}$
	 and  energy-momentum tensor $T^{\alpha\beta}$ satisfying \eqref{eq:admissible}.
	    It is  useful in  the derivation of the moment system as well as the numerical scheme.
\end{remark}


{ Before discussing the moment method, we first non-dimensionalize the relativistic Boltzmann equation \eqref{eq:1DBoltz}. Here we only consider the Anderson-Witting model \eqref{eq:colAW}.
If  setting
$$
x=L\hat{x}, ~ p=c\hat{p}, ~ p^{0}=c\hat{p^{0}},~ t=\frac{L}{c}\hat{t},~
g=c\hat{g},~ f=\frac{n_{0}}{c^3}\hat{f},
$$
where  $L$ denotes the macroscopic characteristic length,
 $n_{0}$ and $\theta_{0}=m c^2/k_{B}$ are the reference particle number and  temperature,
 respectively,
then   the 1D relativistic Boltzmann equation \eqref{eq:1DBoltz} with \eqref{eq:colAW}
is  non-dimensionalized as follows
\[
\frac{n_{0}}{c^2L}\left(p^{0}\frac{\partial \hat{f}}{\partial \hat{t}}+p^{1}\frac{\partial \hat{f}}{\partial \hat{x}}\right)
=\frac{n_{0}^{2}\pi d^2}{c^2}\hat{U}_{\alpha}\hat{p}^{\alpha}\hat{\rho}\left(\hat{f}^{(0)}-\hat{f}\right),
\]
or
\[
\hat{p}^{0}\frac{\partial \hat{f}}{\partial \hat{t}}+\hat{p}^{1}\frac{\partial \hat{f}}{\partial \hat{x}}
=n_{0}L\pi d^2\hat{U}_{\alpha}\hat{p}^{\alpha}\hat{\rho}\left(\hat{f}^{(0)}-\hat{f}\right).
\]
Thanks to $K_{n}=\frac{\lambda}{L}=\frac{\tau_{0}c}{L}=\frac{1}{n_{0}L\pi d^2}$, 
the above   equation is
rewritten as
\begin{equation}\label{EQ:0000000001}
\hat{p}^{0}\frac{\partial \hat{f}}{\partial \hat{t}}+\hat{p}^{1}\frac{\partial \hat{f}}{\partial \hat{x}}=
\frac{\hat{\rho}}{K_{n}}\hat{U}_{\alpha}\hat{p}^{\alpha}\left(\hat{f}^{(0)}-
\hat{f}\right).
\end{equation}
Thus, if $\tilde\tau:=\frac{K_{n}}{\hat\rho}$ may be considered as a new ``relaxation time'',
then the  collision term of   relativistic Boltzmann equation \eqref{EQ:0000000001}
has  the same form of non-relativistic BGK model.
For the sake of convenience, in the following, we still use $\tau$, $x$, $t$, $f$, $p$, $p^{0}$, $\rho$ to replace $\tilde\tau$, $\hat x$, $\hat t$, $\hat f$, $\hat p$, $\hat p^{0}$, $\hat \rho$, respectively. }

\section{Two families of orthogonal polynomials}
\label{sec:orth}
This section  introduces two families of orthogonal polynomials dependent on a parameter $\zeta$,
similar to  those given in \cite{GRADP:1974},
 and studies their properties, which will be used in the derivation and discussion of our moment system. {All proofs  are given in the Appendix \ref{sec:App3}.}

If considering
\[
  \omega^{(\ell)}(x;\zeta)=\frac{(x^2-1)^{\ell-\frac{1}{2}}}{K_{1}(\zeta)}\exp(-\zeta x), \ \ell=0,1,
\]
as  the weight functions in the interval $[1,+\infty)$, where  $\zeta\in\mathbb{R}^{+}$ denotes
 a parameter,
then  the inner products with respect to $\omega^{(\ell)}(x;\zeta)$ {can} be introduced as follows
\[
\left(f,g\right)_{\omega^{(\ell)}}:=\int_{1}^{+\infty}f(x)g(x)\omega^{(\ell)}(x;\zeta)dx, \quad f,g\in L^{2}_{\omega^{(\ell)}}[1,+\infty),\ \ell=0,1,
\]
where  $L^{2}_{\omega^{(\ell)}}[1,+\infty):=\left\{f\big|\int_{1}^{+\infty}f(x)^2\omega^{(\ell)}(x;\zeta)dx<+\infty\right\}$.
It is worth noting that
the choice of the weight function $\omega^{(\ell)}(x;\zeta)$
is dependent on the equilibrium distribution {$f^{(0)}(x,p, t)$} in
\eqref{eq:equm1}.

Let $\{P_{n}^{(\ell)}(x;\zeta)\}$, $\ell=0,1$, be two families of standard orthogonal polynomials
with respect to the weight function $\omega^{(\ell)}(x;\zeta)$ in the interval  $[1,+\infty)$, i.e.
\begin{equation}
  \label{eq:P01orth}
     \left(P_{m}^{(\ell)},P_{n}^{(\ell)}\right)_{\omega^{(\ell)}}=\delta_{m,n}, \quad
   \ell=0,1,
\end{equation}
where $\delta_{m,n}$ denotes the Kronecker delta function,
which is equal to 1 if $m=n$, and 0 otherwise.
Obviously,   $\{P_{n}^{(\ell)}(x;\zeta)\}$ satisfies
\begin{equation}
\label{eq:P01orth-2}
\left(P_{n}^{(\ell)},x^{k}\right)_{\omega^{(\ell)}}=0,\ k\leq n-1,
\end{equation}
which imply
\begin{equation}
\label{eq:P01orth-3}
Q(x;\zeta)=\sum_{i=0}^{n}\left(x^{n},P_{n}^{(\ell)}\right)_{\omega^{(\ell)}}P_{n}^{(\ell)}(x;\zeta),
\end{equation}
for any   polynomial $Q(x;\zeta)$ of degree $\leq n$ in $L^{2}_{\omega^{(\ell)}}[1,+\infty)$.

\begin{figure}
  \centering
  \subfigure[{ $x-P_{n}^{(\ell)}(x,\zeta)$ with $\zeta=1$.}]
  {
  \includegraphics[width=2.6in]{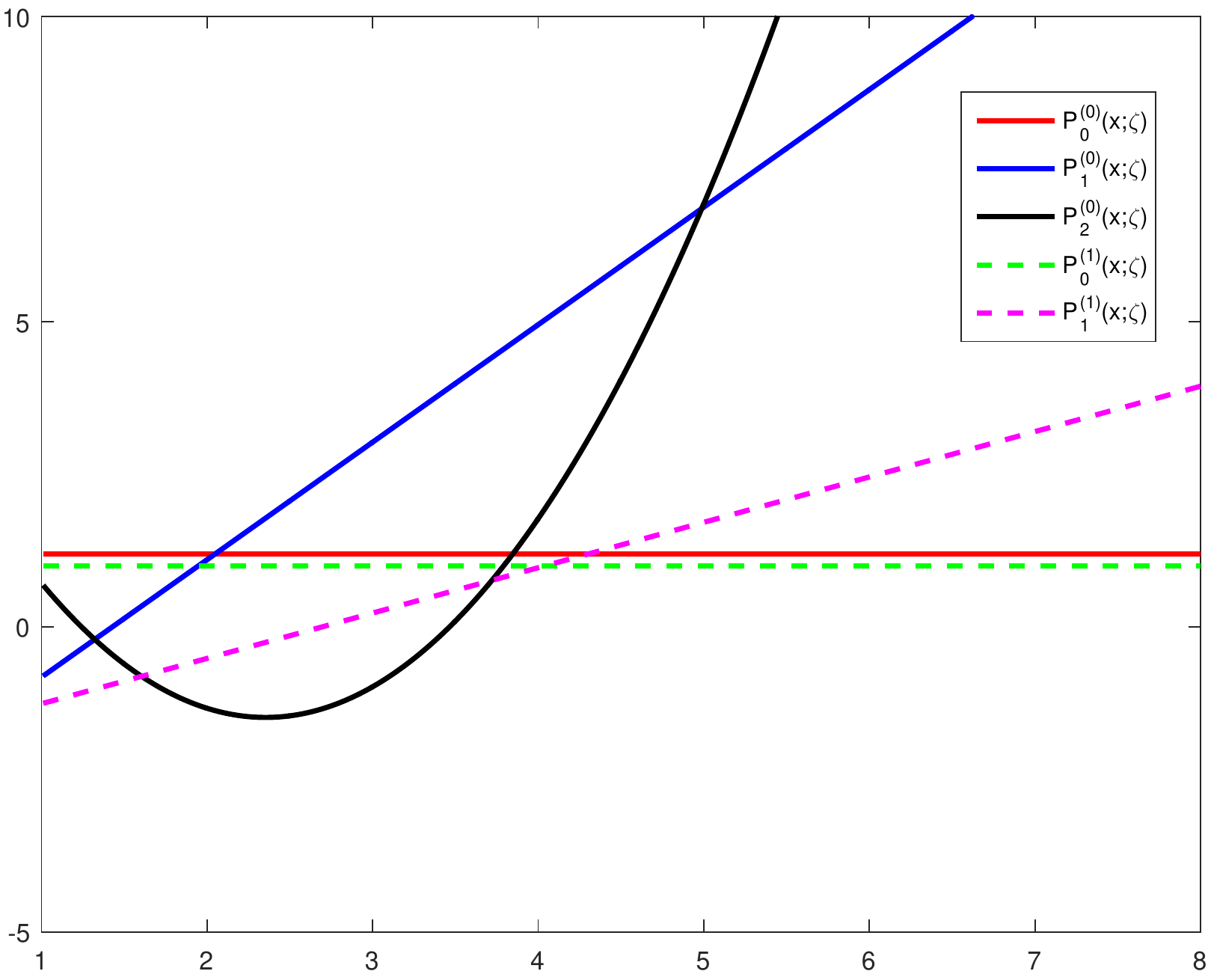}
  }
  \subfigure[{ $\zeta-P_{n}^{(\ell)}(x,\zeta)$ with $x=3$.}]
  {
  \includegraphics[width=2.6in]{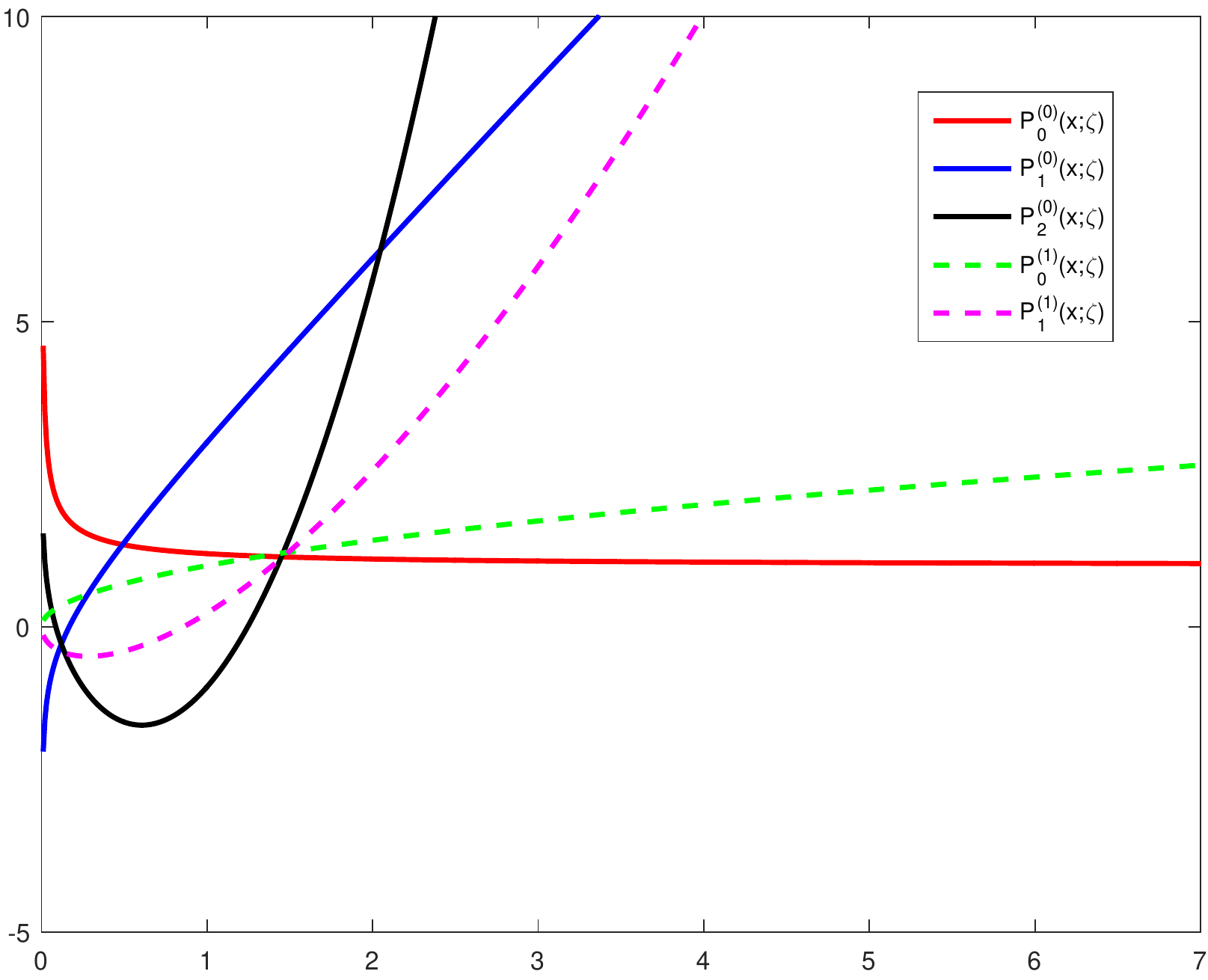}
   }
  \caption{\small Values of the polynomials in \eqref{EQ-3.3aaaaaaa} with respect to $x$ in (a) and $\zeta$ in (b).}
\label{fig:polyp}
\end{figure}

The orthogonal polynomials $\{P_{n}^{(\ell)}(x;\zeta)\}$
{can} be {obtained} by using the Gram-Schmidt  process.
For example, several orthogonal polynomials of lower degree  
are given as follows
\begin{align}\label{EQ-3.3aaaaaaa}
	\begin{aligned}
     P_{0}^{(0)}(x;\zeta)&=\frac{1}{\sqrt{G(\zeta)-2\zeta^{-1}}},\\
     P_{1}^{(0)}(x;\zeta)&=\frac{\sqrt{G(\zeta)-2\zeta^{-1}}}{\sqrt{G(\zeta)^2-3\zeta^{-1}G(\zeta)+2\zeta^{-2}-1}}\left(x-\frac{1}{G(\zeta)-2\zeta^{-1}}\right),\\
     P_{2}^{(0)}(x;\zeta)&=\frac{\zeta\sqrt{G(\zeta)^2-3\zeta^{-1}G(\zeta)+2\zeta^{-2}-1}}{\sqrt{2G(\zeta)^3-7\zeta ^{-1}G(\zeta)^2-2G(\zeta)+6\zeta^{-2}G(\zeta)+\zeta^{-1}}}\\
     &\cdot \left(x^2-\frac{G(\zeta)^2-2\zeta^{-1}G(\zeta)-1}{\zeta\left(G(\zeta)^2-3\zeta^{-1}G(\zeta)+2\zeta^{-2}-1\right)}x
     -\frac{G(\zeta)^2-3\zeta^{-1}G(\zeta)+\zeta^{-2}-1}{G(\zeta)^2-3\zeta^{-1}G(\zeta)+2\zeta^{-2}-1}\right),\\
     P_{0}^{(1)}(x;\zeta)&=\sqrt{\zeta},\\
     P_{1}^{(1)}(x;\zeta)&=\frac{\sqrt{\zeta}}{\sqrt{-G(\zeta)^2+3\zeta^{-1}G(\zeta)+1}}\left(x-G(\zeta)\right),
\end{aligned}\end{align}
{ plotted  in Fig. \ref{fig:polyp} with respect to $x$ and $\zeta$. }

It shows that
the coefficients in those  orthogonal polynomials are so irregular
that it will be very complicate to  study
 the properties of $\{P_{n}^{(\ell)}(x;\zeta)\}$.
Let $c_{n}^{(\ell)}$  be the leading coefficient of
$P_{n}^{(\ell)}(x;\zeta)$, $\ell=0,1$.
Without loss of generality,  assume $c_{n}^{(\ell)}>0$, $\ell=0,1$.
%
Due to  the important result
on the zeros of orthogonal polynomials \cite[Theorem 3.2]{SP:2011},
the polynomial   $P_{n}^{(\ell)}(x;\zeta)$ has   exactly $n$ real simple zeros
in the interval $(1,+\infty)$, $\ell=0,1$.
Thus if   those zeros are denoted by $\{x_{i,n}^{(\ell)}\}_{i=1}^{n}$ in an increasing order,
 then the polynomial $P_{n}^{(\ell)}(x;\zeta)$ {can} be rewritten as follows
\begin{equation}
     \label{eq:Preprezero}
     P_{n}^{(\ell)}(x;\zeta)=c_{n}^{(\ell)}\prod_{i=1}^{n}(x-x_{i,n}^{(\ell)}).
\end{equation}

 In the following, we want to derive the  recurrence relations of $\{P_{n}^{(\ell)}(x;\zeta)\}$,
 calculate their derivatives with respect to $x$ and $\zeta$,  respectively,
 and study the properties of zeros and coefficient matrices in the  recurrence relations.


\subsection{Recurrence relations}
\label{subsec:recurr}
This section presents  the recurrence relations for the orthogonal polynomials
 $\{P_{n}^{(\ell)}(x;\zeta)\}$, $\ell=0,1$,
 the recurrence relations between  $\{P_{n}^{(0)}(x;\zeta)\}$ and  $\{P_{n}^{(1)}(x;\zeta)\}$,
 and the specific forms of the coefficients in those recurrence relations.

Using  the three-term recurrence relation
and the existence theorem of zeros of general orthogonal polynomials in Theorems 3.1 and 3.2 of \cite{SP:2011}
{gives} the following conclusion.

\begin{thm}
\label{thm:rec}
For $\ell=0,1$,
a three-term recurrence relation for
the orthogonal polynomials $\{P_{n}^{(\ell)}(x;\zeta)\}$ {can} be given by
\begin{equation}
  xP_{n}^{(\ell)}=a_{n-1}^{(\ell)}P_{n-1}^{(\ell)}+b_{n}^{(\ell)}P_{n}^{(\ell)}+a_{n}^{(\ell)}P_{n+1}^{(\ell)}, \label{eq:recP0P1}
\end{equation}
or in the matrix-vector form
\begin{equation}
\label{eq:recP0P1mat}
x\vec{P}_{n}^{(\ell)}=\vec{J}_{n}^{(\ell)}\vec{P}_{n}^{(\ell)}+a_{n}^{(\ell)}P_{n+1}^{(\ell)}\vec{e}_{n+1},
\  \vec{P}_{n}^{(\ell)}:=(P_{0}^{(\ell)},\cdots,P_{n}^{(\ell)})^{T},
\end{equation}
where  both coefficients
\begin{equation}
\label{eq:abn}
a_{n}^{(\ell)}:=\left(xP_{n}^{(\ell)},P_{n+1}^{(\ell)}\right)_{\omega^{(\ell)}}=\frac{c_{n}^{(\ell)}}{c_{n+1}^{(\ell)}}, \quad b_{n}^{(\ell)}:=\left(xP_{n}^{(\ell)},P_{n}^{(\ell)}\right)_{\omega^{(\ell)}}=\sum_{i=1}^{n+1}x_{i,n+1}^{(\ell)}-\sum_{i=1}^{n}x_{i,n}^{(\ell)},
\end{equation}
are positive,
 $\vec{e}_{n+1}$ is the last column of the identity matrix of order $(n+1)$,
 and
\[
     \vec{J}_{n}^{(\ell)}:=\begin{pmatrix}
     b_{0}^{(\ell)} &  a_{0}^{(\ell)}  & 0        &            &   & \\
     a_{0}^{(\ell)}& b_{1}^{(\ell)} & a_{1}^{(\ell)}&             &  &   \\
 &          \ddots    &     \ddots     &  \ddots  &      &     \\
                         &    &    &  a_{n-2}^{(\ell)}   & b_{n-1}^{(\ell)} &   a_{n-1}^{(\ell)}   \\
                            &           &    &    0 & a_{n-1}^{(\ell)}   & b_{n}^{(\ell)} \\
      \end{pmatrix}\in {\mathbb R}^{(n+1)\times(n+1)},
\]
which is symmetric positive definite  tridiagonal matrix
with the spectral radius  larger than $1$.
\end{thm}

Besides, the recurrence relations between  $\{P_{n}^{(0)}(x;\zeta)\}$ and  $\{P_{n}^{(1)}(x;\zeta)\}$
{can} also be obtained.

\begin{thm}
\label{thm:huxiang}
{\tt(i)} Two three-term recurrence relations between $\{P_{n}^{(0)}(x;\zeta)\}$ and $\{P_{n}^{(1)}(x;\zeta)\}$ {can} be given
by
\begin{align}
 &(x^2-1)P_{n}^{(1)}=p_{n}P_{n}^{(0)}+q_{n}P_{n+1}^{(0)}+r_{n+1}P_{n+2}^{(0)},\label{eq:recP01}\\
    &P_{n+1}^{(0)}=r_{n}P_{n-1}^{(1)}+q_{n}P_{n}^{(1)}+p_{n+1}P_{n+1}^{(1)},\label{eq:recP10}
\end{align}
or in the matrix-vector form
\begin{align}
\vec{P}_{n+1}^{(0)}&=\vec{J}_{n}^{T}\vec{P}_{n}^{(1)}+p_{n+1}P_{n+1}^{(1)}\vec{e}_{n+2},\label{eq:recPQ}\\
(x^2-1)\vec{P}_{n}^{(1)}&=\vec{J}_{n}\vec{P}_{n+1}^{(0)}+r_{n+1}P_{n+2}^{(0)}\vec{e}_{n+1}\label{eq:recQP},
\end{align}
where
\begin{equation}
\label{eq:pqrn}
p_{n}:=\frac{c_{n}^{(0)}}{c_{n}^{(1)}}, q_{n}:=\frac{c_{n}^{(1)}}{c_{n+1}^{(0)}}\left(\sum_{i=1}^{n+2}x_{i,n+2}^{(0)}-\sum_{i=1}^{n}x_{i,n}^{(1)}\right)
=\frac{c_{n+1}^{(0)}}{c_{n}^{(1)}}
\sum_{i=1}^{n+1}  ( x_{i,n+1}^{(1)}-x_{i,n+1}^{(0)}),
 r_{n}:=\frac{c_{n-1}^{(1)}}{c_{n+1}^{(0)}},
\end{equation}
 and
\[
    \vec{J}_{n}:=\begin{pmatrix}
    p_{0} & q_{0} & r_{1} & 0     & 0 & \cdots & 0\\
    0     & p_{1} & q_{1} & r_{2} & 0 & \cdots & 0\\
    \     & \ddots& \ddots& \ddots& \ & \      & \ \\
    \     & \     &     \ &    \  & 0 & p_{n}  & q_{n}
    \end{pmatrix}
    \in {\mathbb R}^{(n+1)\times (n+2)}.
\]
{\tt (ii)}   Two two-term recurrence relations
 between  $\{P_{n}^{(0)}(x;\zeta)\}$ and  $\{P_{n}^{(1)}(x;\zeta)\}$
 can be derived as follows
\begin{align}
    &(x^2-1)P_{n}^{(1)}=\tilde{p}_{n}(x+\tilde{q}_{n})P_{n+1}^{(0)}+\tilde{r}_{n}P_{n}^{(0)},
    \label{eq:recP01x}\\
    &P_{n+1}^{(0)}=\frac{1}{\tilde{p}_{n}}(x-\tilde{q}_{n})P_{n}^{(1)}-\frac{a_{n-1}^{(1)}}{a_{n}^{(0)}}\tilde{r}_{n}P_{n-1}^{(1)},\label{eq:recP10x}
\end{align}
where
\begin{equation}
\label{eq:pqrnx}
\tilde{p}_{n}:=\frac{c_{n}^{(1)}}{c_{n+1}^{(0)}},\quad
 \tilde{q}_{n}:=\sum_{i=1}^{n+1}x_{i,n+1}^{(0)}-\sum_{i=1}^{n}x_{i,n}^{(1)},
\quad \tilde{r}_{n}:={p_{n}(1-\tilde{p}_{n}^{2})}.
\end{equation}

\end{thm}

\subsection{Partial derivatives}\label{subsec:deriva}
This section
calculates the derivatives of the polynomial $P_{n}^{(\ell)} (x;\zeta)$ with respect to $x$ and $\zeta$,
$\ell=0,1$.
\begin{thm}\label{thm:derivezeta}
	For $\ell=0,1$,
	 the first-order derivative of  the polynomial $P_{n+1}^{(\ell)} (x;\zeta)$
	with respect to the parameter $\zeta$ satisfies
\begin{equation}
 \frac{\partial P_{n+1}^{(\ell)}}{\partial \zeta}=a_{n}^{(\ell)}P_{n}^{(\ell)}-\frac{1}{2}\left(G(\zeta)-\zeta^{-1}-b_{n+1}^{(\ell)}\right)P_{n+1}^{(\ell)}.
 \label{eq:partialPnzeta}
 \end{equation}

\end{thm}

\begin{thm}
\label{thm:derivex}
The first-order derivatives of  the polynomials $\{P_{n}^{(\ell)} (x;\zeta)\}$
with respect to  the variable $x$ satisfy
 \begin{align}
     &\frac{\partial P_{n+1}^{(0)}}{\partial x}=\frac{n+1}{\tilde{p}_{n}}P_{n}^{(1)}+
     \zeta r_{n}P_{n-1}^{(1)}, \label{eq:derivePn0x}\\
     &(x^2-1)\frac{\partial P_{n}^{(1)}}{\partial x}+xP_{n}^{(1)}=(n+1)\tilde{p}_{n}P_{n+1}^{(0)}+
     \zeta p_{n} P_{n}^{(0)}.
     \label{eq:derivePn1x}
\end{align}

\end{thm}

\subsection{Zeros}
\label{subsec:zeros}
Using the separation theorem of zeros of
general orthogonal polynomials 
 \cite{SP:2011} {gives} the following conclusion on
 our  orthogonal polynomials $\{P_{n}^{(\ell)}(x;\zeta)\}$.

\begin{thm}
\label{thm:zerointerself}
For $\ell=0,1$, the zeros $\{x_{i,n}^{(\ell)}\}_{i=1}^{n}$ of $P_{n}^{(\ell)}(x;\zeta)$ and
$\{x_{i,n+1}^{(\ell)}\}_{i=1}^{n+1}$ of $P_{n+1}^{(\ell)}(x;\zeta)$ satisfy
the  separation property
\[
1<x_{1,n+1}^{(\ell)}<x_{1,n}^{(\ell)}<x_{2,n+1}^{(\ell)}<\cdots<x_{n,n}^{(\ell)}<x_{n+1,n+1}^{(\ell)}.
\]
\end{thm}

There is still another important separation property for the
zeros of  the orthogonal polynomials $\{P_{n}^{(\ell)}(x;\zeta), \ell=0,1\}$.

\begin{thm}
\label{thm:zerointer01}
The $n$ zeros $\{x_{i,n}^{(1)}\}_{i=1}^{n}$ of $P_{n}^{(1)}$ and
$n+1$ zeros of $\{x_{i,n+1}^{(0)}\}_{i=1}^{n+1}$ of $P_{n+1}^{(0)}$ satisfy
\[
1<x_{1,n+1}^{(0)}<x_{1,n}^{(1)}<x_{2,n+1}^{(0)}<\cdots<x_{n,n}^{(1)}<x_{n+1,n+1}^{(0)}.
\]

\end{thm}

According to Theorems \ref{thm:zerointerself} and \ref{thm:zerointer01}, we can further
know  the sign of the coefficients of the recurrence relations in Theorem \ref{thm:huxiang}.
\begin{corollary}
\label{col:sign}
All quantities $p_n, q_n, r_n$ in \eqref{eq:pqrn} and
$\tilde{p}_n, \tilde{q}_n, \tilde{r}_n$ in
\eqref{eq:pqrnx} are positive.
\end{corollary}

Using  Corollary \ref{col:sign}, $\tilde{r}_{n}={p_{n}(1-\tilde{p}_{n}^{2})}$, {and $\tilde{p}_{n}=(c_{n+1}^{(0)})^{-1}c_{n}^{(1)}$}
give the following corollary.
\begin{corollary}
\label{corollary3.2}
	The leading coefficient of  $P_{n+1}^{(0)}$ is
	larger than that of $P_{n}^{(1)}$, i.e. $c_{n+1}^{(0)}>c_{n}^{(1)}$.
\end{corollary}

According to Theorems \ref{thm:derivezeta} and \ref{thm:zerointerself},
we have further  the following conclusion.
\begin{corollary}
\label{col:varx}
The zeros $\{x_{i,n}^{(\ell)}\}_{i=1}^{n}$ of $P_{n}^{(\ell)}$ strictly decrease   with respect to
$\zeta$, i.e.
\[
    \frac{\partial x_{i,n}^{(\ell)}}{\partial \zeta}<0.
\]

\end{corollary}

\subsection{Generalized eigenvalues and eigenvectors of coefficient matrices in the recurrence relations}
\label{subsec:eigenorth}
This section discusses 
the generalized eigenvalues and eigenvectors of two $(2n+1)\times (2n+1)$ matrices   $\vec{A}^{0}_{n}$ and $\vec{A}^{1}_{n}$,
defined by
\begin{equation}\label{EQ:3.20-000000}
    \vec{A}^{0}_{n}:=\begin{pmatrix}
    \vec{J}_{n}^{(0)}& \vec{O}\\
    \vec{O}& \vec{J}_{n-1}^{(1)}
    \end{pmatrix}
    ,\quad
    \vec{A}^{1}_{n}:=
      \begin{pmatrix}
    \vec{O}& \vec{J}_{n-1}^{T}\\
    \vec{J}_{n-1}& \vec{O}
    \end{pmatrix},
\end{equation}
where  $\vec{J}_{n}^{(0)}$, $\vec{J}_{n}^{(1)}$, and $\vec{J}_{n}$ appear
in the recurrence relations in Theorems \ref{thm:rec} and \ref{thm:huxiang}.

Consider  the following generalized eigenvalue problem (2nd sense): Find a vector $\vec y$
that obeys $ \vec{A}^{1}_{n} \vec y= \hat{\lambda}   \vec{A}^{0}_{n} \vec y$.
If let $\vec u$ denote the first $n+1$ rows of $\vec y$,
and $\vec v$ be the last $n$ rows of $\vec y$,
then
\begin{equation}
  \hat{\lambda}  \vec{J}_{n}^{(0)}\vec{u}=\vec{J}_{n-1}^{T}\vec{v},\quad
   \hat{\lambda}  \vec{J}_{n-1}^{(1)}\vec{v}=\vec{J}_{n-1}\vec{u}.\label{eq:eigeq1}
\end{equation}

Multiplying \eqref{eq:recP0P1mat}, \eqref{eq:recPQ},  and  \eqref{eq:recQP} by $P_{n}^{(1)}(-x;\zeta)$ with $|x|>1$ gives
\begin{align}
&\vec{P}_{n}^{(0)}(x;\zeta)P_{n}^{(1)}(-x;\zeta)= \frac{1}{x}\vec{J}_{n}^{(0)}\vec{P}_{n}^{(0)}(x;\zeta)P_{n}^{(1)}(-x;\zeta)+\frac{1}{x}a_{n}^{(0)}P_{n+1}^{(0)}(x;\zeta)
P_{n}^{(1)}(-x;\zeta)\vec{e}_{n+1},\label{eq:eqm1}\\
&\vec{P}_{n-1}^{(1)}(x;\zeta)P_{n}^{(1)}(-x;\zeta)= \frac{1}{x}\vec{J}_{n-1}^{(1)}\vec{P}_{n-1}^{(1)}(x;\zeta)P_{n}^{(1)}(-x;\zeta)+\frac{1}{x}a_{n-1}^{(1)}P_{n}^{(1)}(x;\zeta)
P_{n}^{(1)}(-x;\zeta)\vec{e}_{n},\label{eq:eqm2}\\
&\vec{P}_{n}^{(0)}(x;\zeta)P_{n}^{(1)}(-x;\zeta)=\vec{J}_{n-1}^{T}\vec{P}_{n-1}^{(1)}(x;\zeta)P_{n}^{(1)}(-x;\zeta)+p_{n}P_{n}^{(1)}(x;\zeta)P_{n}^{(1)}(-x;\zeta)\vec{e}_{n+1}\label{eq:eqm3},\\
&(x^2-1)\vec{P}_{n-1}^{(1)}(x;\zeta)P_{n}^{(1)}(-x;\zeta)=\vec{J}_{n-1}\vec{P}_{n}^{(0)}(x;\zeta)P_{n}^{(1)}(-x;\zeta)+r_{n}P_{n+1}^{(0)}(x;\zeta)P_{n}^{(1)}(-x;\zeta)\vec{e}_{n}\label{eq:eqm4}.
\end{align}
If substituting \eqref{eq:eqm1} and \eqref{eq:eqm2} into  \eqref{eq:eqm3} and\eqref{eq:eqm4} respectively, then one {obtains}
\begin{align}
\nonumber \vec{J}_{n}^{(0)}&\vec{P}_{n}^{(0)}(x;\zeta)P_{n}^{(1)}(-x;\zeta)=x\vec{J}_{n-1}^{T}\vec{P}_{n-1}^{(1)}(x;\zeta)P_{n}^{(1)}(-x;\zeta)\\
&+\left(xp_{n}P_{n}^{(1)}(x;\zeta)P_{n}^{(1)}(-x;\zeta)-a_{n}^{(0)}P_{n+1}^{(0)}(x;\zeta)
P_{n}^{(1)}(-x;\zeta)\right)\vec{e}_{n+1}\label{eq:eqm11},\\
\nonumber \frac{x^2-1}{x}&\vec{J}_{n-1}^{(1)}\vec{P}_{n-1}^{(1)}(x;\zeta)P_{n}^{(1)}(-x;\zeta)=\vec{J}_{n-1}\vec{P}_{n}^{(0)}(x;\zeta)P_{n}^{(1)}(-x;\zeta)\\
&+\left(r_{n}P_{n+1}^{(0)}(x;\zeta)P_{n}^{(1)}(-x;\zeta)-\frac{x^2-1}{x}a_{n-1}^{(1)}P_{n}^{(1)}(x;\zeta)
P_{n}^{(1)}(-x;\zeta)\right)\vec{e}_{n}\label{eq:eqm12}.
\end{align}
Transforming \eqref{eq:eqm11} and \eqref{eq:eqm12} by  $x$ to $-x$
and then adding them into \eqref{eq:eqm11} and \eqref{eq:eqm12} respectively gives
\begin{align}
& \frac{\sqrt{x^2-1}}{x}  \vec{J}_{n}^{(0)} \vec u(x;\zeta)
=\vec{J}_{n-1}^{T}\vec v(x;\zeta)
 -\frac{\sqrt{x^2-1}}{x}a_{n}^{(0)}
Q_{2n}(x;\zeta) \vec{e}_{n+1},
\label{eq:eqm1final}
\\
& \frac{\sqrt{x^2-1}}{x}   \vec{J}_{n-1}^{(1)}\vec v(x;\zeta)
 =\vec{J}_{n-1}\vec u(x;\zeta)+r_{n}Q_{2n}(x;\zeta)
\vec{e}_{n},
\label{eq:eqm2final}
\end{align}
for $|x|> 1$, where
\begin{align*}
\vec u(x;\zeta)=&\vec{P}_{n}^{(0)}(x;\zeta)P_{n}^{(1)}(-x;\zeta)+\vec{P}_{n}^{(0)}(-x;\zeta)P_{n}^{(1)}(x;\zeta),\\
\vec v(x;\zeta)=&\sqrt{x^2-1}\left(\vec{P}_{n-1}^{(1)}(x;\zeta)P_{n}^{(1)}(-x;\zeta)-\vec{P}_{n-1}^{(1)}(-x;\zeta)P_{n}^{(1)}(x;\zeta)\right),
\end{align*}
and
\begin{equation}
\label{eq:Q2n}
Q_{2n}(x;\zeta):=P_{n+1}^{(0)}(x;\zeta)P_{n}^{(1)}(-x;\zeta)+P_{n+1}^{(0)}(-x;\zeta)P_{n}^{(1)}(x;\zeta).
\end{equation}

It is not difficult to find that
if the second terms at the right-hand sides of \eqref{eq:eqm1final} and \eqref{eq:eqm2final} disappear, then \eqref{eq:eqm1final} and \eqref{eq:eqm2final} reduce to two equations in
\eqref{eq:eigeq1}.
 Thus in order to {obtain} the generalized eigenvalues and eigenvectors of $\vec{A}_{n}^{0}$ and $\vec{A}_{n}^{1}$,
one has to study the zeros of $Q_{2n}(x;\zeta)$.

\begin{lemma}
\label{lem:zerointer}
The function $Q_{2n}(x;\zeta)$ is an even  polynomial of degree $2n$ and has $2n$
real simple zeros $\{z_{i,n}, i=\pm 1,\cdots,\pm n\}$, which satisfy
$z_{-i,n}=-z_{i,n}$ and $z_{i,n}\in (1,+\infty)$
for $i=1,\cdots,n$.

\end{lemma}

\begin{figure}
      \centering
      \includegraphics[width=13.5cm,height=7cm]{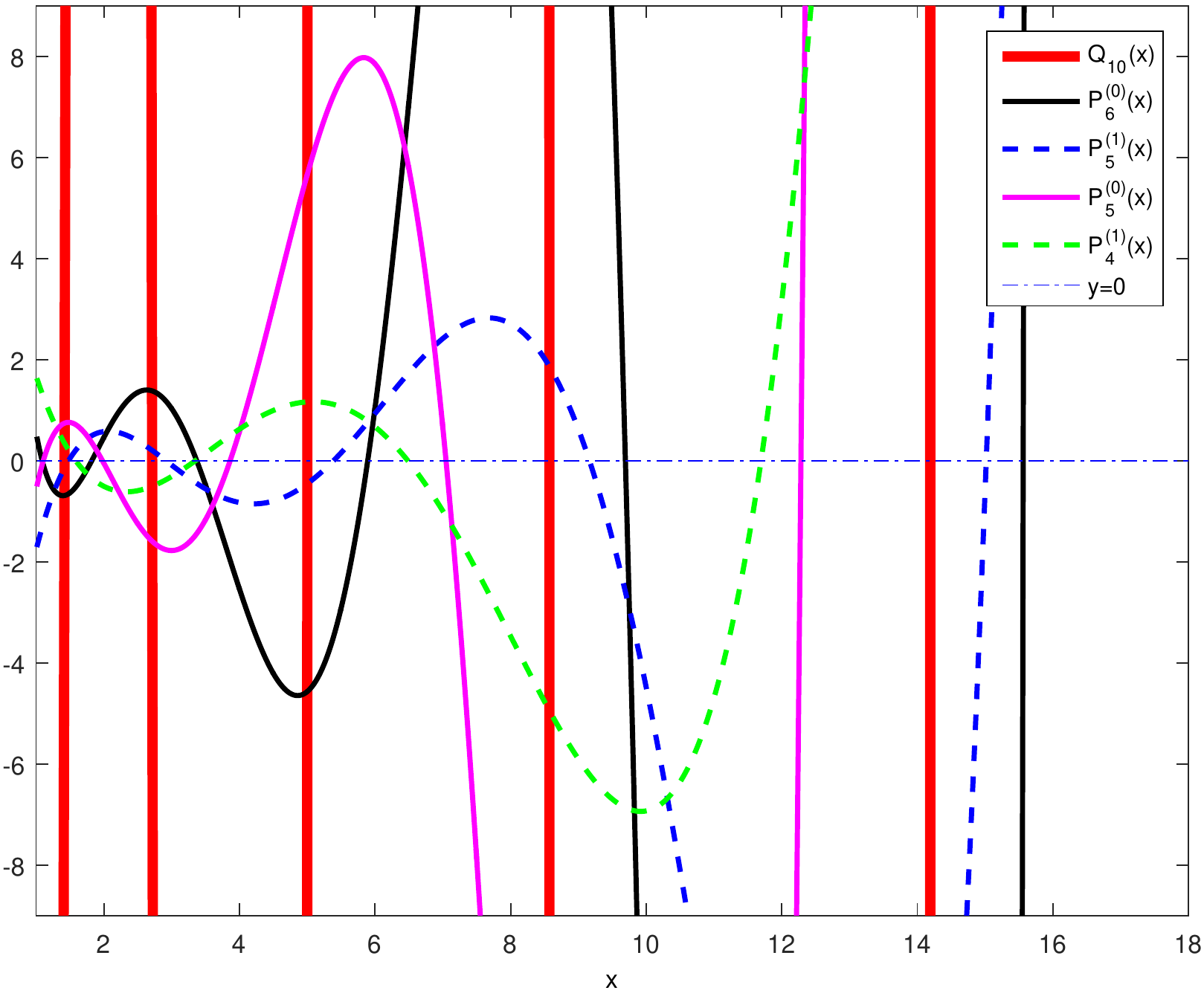}
  \caption{\small Plots of  the polynomials $Q_{10}(x;\zeta)$, $P^{(0)}_{5}(x;\zeta)$,
  $P^{(0)}_{6}(x;\zeta)$, $P^{(1)}_{4}(x;\zeta)$, and $P^{(1)}_{5}(x;\zeta)$ with  $\zeta=1$.}
   \label{fig:poly}
    \end{figure}

The polynomials $Q_{10}(x;\zeta)$, $P^{(0)}_{5}(x;\zeta)$,
  $P^{(0)}_{6}(x;\zeta)$, $P^{(1)}_{4}(x;\zeta)$, and $P^{(1)}_{5}(x;\zeta)$
 with  $\zeta=1$ are plotted in Fig. \ref{fig:poly}, where  the relation between their zeros {can} be clearly observed.

With the aid of Theorems \ref{thm:derivezeta} and \ref{thm:derivex},
we {can} calculate the partial derivatives at $z_{i,n}$ of $Q_{2n}(x;\zeta)$ with respect to
$x$ and $\zeta$.

\begin{lemma}
\label{lem:derivezetax}
At the positive zeros $\{z_{i,n}\}_{i=1}^{n}$, the partial derivatives of $Q_{2n}(x;\zeta)$ satisfy
    \begin{align*}
    \frac{\partial Q_{2n}}{\partial \zeta}(z_{i,n};\zeta)
    =&2\frac{P_{n}^{(1)}(z_{i,n};\zeta)}{P_{n}^{(1)}(-z_{i,n};\zeta)}
    \frac{a_{n}^{(0)}}{\tilde{r}_{n}}\left((\tilde{p}_{n}+
    \tilde{p}_{n}^{-1})z_{i,n}P_{n}^{(1)}(-z_{i,n};\zeta)P_{n+1}^{(0)}(-z_{i,n};\zeta)\right.\\
    &\left.+(z_{i,n}^2-1)P_{n}^{(1)}(-z_{i,n};\zeta)^2
    +P_{n+1}^{(0)}(-z_{i,n};\zeta)^2\right),
    \end{align*}
        \begin{align*}
    \frac{\partial Q_{2n}}{\partial x}(z_{i,n};\zeta)
    =&2\zeta\frac{P_{n}^{(1)}(z_{i,n};\zeta)}{P_{n}^{(1)}(-z_{i,n};\zeta)}
    \frac{a_{n}^{(0)}}{\tilde{r}_{n}}\left((\tilde{p}_{n}+\tilde{p}_{n}^{-1})P_{n+1}^{(0)}(-z_{i,n};\zeta)P_{n}^{(1)}(-z_{i,n};\zeta)\right.\\
    &\left.+z_{i,n}P_{n}^{(1)}(-z_{i,n};\zeta)^{2}
    +z_{i,n}(z_{i,n}^2-1)^{-1}P_{n+1}^{(0)}(-z_{i,n};\zeta)^2 \right).
    \end{align*}
Moreover, one has
\begin{align}\label{EQ3.29zzzz}
    {\rm sign}\left(\frac{\partial Q_{2n}}{\partial \zeta}(z_{i,n};\zeta)\right)={\rm sign}\left(\frac{\partial Q_{2n}}{\partial x}(z_{i,n};\zeta)\right)={\rm sign}\left(\frac{P_{n}^{(1)}(z_{i,n};\zeta)}{P_{n}^{(1)}(-z_{i,n};\zeta)}\right).
\end{align}

\end{lemma}

Similar to Corollary \ref{col:varx}, the following conclusion holds.
\begin{lemma}
\label{lem:varz}
The zeros $\{z_{i,n}, i=\pm 1,\cdots,\pm{n}\}$ of $Q_{2n}(x;\zeta)$ satisfy
\[
 \frac{\partial z_{i,n}}{\partial \zeta}<0,\ i=1,\cdots,n; \quad \frac{\partial z_{i,n}}{\partial \zeta}>0,\ i=-n,\cdots,-1.
\]

\end{lemma}

Thanks to  Lemmas \ref{lem:zerointer} and \ref{lem:varz},
 the generalized eigenvalues and eigenvectors of two \\$(2n+1)\times (2n+1)$
matrices $\vec{A}^{0}_{n}$ and $\vec{A}^{1}_{n}$ {can} be obtained with the aid of the zeros of
$Q_{2n}(x;\zeta)$.

\begin{thm}
\label{thm:eig}
Besides a zero eigenvalue denoted by $\hat{\lambda} _{0,n}$,
the matrix pair $\vec{A}^{0}_{n}$ and $\vec{A}^{1}_{n}$ {has} $2n$ non-zero, real and simple  generalized eigenvalues, which  satisfy
\begin{equation}
\label{eq:eigvalue}
 \hat{\lambda} _{i,n}:=\frac{\sqrt{z_{i,n}^{2}-1}}{z_{i,n}},\quad | \hat{\lambda} _{i,n}|<1,
\ i=\pm 1,\cdots,\pm n,
\end{equation}
and
\begin{equation}
\label{eq:eigparzeta}
\frac{\partial  \hat{\lambda} _{i,n}}{\partial \zeta}<0, \ i=1,\cdots,n;
\quad \frac{\partial \hat{\lambda} _{i,n}}{\partial \zeta}>0,\ i=-n,\cdots,-1.
\end{equation}
Corresponding $(2n+1)$ generalized eigenvectors can be expressed as
\begin{equation}
\vec{y}_{i,n}:=\left(\vec{u}_{i,n}^{T},\vec{v}_{i,n}^{T}\right)^{T},
\label{eq:eigvec}
\end{equation}
with
\begin{align}
\label{EQ:3.29bbbb}
\begin{aligned}
   &\vec{u}_{i,n}=\vec{P}_{n}^{(0)}(z_{i,n};\zeta)P_{n}^{(1)}(-z_{i,n});\zeta
    +\vec{P}_{n}^{(0)}(-z_{i,n};\zeta)P_{n}^{(1)}(z_{i,n};\zeta),\\
    &\vec{v}_{i,n}=\sqrt{z_{i,n}^2-1}\left(\vec{P}_{n-1}^{(1)}(z_{i,n};\zeta)P_{n}^{(1)}(-z_{i,n};\zeta)
    -\vec{P}_{n-1}^{(1)}(-z_{i,n};\zeta)P_{n}^{(1)}(z_{i,n};\zeta)\right).
\end{aligned}\end{align}
for $i=\pm 1,\cdots,\pm n$,
and
\begin{align}
\label{EQ:3.29bbbb--zzz}
    \vec{u}_{0,n}=\vec{P}_{n}^{(0)}(1;\zeta)P_{n+1}^{(0)}(-1;\zeta)-
    \vec{P}_{n}^{(0)}(-1;\zeta)P_{n+1}^{(0)}(1;\zeta),\quad
    \vec{v}_{0,n}={0}.
\end{align}

\end{thm}

\section{Moment method by operator projection} 
\label{sec:moment}
This section begins to extend the moment method by operator projection  \cite{MR:2014}
to the one-dimensional  relativistic Boltzmann equation \eqref{eq:1DBoltz} and derive its arbitrary order  hyperbolic moment model. For the sake of convenience, without loss of generality, units in which both the speed of light $c$  and  rest mass $m$ of particle are equal to one will be used in the following.
{ All proofs are given in the Appendix \ref{sec:App4}.}

\subsection{Weighted polynomial space}
\label{subsec:expand}
In order to use the moment method by the operator projection
to derive the hyperbolic moment model of the kinetic equation,
 we should define weighted polynomial spaces and  norms
 as well as the projection operator.
 Thanks to the equilibrium distribution $f^{(0)} $ in \eqref{eq:equm1},
 the weight function is chosen as $g^{(0)}$, which will be
replaced with  the new notation $g^{(0)}_{[u,\theta]}$, considering
the dependence of $g^{(0)}$  on the  macroscopic fluid velocity $u$
and $\theta=k_BT/m=\zeta^{-1}$, that is
\begin{equation}\label{eq:equm1-zzzzz}
g^{(0)}_{[u,\theta]}=\frac{1}{2K_{1}(\zeta)}\exp\left(-\frac{E}{\theta}\right),\ E=U_{\alpha}p^{\alpha}.
\end{equation}

Associated with the weight function $g^{(0)}_{[u,\theta]}$,  our weighted
polynomial space is defined by
    \[
    \mathbb{H}^{g^{(0)}_{[u,\theta]}}:
    ={\rm span}\left\{p^{\mu_{1}}p^{\mu_{2}}\cdots p^{\mu_{\ell}}g^{(0)}_{[u,\theta]}:
     \ \ \mu_{i}=0,1, \  \ell\in\mathbb{N}\right\},
    \]
which is an infinite-dimensional linear space equipped with the inner product
\[
   <f,g>_{g^{(0)}_{[u,\theta]}}:=\int_{\mathbb{R}}\frac{1}{g^{(0)}_{[u,\theta]}}f(p)g(p)\frac{dp}{p^{0}}, \quad f,g\in \mathbb{H}^{g^{(0)}_{[u,\theta]}}.
\]
Similarly, for a finite positive integer $M\in\mathbb{N}$,   a finite-dimensional  weighted
polynomial space   can  be defined by
\[
    \mathbb{H}^{g^{(0)}_{[u,\theta]}}_{M}:={\rm span}\left\{p^{\mu_{1}}p^{\mu_{2}}\cdots p^{\mu_{\ell}}g^{(0)}_{[u,\theta]}: \ \ \mu_{i}=0,1, \ \ell=0,1,\cdots, M\right\},
    \]
    which is  a closed subspace of $\mathbb{H}^{g^{(0)}_{[u,\theta]}}$ obviously.

Thanks to Theorem \ref{lem:admissible1},
for all physically admissible $u$ and $\theta$ satisfying $|u|<1$ and $\theta>0$,
introduce two  notations
\begin{align}
\label{eq:basis}
\vec{\mathcal{P}}_{\infty}[u,\theta]:=&(\tilde{P}_{0}^{(0)}[u,\theta],\tilde{P}_{1}^{(0)}[u,\theta],\tilde{P}_{0}^{(1)}[u,\theta],
\cdots,\tilde{P}_{M}^{(0)}[u,\theta],\tilde{P}_{M-1}^{(1)}[u,\theta] ,\cdots)^{T},
\\
\label{eq:huap}
\vec{\mathcal{P}}_{M}[u,\theta]:=&(\tilde{P}_{0}^{(0)}[u,\theta],\tilde{P}_{1}^{(0)}[u,\theta],\tilde{P}_{0}^{(1)}[u,\theta],
\cdots,\tilde{P}_{M}^{(0)}[u,\theta],\tilde{P}_{M-1}^{(1)}[u,\theta] )^{T},
\end{align}
 where $\tilde{P}_{k}^{(0)}[u,\theta]=g^{(0)}_{[u,\theta]}P_{k}^{(0)}$ and
$ \tilde{P}_{k}^{(1)}[u,\theta]= g^{(0)}_{[u,\theta]}(U^{0})^{-1}P_{k}^{(1)}p_{<1>}$.


\begin{lemma}
\label{lem:basis}
The set of all components of $\vec{\mathcal{P}}_{\infty}[u,\theta]$ (resp. $\vec{\mathcal{P}}_{M}[u,\theta]$)
form a standard orthogonal basis of $\mathbb{H}^{g^{(0)}_{[u,\theta]}}$
 (resp. $\mathbb{H}^{g^{(0)}_{[u,\theta]}}_{M}$).

\end{lemma}

{
\begin{remark}
In the non-relativistic limit, $E=U_\alpha p^\alpha$, $p_{<1>}$
 and $g^{(0)}_{[u,\theta]}$ reduce to  $p^2$,   $-p$ and
$
\frac{1}{\sqrt{2\pi\theta}}\exp(-\frac{p^2}{\theta})$, respectively,
thus the basis become the generalized Hermite polynomials \cite{1DB:2013}.
\end{remark}
}

Since $\mathbb{H}^{g^{(0)}_{[u,\theta]}}_{M}$ is a subspace of $\mathbb{H}^{g^{(0)}_{[u,\theta]}}_{N}$ when $ M<N<+\infty$,
there exists a matrix \\$P_{M,N}\in \mathbb{R}^{(2M+1)\times(2N+1)}$ with full row rank
such that
$\vec{\mathcal{P}}_{M}[u,\theta]=P_{M,N}\vec{\mathcal{P}}_{N}[u,\theta]$, where
\[
\vec{P}_{M,N}:={\rm diag}\{\vec{I}_{2M+1,2M+1},\vec{O}_{2M+1,2N-2M}\}.
\]
Using the properties of the orthogonal polynomials  $\{P_{n}^{(\ell)}(x;\zeta), \ell=0,1, n\geq 0\}$ in Section  \ref{sec:orth}
{can} further give calculation of the partial derivatives and recurrence relations of the basis functions
$ \{\tilde{P}_{n}^{(0)}[u,\theta], n\geq 0\}$ and $\{\tilde{P}_{n-1}^{(1)}[u,\theta],  n\geq 1\}$.
\begin{lemma}[Derivative relations]
\label{lem:derive}
The partial derivatives   of  basis functions\\
$ \{\tilde{P}_{n}^{(0)}[u,\theta], n\geq 0\}$ and $\{\tilde{P}_{n-1}^{(1)}[u,\theta],  n\geq 1\}$ {can}  be calculated by
\begin{align*}
  \frac{\partial \tilde{P}_{n}^{(0)}[u,\theta]}{\partial s}=&-\frac{\partial \theta}{\partial s}\zeta^2\left(\frac{1}{2}\left(G(\zeta)-\zeta^{-1}-b_{n}^{(0)}\right)\tilde{P}_{n}^{(0)}[u,\theta]-a_{n}^{(0)}\tilde{P}_{n+1}^{(0)}[u,\theta]\right)\\
   &+\frac{\partial u}{\partial s}\frac{1}{(1-u^2)}\left(\left(n\tilde{p}_{n-1}^{-1}-\zeta q_{n-1}\right)\tilde{P}_{n-1}^{(1)}[u,\theta]
   -\zeta p_{n}\tilde{P}_{n}^{(1)}[u,\theta]\right),
   \\
  \frac{\partial \tilde{P}_{n-1}^{(1)}[u,\theta]}{\partial s}=&-\frac{\partial \theta}{\partial s}\zeta^2\left(\frac{1}{2}\left(G(\zeta)-\zeta^{-1}-b_{n-1}^{(1)}\right)\tilde{P}_{n-1}^{(1)}[u,\theta]+a_{n-1}^{(1)}\tilde{P}_{n}^{(1)}[u,\theta]\right)\\
  &+\frac{\partial u}{\partial s}\frac{1}{1-u^2}\left(\left(n\tilde{p}_{n-1}-\zeta q_{n-1}\right)\tilde{P}_{n}^{(0)}[u,\theta]-\zeta r_{n}\tilde{P}_{n+1}^{(0)}[u,\theta]\right),
\end{align*}
for $s=t$ and $x$.
It indicates that $\frac{\partial \tilde{P}_{M}^{(0)}}{\partial s}$ and $\frac{\partial \tilde{P}_{M-1}^{(0)}}{\partial s}$
$\in \mathbb{H}^{g^{(0)}_{[u,\theta]}}_{M+1}$.

\end{lemma}

\begin{lemma}[Recurrence relations]
\label{lem:rec}
The basis functions \\
{$\{\tilde{P}_{n}^{(0)}[u,\theta], n\geq 0\}$ and
$\{\tilde{P}_{n-1}^{(1)}[u,\theta]$, $n\geq 1\}$} satisfy the following recurrence relations
\begin{equation}
    \label{eq:rec}
    \begin{aligned}
      p^{0}\vec{\mathcal{P}}_{M}[u,\theta]=&\vec{M}_{M}^{t}\vec{\mathcal{P}}_{M}[u,\theta]\\
      &+\left(-U^{1}p_{M}\tilde{P}_{M}^{(1)}[u,\theta]+U^{0}a_{M}^{(0)}\tilde{P}_{M+1}^{(0)}[u,\theta]\right)\vec{e}_{2M+1}^{1}\\
      &+\left(-U^{1}r_{M}\tilde{P}_{M+1}^{(0)}[u,\theta]+U^{0}a_{M-1}^{(1)}\tilde{P}_{M}^{(1)}[u,\theta]\right)\vec{e}_{2M+1}^{2},\\
      p\vec{\mathcal{P}}_{M}[u,\theta]=&\vec{M}_{M}^{x}\vec{\mathcal{P}}_{M}[u,\theta]\\
      &+\left(-U^{0}p_{M}\tilde{P}_{M}^{(1)}[u,\theta]+U^{1}a_{M}^{(0)}\tilde{P}_{M+1}^{(0)}[u,\theta]\right)\vec{e}_{2M+1}^{1}\\
      &+\left(-U^{0}r_{M}\tilde{P}_{M+1}^{(0)}[u,\theta]+U^{1}a_{M-1}^{(1)}\tilde{P}_{M}^{(1)}[u,\theta]\right)\vec{e}_{2M+1}^{2},
    \end{aligned}
    \end{equation}
where
$\vec{e}_{2M+1}^{1}$ and $\vec{e}_{2M+1}^{2}$ are the penultimate and the last column of the identity matrix of order $(2M+1)$, respectively, and
     \begin{equation}
    \label{eq:MtMx}
    \begin{aligned}
    \vec{M}_{M}^{t}:&=-U^{1}\vec{P}_{M}^{p}\vec{A}_{M}^{1}(\vec{P}_{M}^{p})^{T}+ U^{0}\vec{P}_{M}^{p}\vec{A}_{M}^{0}(\vec{P}_{M}^{p})^{T} ,\\ \vec{M}_{M}^{x}:&=-U^{0}\vec{P}_{M}^{p}\vec{A}_{M}^{1}(\vec{P}_{M}^{p})^{T}+U^{1}\vec{P}_{M}^{p}\vec{A}_{M}^{0}(\vec{P}_{M}^{p})^{T},
    \end{aligned}
    \end{equation}
in which   $\vec{P}_{M}^{p}$ is a permutation matrix  making
     \begin{equation}
     \label{eq:MtMx-22222}
    \vec{P}_{M}^{p}\vec{\mathcal{\tilde{P}}}_{M}[u,\theta]=\vec{\mathcal{P}}_{M}[u,\theta], \quad 
    \end{equation}
with
\[
    \vec{\mathcal{\tilde{P}}}_{M}[u,\theta]:=(\tilde{P}_{0}^{(0)}[u,\theta],
    \cdots,\tilde{P}_{M}^{(0)}[u,\theta],\tilde{P}_{0}^{(1)}[u,\theta],\cdots,\tilde{P}_{M-1}^{(1)}[u,\theta])^{T}.
\]

\end{lemma}

For a finite integer $M\geq1$,  define an operator
$\Pi_{M}[u,\theta]: \mathbb{H}^{g^{(0)}_{[u,\theta]}}\rightarrow \mathbb{H}^{g^{(0)}_{[u,\theta]}}_{M}$
 by
\begin{align}\label{EQ-projection-aaaa}
\Pi_{M}[u,\theta]f:=\sum_{i=0}^{M}f_{i}^{0}\tilde{P}_{i}^{(0)}[u,\theta]+
\sum_{j=0}^{M-1}f_{j}^{1}\tilde{P}_{j}^{(1)}[u,\theta],
\end{align}
or  in a compact form
\begin{align}\label{EQ-projection-aaaa2}
\Pi_{M}[u,\theta]f=[\vec{\mathcal{P}}_{M}[u,\theta],\vec{f}_{M}]_{M},
\end{align}
where
\begin{align}
f_{i}^{0}&=<f,\tilde{P}_{i}^{(0)}[u,\theta]>_{g^{(0)}_{[u,\theta]}},\  i\leq M,\ \
f_{j}^{1}=<f,\tilde{P}_{j}^{(1)}[u,\theta]>_{g^{(0)}_{[u,\theta]}},\  j\leq M-1
\label{eq:deff1},\\
\vec{f}_{M}&=(f_{0}^{0},f_{1}^{0},f_{0}^{1},\cdots,f_{M}^{0},f_{M-1}^{1})^{T}.
\label{EQ-projection-bbbb}
\end{align}
and the symbol $[\cdot,\cdot]_{M}$ denotes the common inner product of two $(2M+1)$-dimensional vectors.

\begin{lemma}
\label{lem:project}
The operator $\Pi_{M}[u,\theta]$ is
 linear bounded and projection operator in sense that
\begin{description}
	\item[(i)]   $\Pi_{M}[u,\theta]f  \in\mathbb{H}^{g^{(0)}_{[u,\theta]}}_{M}$ for all
	$f\in\mathbb{H}^{g^{(0)}_{[u,\theta]}}$,
   \item[(ii)] $\Pi_{M}[u,\theta]f=f$ for all $f\in\mathbb{H}^{g^{(0)}_{[u,\theta]}}_{M}$.
\end{description}

\end{lemma}

\begin{remark}
	 The so-called  Grad type expansion is to expand the distribution function $f(x,p,t)$
	 in the weighted polynomial space $\mathbb{H}^{g^{(0)}_{[u,\theta]}}$ as follows
	  \[
	  f(x,p,t)=\left[\vec{\mathcal{P}}_{\infty}[u,\theta],\vec{f}_{\infty}\right]_{\infty},
	  \]
	  where the symbol $[\cdot,\cdot]_{\infty}$ denotes the common inner product of two infinite-dimensional vectors,
	  and $\vec{f}_{\infty}=(f_{0}^{0},f_{1}^{0},f_{0}^{1},\cdots,f_{M}^{0},f_{M-1}^{1}, \cdots)^{T}$.
\end{remark}

\subsection{Derivation of the moment model}
\label{subsec:deduction}

Based on the weighted polynomial spaces $\mathbb{H}^{g^{(0)}_{[u,\theta]}}$
and $\mathbb{H}^{g^{(0)}_{[u,\theta]}}_M$ in Section \ref{subsec:expand}
and the projection operator  $\Pi_{M}[u,\theta]$ defined in \eqref{EQ-projection-aaaa},
the moment method by the operator projection \cite{MR:2014} {can} be
implemented  for the  1D special relativistic Boltzmann equation \eqref{eq:1DBoltz}.
In view of the fact that  the  variables $\{\rho, u, \theta, \Pi,  {n}^{1}\}$
are several physical quantities of practical interest
and the first three are required in calculating the equilibrium distribution $f^{(0)}$.

The $(2M+1)$-dimensional vector
 \[\vec{W_{M}}=(\rho,u,\theta,\Pi,\tilde{n}^{1},f_{3}^{0},f_{2}^{1},\cdots,f_{M}^{0},f_{M-1}^{1})^{T},
\]
will be considered as the dependent variable vector,  instead of $\vec{f}_{M}$ defined in \eqref{EQ-projection-bbbb},
where $\tilde{n}^{1}:=n^{1}\sqrt{1-u^2}$. The relations between $\vec{W_{M}}$ and
$\vec{f_{M}}$ is
\begin{equation}
\label{eq:constraint}
\vec{f}_{M}=\vec{D}_{M}^{W}\vec{W}_{M},
\end{equation}
where  the square matrix $\vec{D}_M^{\vec{W}}$ depends on $\theta$ and  is of the following explicit form
\begin{equation*}
\vec{D}_{1}^{W}=
\begin{pmatrix}
(c_{0}^{(0)})^{-1}& 0  & 0 \\
0                  & 0  & 0 \\
0                  & 0  & 0
\end{pmatrix},\quad
\vec{D}_{2}^{W}=
\begin{pmatrix}
(c_{0}^{(0)})^{-1}& 0  & 0 & {c_{0}^{(0)}}                  & 0                            \\
0                  & 0  & 0 & c_{1}^{(0)}x_{1,1}^{(0)}          & 0                            \\
0                  & 0  & 0 &      0                      & -c_{0}^{(1)}                   \\
0                  & 0  & 0 & {-}c_{2}^{(0)}x_{1,2}^{(0)}x_{2,2}^{(0)} & 0                            \\
0                  & 0  & 0 & 0                           & {c_{1}^{(1)}}x_{1,1}^{(1)}
\end{pmatrix},
\end{equation*}
and  $\vec{D}_{M}^{W}={\rm diag}\{\vec{D}_{2}^{W},\vec{I}_{2M-4}\}$ for $M\geq3$, {which is derived from \eqref{eq:condition-222} and \eqref{eq:condition}}.
\begin{figure}[h!]
      \centering
      \includegraphics[width=14cm,height=5.2cm]{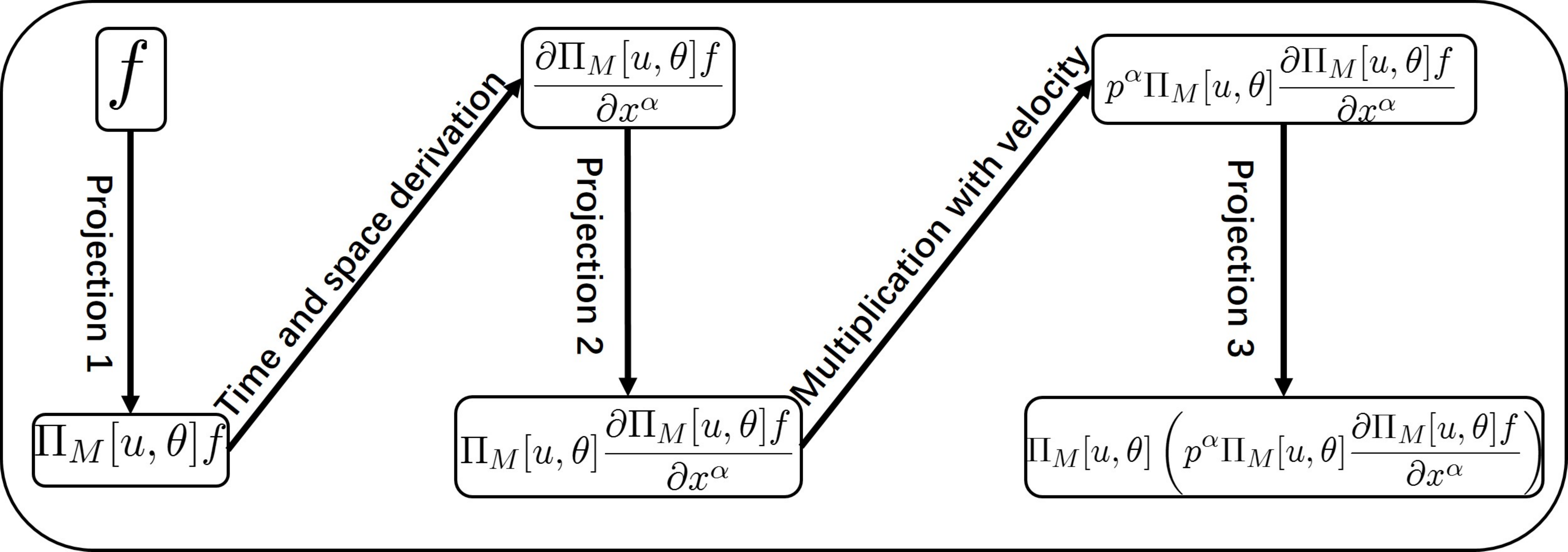}
  \caption{\small Schematic diagram of the moment method by the operator projection
  for the  1D special relativistic Boltzmann equation.}
    \label{fig:step}
\end{figure}
Referring to the schematic diagram shown in  Fig. \ref{fig:step},
the arbitrary order moment system for the  Boltzmann equation \eqref{eq:1DBoltz}
{can} be derived by the operator projection as follows:

\noindent
{\tt Step 1 (Projection 1):}
Projecting the distribution function $f$
  into space $\mathbb{H}^{g^{(0)}_{[u,\theta]}}_{M}$  by the operator $  \Pi_{M}[u,\theta]$ defined in \eqref{EQ-projection-aaaa2}.

\noindent
{\tt Step 2: }
Calculating the partial derivatives in time and space provides
    \begin{align}
    \nonumber \frac{\partial \Pi_{M}[u,\theta]f}{\partial s}&=\left[\frac{\partial \vec{\mathcal{P}}_{M}[u,\theta]}{\partial s},\vec{f}_{M}\right]_{M}
    + \left[\vec{\mathcal{P}}_{M}[u,\theta],\frac{\partial \vec{f}_{M}}{\partial s}\right]_{M}\\
    &=\left[\vec{C}_{M+1}\vec{P}_{M,M+1}^{T}\vec{\mathcal{P}}_{M}[u,\theta],\vec{P}_{M,M+1}^{T}\vec{f}_{M}\right]_{M+1}
    + \left[\vec{\mathcal{P}}_{M}[u,\theta],\frac{\partial \vec{f}_{M}}{\partial s}\right]_{M},\label{eq:diff}
    \end{align}
   for $s=t$ and $x$,
   where $\vec{C}_{M+1}$ is a square matrix of order $(2M+3)$ and directly derived  with the aid of the derivative relations of the
    basis functions in Lemma \ref{lem:derive}.

\noindent
{\tt Step 3 (Projection 2):}
Projecting the partial derivatives in \eqref{eq:diff} into
  the  space $\mathbb{H}^{g^{(0)}_{[u,\theta]}}_{M}$ gives
    \begin{align}
    \nonumber \Pi_{M}[u,\theta]\frac{\partial \Pi_{M}[u,\theta]f}{\partial s}&=
    \left[\vec{\mathcal{P}}_{M}[u,\theta],\vec{C}_{M}^{T}\vec{f}_{M}\right]_{M}
    + \left[\vec{\mathcal{P}}_{M}[u,\theta],\frac{\partial \vec{f}_{M}}{\partial s}\right]_{M}\\
    \nonumber &=\left[\vec{\mathcal{P}}_{M}[u,\theta],\vec{C}_{M}^{T}\vec{D}_{M}^{W}\vec{W}_{M}+\frac{\partial \left(\vec{D}_{M}^{W}\vec{W}_{M}\right)}{\partial s}\right]_{M}\\
    &=:\left[\vec{\mathcal{P}}_{M}[u,\theta],\vec{D}_{M}\frac{\partial \vec{W}_{M}}{\partial s}\right]_{M},
    \label{eq:proj1}
    \end{align}
    where the $(2M+1)$-by-$(2M+1)$  matrix  $\vec{D}_{M}$  {can} be obtained from $\vec{C}_{M}$ and $\vec{D}_{M}^{\vec{W}}$
    and is   of the following form
        \begin{equation}
        \label{eq:Dmatrix}
        \vec{D}_{M}=\begin{pmatrix}
        &\ &\ &D_{2} &\ &\ & O &\ &\ \\
        &0 &* &* &0 &\ &\ &\ &\  \\
        &\vdots&\vdots&\vdots&\ &\  &\vec{I}_{2M-4}&\ &\\\
        &0 &* &* &0 &\ &\  &\ &\
        \end{pmatrix}, \      M\geq3,
        \end{equation}
and
        \begin{align*}
     \vec{D}_{2}&=\begin{pmatrix}
    (c_{0}^{(0)})^{-1}& 0                                     & -\rho\zeta^2(c_{1}^{(0)})^{-2}(c_{0}^{(0)})^{-1}& {c_{0}^{(0)}}                    & 0   \\
    0                  & (1-u^2)^{-1}c_{1}^{(0)}\tilde{n}^{1}  & \rho\zeta^2(c_{1}^{(0)})^{-1}             & c_{1}^{(0)}x_{1,1}^{(0)}           & 0                            \\
    0                  & -(1-u^2)^{-1}c_{0}^{(1)}\rho              & 0                                          &      0                        &-c_{0}^{(1)}                   \\
    0                  & -(1-u^2)^{-1}c_{2}^{(0)}\tilde{n}^{1}(x_{1,2}^{(0)}+x_{2,2}^{(0)})  & 0 & {-}c_{2}^{(0)}x_{1,2}^{(0)}x_{2,2}^{(0)} & 0                            \\
    0                  & -(1-u^2)^{-1}c_{1}^{(1)}\Pi & 0 & 0                           & {c_{1}^{(1)}}x_{1,1}^{(1)}
    \end{pmatrix},\\
     \vec{D}_{1}&=\begin{pmatrix}
    (c_{0}^{(0)})^{-1}& 0                                     & -\rho\zeta^2(c_{1}^{(0)})^{-2}(c_{0}^{(0)})^{-1}
    \\
    0                  & 0 
    & \rho\zeta^2(c_{1}^{(0)})^{-1}
    \\
    0                  & -(1-u^2)^{-1}c_{0}^{(1)}\rho            & 0
    \end{pmatrix},
    \end{align*}
{ where  the
elements ``$*$'' of $\vec D_M$ in \eqref{eq:Dmatrix} are explicitly given by
\begin{align*}
\vec{D}_{M}(2n+1,2)&=\frac{1}{(1-u^2)}\left(\left(n\tilde{p}_{n-1}-\zeta q_{n-1}\right)f_{n-1}^{1}-\zeta r_{n-1}f_{n-2}^{1}\right),\\
\vec{D}_{M}(2n+2,2)&=\frac{1}{(1-u^2)}\left(\left((n+1)\tilde{p}_{n}^{-1}-\zeta q_{n}\right)f_{n+1}^{0}-\zeta p_{n}f_{n}^{0}\right),\\
\vec{D}_{M}(2n+1,3)&=-\zeta^2\left(\frac{1}{2}\left(G(\zeta)-\zeta^{-1}-b_{n}^{(0)}\right)f_{n}^{0}-a_{n-1}^{(0)}f_{n-1}^{0}\right),\\
\vec{D}_{M}(2n+1,3)&=-\zeta^2\left(\frac{1}{2}\left(G(\zeta)-\zeta^{-1}-b_{n}^{(1)}\right)f_{n}^{1}+a_{n-1}^{(1)}f_{n-1}^{1}\right).
\end{align*}
}
 \noindent
 {\tt Step 4:}
 Multiplying  \eqref{eq:proj1} by the particle velocity $(p^\alpha)$ yields
    \begin{align}
    \nonumber p^{0}\Pi_{M}[u,\theta]\frac{\partial \Pi_{M}[u,\theta]f}{\partial {t}}
    &:=[p^{0}\vec{\mathcal{P}}_{M}[u,\theta],\vec{D}_{M}\frac{\partial \vec{W}_{M}}{\partial {t}}]_{M}\\
    &=[\vec{M}_{M+1}^{t}\vec{P}_{M,M+1}^{T}\vec{\mathcal{P}}_{M}[u,\theta],\vec{P}_{M,M+1}^{T}\vec{D}_{M}\frac{\partial \vec{W}_{M}}{\partial {t}}]_{M+1},    \label{eq:MtREC}
    \end{align}
    \begin{align}
    \nonumber p\Pi_{M}[u,\theta]\frac{\partial \Pi_{M}[u,\theta]f}{\partial {x}}
    &:=[p\vec{\mathcal{P}}_{M}[u,\theta],\vec{D}_{M}\frac{\partial \vec{W}_{M}}{\partial {x}}]_{M}\\
    &=[\vec{M}_{M+1}^{x}\vec{P}_{M,M+1}^{T}\vec{\mathcal{P}}_{M}[u,\theta],\vec{P}_{M,M+1}^{T}\vec{D}_{M}\frac{\partial \vec{W}_{M}}{\partial {x}}]_{M+1}.    \label{eq:MxREC}
    \end{align}

\noindent
{\tt Step 5 (Projection 3):}
Projecting \eqref{eq:MtREC} and \eqref{eq:MxREC} into the space $\mathbb{H}^{g^{(0)}_{[u,\theta]}}_{M}$ gives
      \begin{equation}
    \label{eq:pMtREC}
    \Pi_{M}[u,\theta]\left(p^{0}\Pi_{M}[u,\theta]\frac{\partial \Pi_{M}[u,\theta]f}{\partial {t}}\right)=[\vec{\mathcal{P}}_{M}[u,\theta],\vec{M}_{M}^{t}\vec{D}_{M}\frac{\partial \vec{W}_{M}}{\partial {t}}]_{M},
    \end{equation}
    \begin{equation}
    \label{eq:pMxREC}
    \Pi_{M}[u,\theta]\left(p\Pi_{M}[u,\theta]\frac{\partial \Pi_{M}[u,\theta]f}{\partial {x}}\right)=[\vec{\mathcal{P}}_{M}[u,\theta],\vec{M}_{M}^{x}\vec{D}_{M}\frac{\partial \vec{W}_{M}}{\partial {x}}]_{M}.
    \end{equation}

\noindent
{\tt Step 6:}
   Substituting them into  the  1D special relativistic Boltzmann equation \eqref{eq:1DBoltz} derives
   the abstract form of the moment system
 \begin{align}
 \Pi_{M}[u,\theta]\left(p^{\alpha}\Pi_{M}[u,\theta]\left( \frac{\partial \Pi_{M}[u,\theta]f}{\partial x^{\alpha}}\right)\right)
        =\Pi_{M}[u,\theta]Q(\Pi_{M}[u,\theta]f,\Pi_{M}[u,\theta]f),
        \label{eq:moment0000}
        \end{align}
and then       matching the coefficients in front of the basis functions $\{ \tilde{P}_{k}^{(\ell)}[u,\theta]\}$
       leads to   an ``explicit'' matrix-vector form of the moment system
         \begin{align}
          \vec{B}_{M}^{0}\frac{\partial \vec{W}_{M}}{\partial t} +  \vec{B}_{M}^{1}\frac{\partial \vec{W}_{M}}{\partial x}=\vec{S}(\vec{W}_{M}), \label{eq:moment1}
        \end{align}
        which consists of  $(2M+1)$ equations,
        where $ \vec{B}_{M}^{0}=\vec{M}^{t}_{M}\vec{D}_{M}$ and $\vec{B}_{M}^{1}=\vec{M}^{x}_{M}\vec{D}_{M}$.
 For a general collision term $Q(f,f)$,
  it is difficult to {obtain} an explicit expression of
    the source term  $\vec{S}(\vec{W}_{M})$ in \eqref{eq:moment1}.
  For the Anderson-Witting model  \eqref{eq:colAW},
  the right-hand side of \eqref{eq:moment0000} becomes
    \begin{align*}
    \frac{1}{\tau}&\Pi_{M}[u,\theta]Q(\Pi_{M}[u,\theta]f,\Pi_{M}[u,\theta]f)=-\frac{1}{\tau}\Pi_{M}[u,\theta]E\Pi_{M}[u,\theta]{(f-f^{(0)})}\\
    &=-\frac{1}{\tau}\Pi_{M}[u,\theta][\vec{P}_{M+1}^{p}\vec{A}_{M+1}^{0}(\vec{P}_{M+1}^{p})^{T}\vec{P}_{M,M+1}^{T}\vec{\mathcal{P}}_{M}[u,\theta],
     \vec{P}_{M,M+1}^{T}(\vec{f}_{M}-\vec{f}^{(0)}_{M})]_{M+1}\\
    &=-\frac{1}{\tau}[\vec{\mathcal{P}}_{M}[u,\theta],\vec{P}_{M}^{p}\vec{A}_{M}^{0}(\vec{P}_{M}^{p})^{T}{\tilde{\vec{D}}_{M}^{W}\vec{W}_{M}}]_{M},
    \end{align*}
which implies that    the source term  $\vec{S}(\vec{W}_{M})$  {can} be explicitly given by
    \begin{equation}
    \label{eq:colAWfinal}
    \vec{S}(\vec{W}_{M})=-\frac{1}{\tau}\vec{P}_{M}^{p}\vec{A}_{M}^{0}(\vec{P}_{M}^{p})^{T}{\tilde{\vec{D}}_{M}^{W}\vec{W}_{M}}
    =-\frac{1}{\tau}\left(U^{0}\vec{M}_{M}^{t}-U^{1}\vec{M}_{M}^{x}\right)\tilde{\vec{D}}_{M}^{W}\vec{W}_{M},
    \end{equation}
where  $\vec{f}^{(0)}_{M}=\left(\rho\sqrt{G(\zeta)-{2}\zeta^{-1}},0,\cdots,0\right)^{T}$,
and  the matrix $\tilde{\vec{D}}_{M}^{W}$ is the same as $\vec{D}_{M}^{W}$ except that the component of the upper left  corner
 is zero.
 It is worth noting that  the first three components
 of $\vec{S}(\vec{W}_{M})$ are zero due to  \eqref{eq:condition-222} and  \eqref{eq:condition}.

{\begin{remark}
With aid of the explicit forms of $\vec{D}_{1}^{W}$,
   $\vec{D}_{1}$, $\vec{D}_{2}^{W}$, and $\vec{D}_{2}$,
     the explicit form of the moment equations with $M=1$ or 2
     are very easily given.
For example, when $M=1$,
the moment system is written as follows
\[
\vec{B}_{1}^{0}\frac{\partial \vec{W}_{1}}{\partial t} +  \vec{B}_{1}^{1}\frac{\partial \vec{W}_{1}}{\partial x}=0,
\]
where
\[
\vec{B}_{1}^{0}=\begin{pmatrix}
c_{0}^{(0)}U^{0}&c_{0}^{(0)}U^{0}U^{1}\rho & 0\\
(c_{1}^{(0)})^{-1}U^{0}&
\zeta(c_{1}^{(0)})^{-1}(x_{1,2}^{(0)}+x_{2,2}^{(0)})U^{1}(U^{0})^2\rho&
-\zeta^2x_{1,2}^{(0)}x_{2,2}^{(0)}U^{0}(c_{0}^{(0)})^{-2}\rho\\
-\sqrt{\zeta}^{-1}U^{1}&-x_{1,1}^{(1)}\sqrt{\zeta}(U^{0})^3& -\sqrt{\zeta}U^{1}\rho
\end{pmatrix},
\]
\[
\vec{B}_{1}^{1}=\begin{pmatrix}
c_{0}^{(0)}U^{1}&c_{0}^{(0)}(U^{0})^3\rho & 0\\
(c_{1}^{(0)})^{-1}U^{1}&
\zeta(c_{1}^{(0)})^{-1}(x_{1,2}^{(0)}+x_{2,2}^{(0)})(U^{0})^3\rho&
-\zeta^2x_{1,2}^{(0)}x_{2,2}^{(0)}U^{1}(c_{0}^{(0)})^{-2}\rho\\
-\sqrt{\zeta}^{-1}U^{0}&-x_{1,1}^{(1)}\sqrt{\zeta}(U^{0})^2U^{1}& -\sqrt{\zeta}U^{0}\rho
\end{pmatrix}.
\]
It is shown that  those equations become the macroscopic RHD equations
\eqref{eq:conser}
by multiplying those equations by $(B_{1}^{0})^{-1}$.
Thus, the conservation laws are a subset of the equations.
\end{remark}


}
\section{Properties of the moment system}
\label{sec:prop}
This section studies some mathematical and physical properties of  moment system \eqref{eq:moment0000} or \eqref{eq:moment1}.
{ All proofs are given in the Appendix \ref{sec:App5}.}

\subsection{Hyperbolicity, eigenvalues, and eigenvectors}
\label{subsec:Hyper}	
In order to prove the hyperbolicity of the moment system \eqref{eq:moment1},
one has to verify that   $\vec{B}_{M}^{0}$ to be invertible and $\vec{B}_{M}:=(\vec{B}_{M}^{0})^{-1}\vec{B}_{M}^{1}$
to be real diagonalizable. In the following, we always assume that
the first three components of  $\vec{W}_{M}$ satisfy $\rho>0$, $|u|<1$, and $\theta>0$.

\begin{lemma}
 \label{lem:Dinv}
{ If the macroscopic variables satisfy $\rho>0$, $|u|<1$, $\theta>0$ and $\Pi>-\rho\theta$, then} the matrix $\vec{D}_{M}$ is invertible for $M\geq1$.

 \end{lemma}
%
\begin{thm}[Eigenvalues and eigenvectors]
\label{thm:eigmoment}
The $(2M+1)$ eigenvalues of the moment system \eqref{eq:moment1} are
given by
\begin{equation}
\label{eq:lambda}
    \lambda_{i,M}=\frac{u- \hat{\lambda} _{i,M}}{1-u  \hat{\lambda}  _{i,M}}, \ i=-M,\cdots, M,
\end{equation}
satisfying  $|\lambda_{i,M}|<1$,
and  corresponding eigenvectors are
\begin{equation}
\label{eq:r}
 \vec{r}_{i,M}=\vec{D}_{M}^{-1}\vec{P}_{M}^{p}\vec{y}_{i,M},\  i=-M,\cdots,M,
 \end{equation}
where $ \hat{\lambda} _{i,M}$ and $\vec{y}_{i,M}$ are given in Theorem \ref{thm:eig}.

\end{thm}

\begin{lemma}
 \label{lem:Mpositive}
Both real matrices  $U^{0}\vec{M}_{M}^{t}-U^{1}\vec{M}_{M}^{x}$ and $\vec{M}_{M}^{t}$  are positive definite.

 \end{lemma}

\begin{thm}[Hyperbolicity]
\label{thm:hyper}
The moment system  \eqref{eq:moment1} is strictly hyperbolic, and the spectral radius of $\vec{B}_{M}$ is less than one.

\end{thm}


\subsection{Characteristic fields}
\label{subsec:Character}
This section further discusses whether  there exists the
genuinely nonlinear or linearly degenerate
characteristic field of the quasilinear moment system.

\begin{thm}
\label{thm:Deg}
For the moment system  \eqref{eq:moment1},
$\lambda_{0,M}$-characteristic field is
linearly degenerate, i.e.
$$
\nabla_{\vec{W}_M}\lambda_{0,M} (\vec{W}_M)\cdot\vec{r}_{0,M}(\vec{W}_M)=0, \ \forall \vec{W}_M.
$$

\end{thm}

\begin{remark}
\label{remark5.2}
With the aid of numerical experiments, 	we {can} conclude that
for the moment system  \eqref{eq:moment1} with $M\geq4$,
there exist at least two  characteristic fields of
which are neither linearly degenerate  nor genuinely nonlinear, { see
 Appendix \ref{p:remark5.2} for more explanation.}
\end{remark}

\subsection{Linear stability}
\label{subsec:LinearS}
It is obvious that the moment system \eqref{eq:moment1}-\eqref{eq:colAWfinal} has
 the local equilibrium solution  $\vec{W}_{M}^{(0)}=(\rho_{0},u_{0},\theta_{0},0,\cdots,0)^{T}$, where $\rho_{0}$, $u_{0}$,
 and $\theta_{0}$   are   constant and satisfy $\rho_{0}>0$, $|u_{0}|<1$, and $\theta_{0}>0$.
 Similar to the non-relativistic case {\cite{LS:2016}},  { let us linearize the moment system \eqref{eq:moment1}--\eqref{eq:colAWfinal} at $\vec{W}_{M}^{(0)}$. Assuming that $\vec{W}_{M}=\vec{W}_{M}^{(0)}(1+\vec{\bar{W}}_{M})$ and each component  of $\vec{\bar{W}}_{M}$ is small, then the linearized moment system is
\begin{equation}
\label{eq:LME}
\vec{B}_{M}^{0}\big|_{\vec{W}_{M}^{(0)}}\frac{\partial \bar{\vec{W}}_{M}}{\partial t}+\vec{B}_{M}^{1}\big|_{\vec{W}_{M}^{(0)}}\frac{\partial \bar{\vec{W}}_{M}}{\partial x}=\vec{Q}_{M}\big|_{\vec{W}_{M}^{0}}\bar{\vec{W}}_{M},
\end{equation}
where
\[
\vec{Q}_{M}=-\frac{1}{\tau}\left(U^{0}\vec{M}_{M}^{t}-U^{1}\vec{M}_{M}^{x}\right)\tilde{\vec{D}}_{M}^{W}.
\]

Following \cite{LS:2016},   $\bar{\vec{W}}_{M}$ is assumed to be
\[
\bar{\vec{W}}_{M}=\tilde{\vec{W}}_{M}\exp(i(\omega t-kx)),
\]
where $i$ is the imaginary unit, $\tilde{\vec{W}}_{M}$ is the nonzero amplitude, and $\omega$ and $k$ denote
  the frequency and  wave number, respectively.
Substituting the above plane waves into \eqref{eq:LME} gives
\[
\left(i\omega\vec{B}_{M}^{0}-ik\vec{B}_{M}^{1}-\vec{Q}_{M}\right)\big|_{\vec{W}_{M}^{(0)}}\tilde{\vec{W}}_{M}=0.
\]
Because the amplitude $\tilde{\vec{W}}_{M}$ is nonzero,  the above coefficient matrix  is singular, i.e.
\begin{equation}
\label{eq:LS}
\det\left(i\omega\vec{B}_{M}^{0}-ik\vec{B}_{M}^{1}-\vec{Q}_{M}\right)\big|_{\vec{W}_{M}^{(0)}}=0,
\end{equation}
which implies   the dispersion relation between $\omega$ and $k$.

The following linear stability result holds for
 the moment system  \eqref{eq:moment1}--\eqref{eq:colAWfinal}.
\begin{thm}
\label{thm:LS}
The  moment system \eqref{eq:moment1} with the source term \eqref{eq:colAWfinal}
is linearly stable both in space and in time at the local equilibrium, that is,
the linearized moment system \eqref{eq:LME} is stable both in time and in space, that is,
$Im(\omega(k))\geq0$ for each $k\in\mathbb{R}$ and $Re(k(\omega))Im(k(\omega))\leq0$ for each
$\omega\in\mathbb{R}^{+}$, respectively.

\end{thm}}

\subsection{Lorentz covariance}
\label{subsec:Lorentz}
In physics, the  Lorentz covariance is a key property of space-time following from the special theory of relativity, see e.g. \cite{SR:1961}.
This section studies the Lorentz covariance of the  moment system \eqref{eq:moment1}.
{
 Besides the truncations or projection
of distribution function, there are the truncations or projections
of equation in the current moment method.
It is   nontrivial to know  which parts of the expansion of the equation
we have removed in the truncation or projection procedure, and whether they are
Lorentz invariant or not.
}



Some   Lorentz covariant quantities are first pointed out below.

\begin{lemma}
	\label{lem:lorentz}
(i)	Each component of $\vec{D}_{M}^{u} d\vec{W}_{M}$ is  Lorentz invariant,
	where \\$\vec{D}_{M}^{u}:={\rm diag}\{1,(1-u^2)^{-1},1,\cdots,{1}\}$ and
	$d\vec{W}_{M}$ denotes the total differential of $\vec{W}_{M}$.

(ii)	The matrices $\vec{A}_{M}^{0}$, $\vec{A}_{M}^{1}$ and the source term $\vec{S}(\vec{W}_{M})$  defined in \eqref{eq:colAWfinal}  are Lorentz invariant.

\end{lemma}
%

\begin{thm}[Lorentz covariance]
	\label{thm:lorentz}
	The moment system  \eqref{eq:moment1} with the source term \eqref{eq:colAWfinal}
	 is Lorentz covariant.

\end{thm}



\section{Numerical experiment}
\label{sec:NE}
This section conducts a numerical experiment to check the behavior
of  our hyperbolic moment equations (HME) \eqref{eq:moment0000} or \eqref{eq:moment1} with \eqref{eq:colAWfinal}
by solving the Cauchy problem with initial data
\begin{equation}
    \label{eq:riemann}
\vec{W}_{M}(x,0)=
                   \begin{cases}
                    \vec{W}_{M}^{L}, & x<0, \\
                     \vec{W}_{M}^{R}, & x>0,
                   \end{cases}
\end{equation}
where $\vec{W}_{M}^{L}=(7,0,1,0,\cdots,0)^{T}$ and $\vec{W}_{M}^{R}=(1,0,1,0,\cdots,0)^{T}$.
It is similar to the problem for the moment system of the non-relativistic BGK equation used in \cite{1DB:2013}.

\subsection{Numerical scheme}

 The spatial grid $\{x_i, i\in \mathbb Z\}$ considered here is uniform
 so that the stepsize $\Delta x=x_{i+1}-x_{i}$ is constant.
 Thanks to  Theorem \ref{thm:eigmoment},
 the grid in $t$-direction $\{t_{n+1}=t_n+\Delta t, n\in \mathbb N\}$
 {can} be  given with the stepsize  $\Delta t=C_{\mbox{\tiny CFL}} \Delta x$,
where $C_{\mbox{\tiny CFL}} $ denotes the
CFL (Courant-Friedrichs-Lewy) number.  Use $f_i^n$ and $\rho_i^n$  to denote the approximations of $f(x_i, p, t_n)$ and $\rho(x_i,t_n)$ respectively.
 For the purpose of  checking the behavior of our  hyperbolic moment system,
 similar to  \cite{NM:2010},
 we only consider a first-order accurate semi-implicit operator-splitting type  numerical scheme
 for the system \eqref{eq:moment0000} or \eqref{eq:moment1},
 which is formed  into the   convection  and collision steps:     
\begin{equation}
	\label{eq:advfinal}
	\Pi_{M}[u_{i}^{n},\theta_{i}^{n}]\left(p^{0}\Pi_{M}[u_{i}^{n},\theta_{i}^{n}](\Pi f)_{i}^{\ast}\right)
	=\Pi_{M}[u_{i}^{n},\theta_{i}^{n}]\left(p^{0}(\Pi f)_{i}^{n}\right)-\frac{\Delta t}{\Delta x}\left[
	(\Pi F^{-})_{i+\frac{1}{2}}^{n}-(\Pi F^{+})_{i-\frac{1}{2}}^{n}\right],
\end{equation}
and
\begin{align}
	\nonumber &\Pi_{M}[u_{i}^{\ast},\theta_{i}^{\ast}]\left(p^{0}\Pi_{M}[u_{i}^{\ast},\theta_{i}^{\ast}]\frac{(\Pi f)_{i}^{n+1}-(\Pi f)_{i}^{*}}{\Delta t}\right)\\
	=&-\frac{1}{\tau_{i}^{\ast}}\Pi_{M}[u_{i}^{\ast},\theta_{i}^{\ast}](U_{i}^{0\ast}p^{0}-U_{i}^{1\ast}p^{1})\left(I-\Pi_{f\rightarrow f^{(0)}}[u_{i}^{\ast},\theta_{i}^{\ast}]\right)\left(\Pi_{M}[u_{i}^{\ast},\theta_{i}^{\ast}](\Pi f)_{i}^{n+1}\right),
	\label{eq:semicol}
\end{align}
where $(\Pi f)_{i}^{n}:=\Pi_{M}[u_{i}^{n},\theta_{i}^{n}]f_{i}^{n}$  and
the ``numerical fluxes'' $(\Pi F^{-})_{i+\frac{1}{2}}^{n}$ and $(\Pi F^{+})_{i-\frac{1}{2}}^{n}$ are
 derived based on the  nonconservative version of the HLL (Harten-Lax-van Leer) scheme \cite{HLL:2008} and given by
\begin{align*}
	&(\Pi F^{-})_{i+\frac{1}{2}}^{n}=\\
	&\left\{
	\begin{array}{ll}
		\Pi_{M}[u_{i}^{n},\theta_{i}^{n}](p\left(\Pi f)_{i}^{n}\right), & \hbox{$0\leq\lambda_{i+\frac{1}{2}}^{L}$,} \\
		\begin{aligned}
			&\frac{\lambda_{i+\frac{1}{2}}^{R}\Pi_{M}[u_{i}^{n},\theta_{i}^{n}]\left(p(\Pi f)_{i}^{n}\right)-
				\lambda_{i+\frac{1}{2}}^{L}\Pi_{M}[u_{i}^{n},\theta_{i}^{n}]\left(p\Pi_{M}[u_{i}^{n},\theta_{i}^{n}](\Pi f)_{i+1}^{n}\right)}
			{\lambda_{i+\frac{1}{2}}^{R}-\lambda_{i+\frac{1}{2}}^{L}}\\
			&+
			 \frac{\lambda_{i+\frac{1}{2}}^{L}\lambda_{i+\frac{1}{2}}^{R}\left(\Pi_{M}[u_{i}^{n},\theta_{i}^{n}]\left(p^{0}\Pi_{M}[u_{i}^{n},\theta_{i}^{n}](\Pi f)_{i+1}^{n}\right)
				-\Pi_{M}[u_{i}^{n},\theta_{i}^{n}]\left(p^{0}(\Pi f)_{i}^{n}\right)\right)}{\lambda_{i+\frac{1}{2}}^{R}-\lambda_{i+\frac{1}{2}}^{L}},
		\end{aligned} & \hbox{$\lambda_{i+\frac{1}{2}}^{L}<0<\lambda_{i+\frac{1}{2}}^{R}$,} \\
		\Pi_{M}[u_{i}^{n},\theta_{i}^{n}]\left(p\Pi_{M}[u_{i}^{n},\theta_{i}^{n}](\Pi f)_{i+1}^{n}\right), & \hbox{$0\geq\lambda_{i+\frac{1}{2}}^{R}$,}
	\end{array}
	\right.
\end{align*}
and
\begin{align*}
	&(\Pi F^{+})_{i-\frac{1}{2}}^{n}=\\
	&\left\{
	\begin{array}{ll}
		\Pi_{M}[u_{i}^{n},\theta_{i}^{n}]\left(p\Pi_{M}[u_{i}^{n},\theta_{i}^{n}](\Pi f)_{i-1}^{n}\right), & \hbox{$0\leq\lambda_{i-\frac{1}{2}}^{L}$,}\\
		\begin{aligned}
			&\frac{\lambda_{i-\frac{1}{2}}^{R}\Pi_{M}[u_{i}^{n},\theta_{i}^{n}]\left(p\Pi_{M}[u_{i}^{n},\theta_{i}^{n}](\Pi f)_{i-1}^{n}\right)-
				\lambda_{i-\frac{1}{2}}^{L}\Pi_{M}[u_{i}^{n},\theta_{i}^{n}]\left(p(\Pi f)_{i}^{n}\right)}
			{\lambda_{i-\frac{1}{2}}^{R}-\lambda_{i-\frac{1}{2}}^{L}}\\
			&+
			\frac{\lambda_{i-\frac{1}{2}}^{L}\lambda_{i-\frac{1}{2}}^{R}\left(\Pi_{M}[u_{i}^{n},\theta_{i}^{n}]\left(p^{0}(\Pi f)_{i}^{n}\right)\right)
				-\Pi_{M}[u_{i}^{n},\theta_{i}^{n}]\left(p^{0}\Pi_{M}[u_{i}^{n},\theta_{i}^{n}](\Pi f)_{i-1}^{n}\right)}{\lambda_{i-\frac{1}{2}}^{R}-\lambda_{i-\frac{1}{2}}^{L}},
		\end{aligned} & \hbox{$\lambda_{i-\frac{1}{2}}^{L}<0<\lambda_{i-\frac{1}{2}}^{R}$,} \\
		\Pi_{M}[u_{i}^{n},\theta_{i}^{n}](p\left(\Pi f)_{i}^{n}\right), & \hbox{$0\geq\lambda_{i-\frac{1}{2}}^{R}$.}
	\end{array}
	\right.
\end{align*}
Here
$\lambda_{i\pm \frac{1}{2}}^{L}=\min\{\lambda_{i}^{\min},\lambda_{i\pm 1}^{\min}\}$ and $\lambda_{i\pm \frac{1}{2}}^{R}=\max\{\lambda_{i}^{\max},\lambda_{i\pm 1}^{\max}\}$,
where $\lambda_{i}^{\min}$ and $\lambda_{i}^{\max}$ denote the minimum and maximum eigenvalues of the moment system \eqref{eq:moment1}
at the grid point $x_{i}$ respectively, see Theorem \ref{thm:eigmoment}.
In Eq. \eqref{eq:semicol},
 the subscript  $f\rightarrow f^{(0)}$ denotes the transformation from $f$ to $f^{(0)}$ defined by
$
\Pi_{f\rightarrow f^{(0)}}[u_{i}^{\ast},\theta_{i}^{\ast}](\Pi f)_{i}^{*}=f_{i}^{(0)\ast}
$ or  $ \vec{f}_{i,M}^{(0)\ast}=\vec{D}_{M}^{f_{i}^{(0)}}\vec{f}_{i,M}^{\ast}$,
where
\begin{align}\label{EQ.6.3b}
\vec{D}_{M}^{f_{i}^{(0)}}=\left(c_{0,i}^{(0)\ast}\right)^{-2}{\rm diag}\{1,0,\cdots,0\}\left(U_{i}^{0\ast}\vec{M}_{M}^{t*}-U_{i}^{1\ast}\vec{M}_{M}^{x*}\right),
\end{align}
whose nonzero components are  only in the first row and the component in the upper left corner is one.

 The above scheme  \eqref{eq:advfinal}  and  \eqref{eq:semicol} is  implemented as follows:
\begin{itemize}
	\item[(i).] Perform the convection step  \eqref{eq:advfinal} to obtain
	$\Pi_{M}[u_{i}^{n},\theta_{i}^{n}]\left(p^{0}\Pi_{M}[u_{i}^{n},\theta_{i}^{n}](\Pi f)_{i}^{\ast}\right)$,
	and then {obtain} $\Pi_{M}[u^{n}_i,\theta^{n}_i](\Pi f)_{i}^{*}$.
	\item[(ii).] Calculate $u^{*}_i$ and $\theta^{*}_i$ by  solving \eqref{eq:admu} and \eqref{eq:admT},
	and then give $(\Pi f)_{i}^{*}$.
	\item[(iii).] Perform the collision step \eqref{eq:semicol} to obtain
		$\Pi_{M}[u_{i}^{*},\theta_{i}^{*}]\left(p^{0}\Pi_{M}[u_{i}^{*},\theta_{i}^{*}](\Pi f)_{i}^{n+1}\right)$,
and then have	 $\Pi_{M}[u^{*}_i,\theta^{*}_i](\Pi f)_{i}^{n+1}$.
	\item[(iv).]  Calculate $u^{n+1}_i$ and $\theta^{n+1}_i$ by  solving \eqref{eq:admu} and \eqref{eq:admT}, and then {obtain}
	$(\Pi f)_{i}^{n+1}$. Set $n=n+1$ and turn to Step (i).
\end{itemize}
It is worth noting that  when $\Pi_{M}[u_{i}^{n},\theta_{i}^{n}]\left(p^{0}\Pi_{M}[u_{i}^{n},\theta_{i}^{n}](\Pi f)_{i}^{\ast}\right)$ is known,
it is easy to {obtain} $\Pi_{M}[u_{i}^{n},\theta_{i}^{n}]  (\Pi f)_{i}^{\ast}$  in Step (i),
but  it  is more technical to calculate  $ (\Pi f)_{i}^{\ast}$ from the known value of  $\Pi_{M}[u_{i}^{n},\theta_{i}^{n}]  (\Pi f)_{i}^{\ast}$ in Step (ii),
see the following discussion (Lemma \ref{lem:ut}).
The other steps are similar to them.

\begin{lemma}
	\label{lem:projectn}
If $u\in(-1,1)$, $\theta\in\mathbb{R}^{+}$, $M\geq 1$, and $0\leq f(x,p,t)\leq+\infty$,  then
for any polynomial $\tilde{f}$ satisfying $\tilde{f} g^{(0)}_{[u,\theta]}\in\mathbb{H}_{M}^{g^{(0)}_{[u,\theta]}}$,
equivalently $\tilde{f}f\in\mathbb{H}_{M}^{f}$, one has
	\begin{equation}
	\label{eq:projectthm}
	<\tilde{f}f,f>_{f}=<\tilde{f}f,\Pi_{M}[u,\theta]f>_{f}=<\tilde{f}g^{(0)}_{[u,\theta]},\Pi_{M}[u,\theta]f>_{g^{(0)}_{[u,\theta]}}.
	\end{equation}

\end{lemma}

\begin{lemma}
	\label{lem:ut}
If $u_{1},u_{2}\in(-1,1)$, $\theta_{1},\theta_{2}\in\mathbb{R}^{+}$, $M\geq 1$, and $0\leq f(x,p,t)\leq+\infty$, then the identity
	\[
\Pi_{M}[u_{1},\theta_{1}]f=	\Pi_{M}[u_{1},\theta_{1}]\Pi_{M}[u_{2},\theta_{2}]f,
	\]
	holds.

\end{lemma}

 Lemma \ref{lem:ut} implies that in order to calculate
\begin{equation}
\label{eq:utstar}
(\Pi f)_{i}^{\ast}=\Pi_{M}[u_{i}^{\ast},\theta_{i}^{\ast}] (\Pi f)_{i}^{\ast}
=\Pi_{M}[u_{i}^{\ast},\theta_{i}^{\ast}]\Pi_{M}[u_{i}^{n},\theta_{i}^{n}](\Pi f)_{i}^{\ast},
\end{equation}
only $u_{i}^{\ast}$ and $\theta_{i}^{\ast}$  have to be {obtained}. 
It can be done the following procedure.
For the given ``distribution function'' $\Pi_{M}[u_{i}^{n},\theta_{i}^{n}](\Pi f)_{i}^{\ast}$,
calculate corresponding partial {particle flow} $N^\alpha$ and
partial energy-momentum tensor $T^{\alpha\beta}$,
and then solve directly \eqref{eq:admu} to {obtain}  $u_{i}^{\ast}$
and solve iteratively \eqref{eq:admT} to obtain $\theta_{i}^{\ast}$  by using Newton-Raphson method.
\begin{remark}
	The function $G(\theta^{-1})-\theta$ in \eqref{eq:admT} is a strictly monotonic and convex function  of $\theta$ in the
	interval $(0,+\infty)$,  because
	\[
	\frac{\partial^2 \left(G(\theta^{-1})-\theta\right)}{\partial \theta^2}
	= \zeta^2\big(2G(\zeta)^3\zeta^2 - 7G(\zeta)^2\zeta - 2G(\zeta)\zeta^2 + 6G(\zeta) + \zeta\big)
	=\zeta^6 (c_{2}^{(0)}c_{1}^{(0)}c_{0}^{(0)})^{-2}>0,
	\]
	where $c_{i}^{(0)}$ is the leading coefficient of the polynomial
	$P_i ^{(0)}(x;\zeta)$ defined in  \eqref{EQ-3.3aaaaaaa}, $i=0,1,2$.
	It means that
	the Newton-Raphson method for solving  \eqref{eq:admT}
	{converges} with any positive initial guess.
\end{remark}


Before ending this subsection, we discuss the stability of the collision step \eqref{eq:semicol} even though $\tau$ is very small.

\begin{thm}
	\label{thm:colS}
	Semi-implicit scheme \eqref{eq:semicol} is unconditionally stable.

\end{thm}

{ All proofs have been given in the Appendix \ref{sec:App6}.}

\subsection{Numerical results}

In our numerical experiment, the Knudsen number $Kn$ is chosen as $0.05$ and $0.5$, respectively,
the spatial domain  $[-1.5,1.5]$ is divided into a uniform grid of 1000 grid points,
and  $C_{\mbox{\tiny CFL}}=0.9$.
In order to verify  our results, the reference
solutions are provided by using the discrete velocity model (DVM) \cite{DVM:2000}
with a fine spatial grid of $10000$ grid points and 50 Gaussian points in the velocity space.

\begin{figure}
\label{fig:005}
  \centering
  \subfigure[M=1 RHD]
  {
  \includegraphics[width=1.7in]{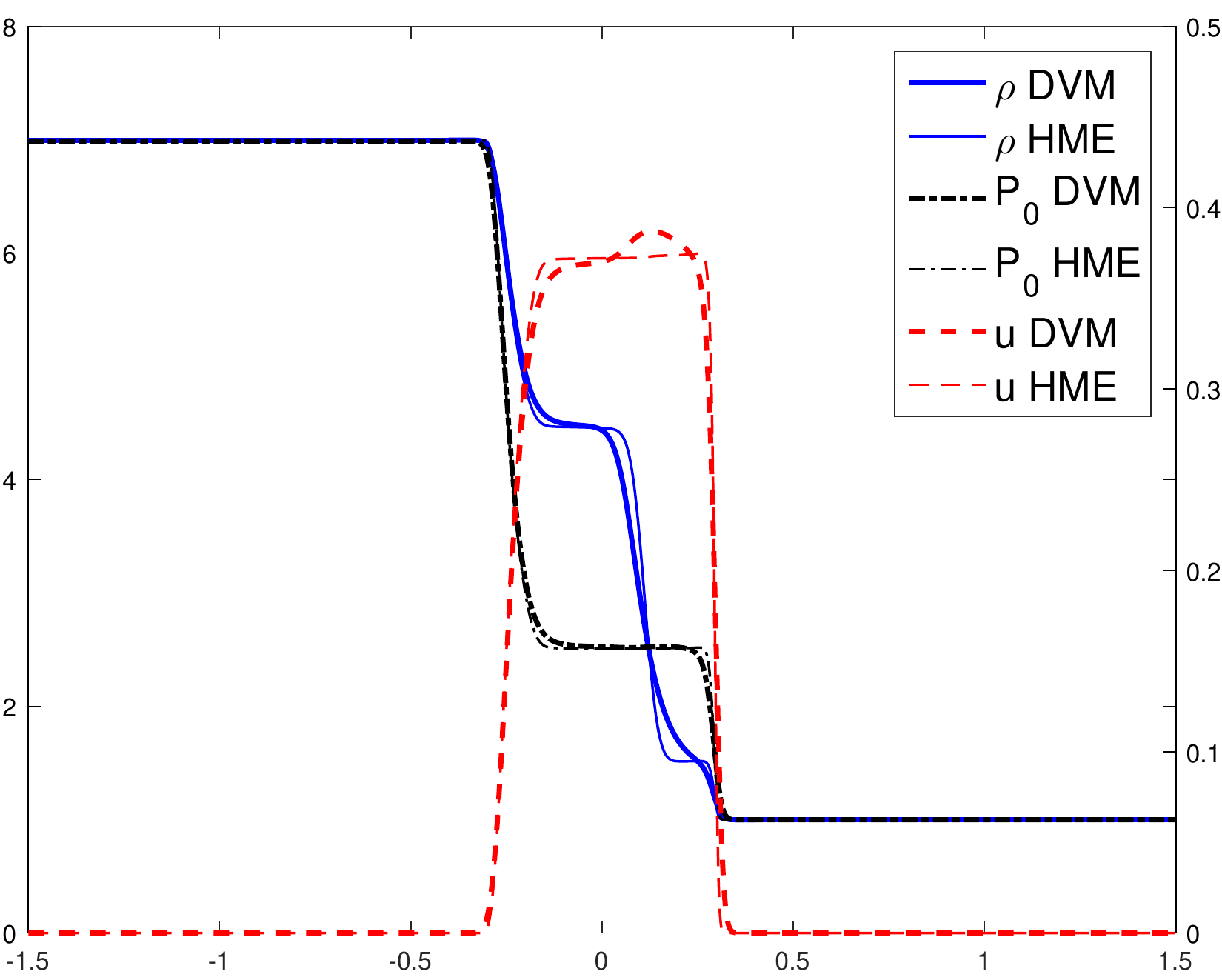}
  }
  \subfigure[M=2]
  {
  \includegraphics[width=1.7in]{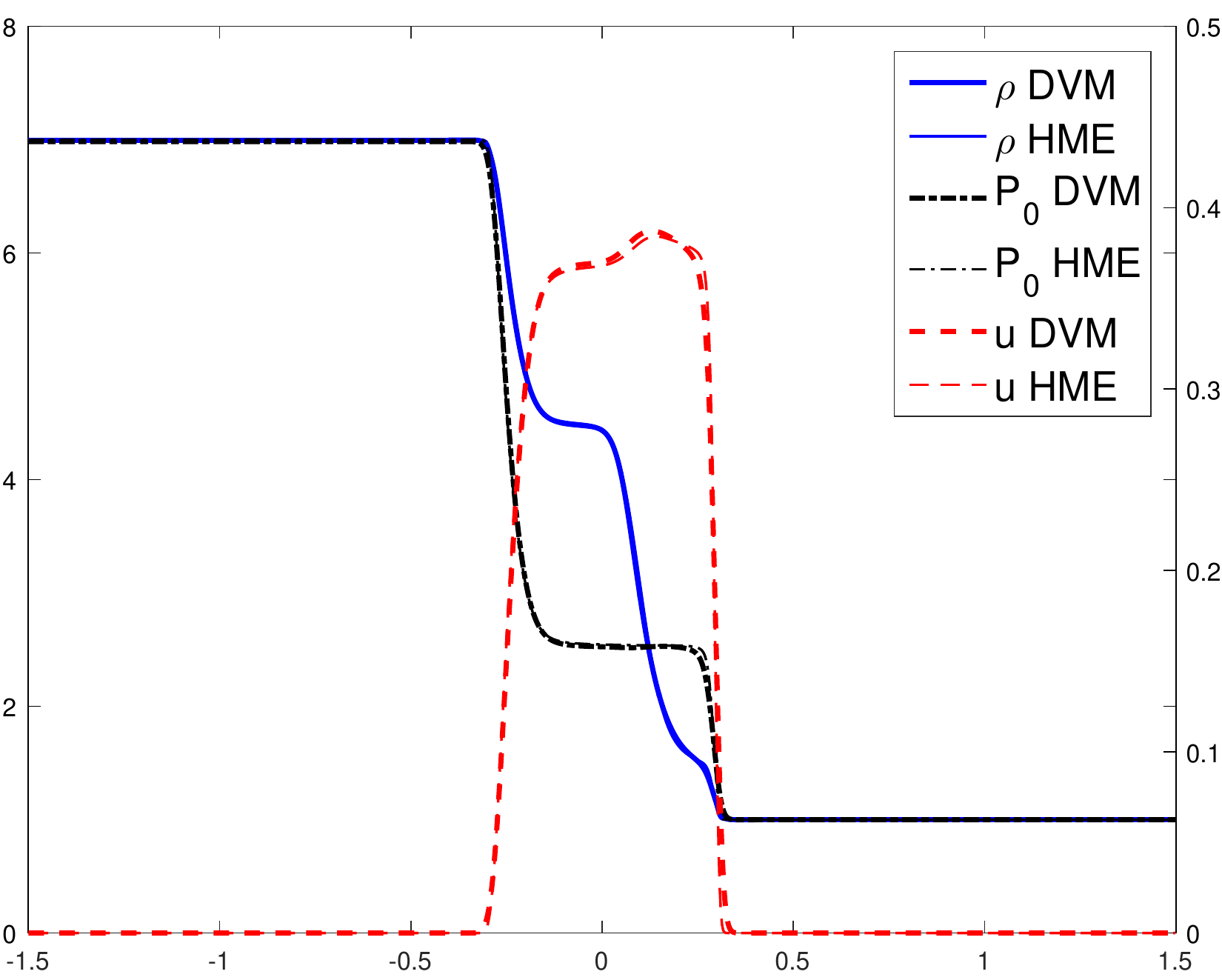}
   }
  \subfigure[M=3]
  {
  \includegraphics[width=1.7in]{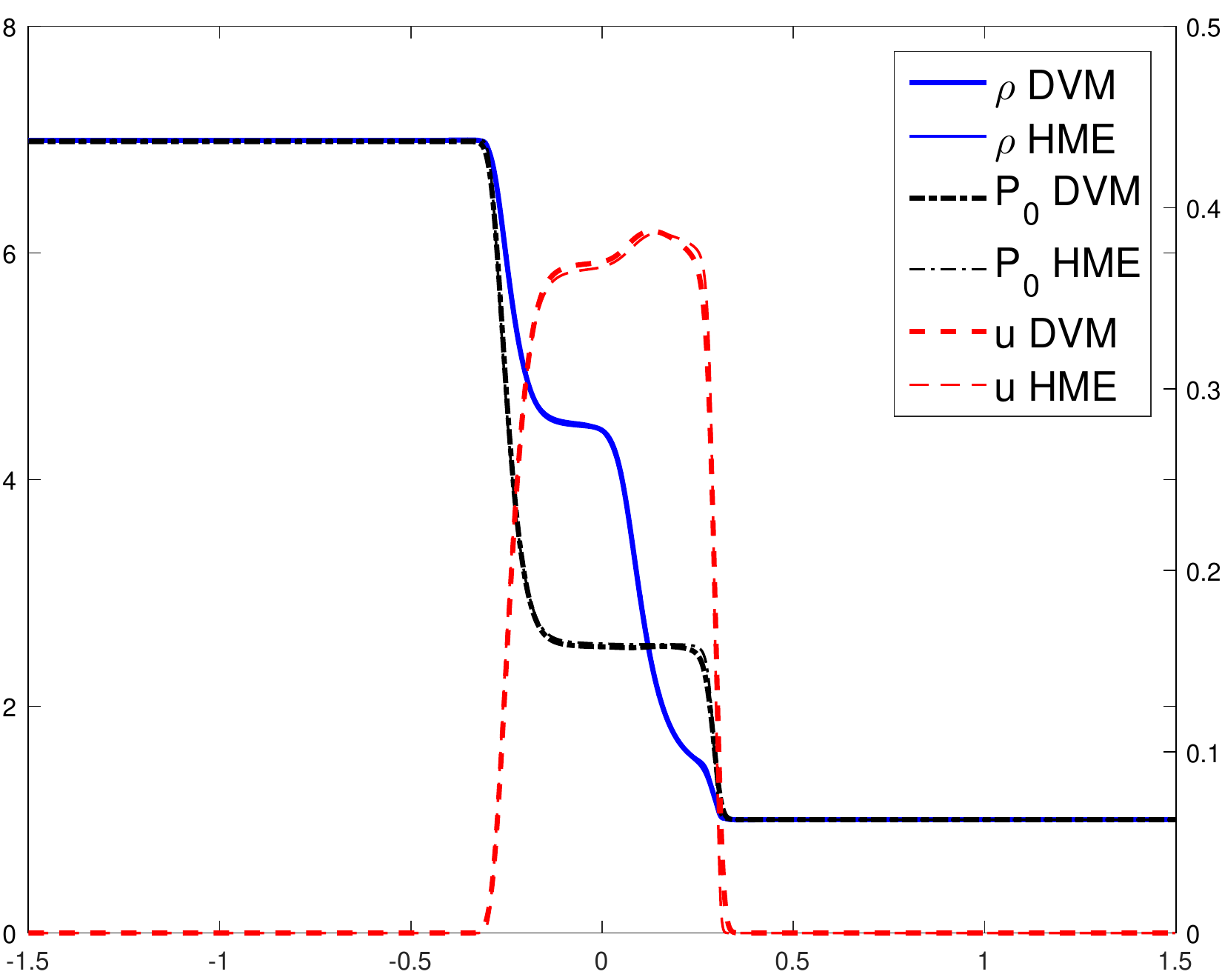}\
  }

  \centering
   \subfigure[M=4]
  {
  \includegraphics[width=1.7in]{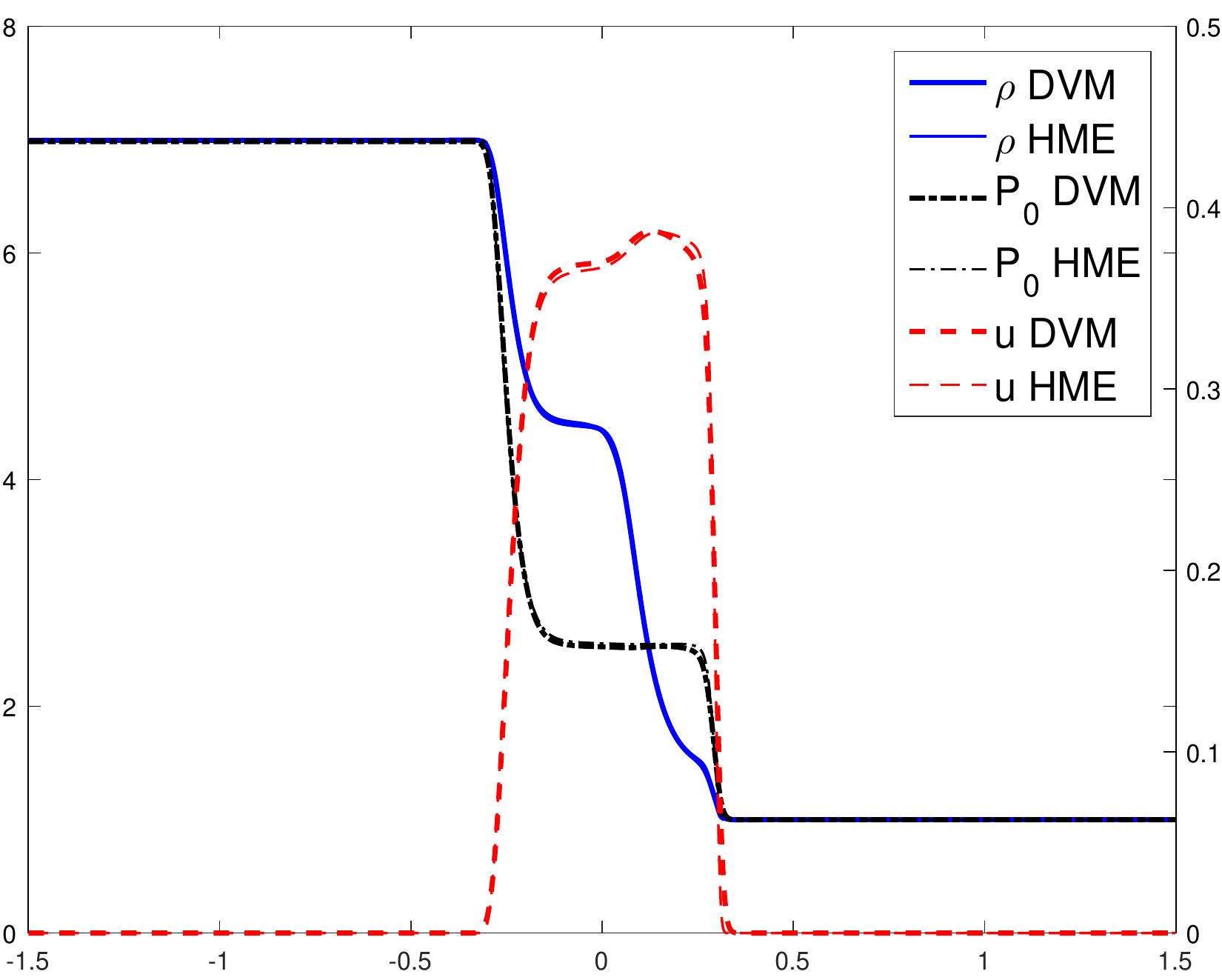}
  }
  \subfigure[M=5]
  {
  \includegraphics[width=1.7in]{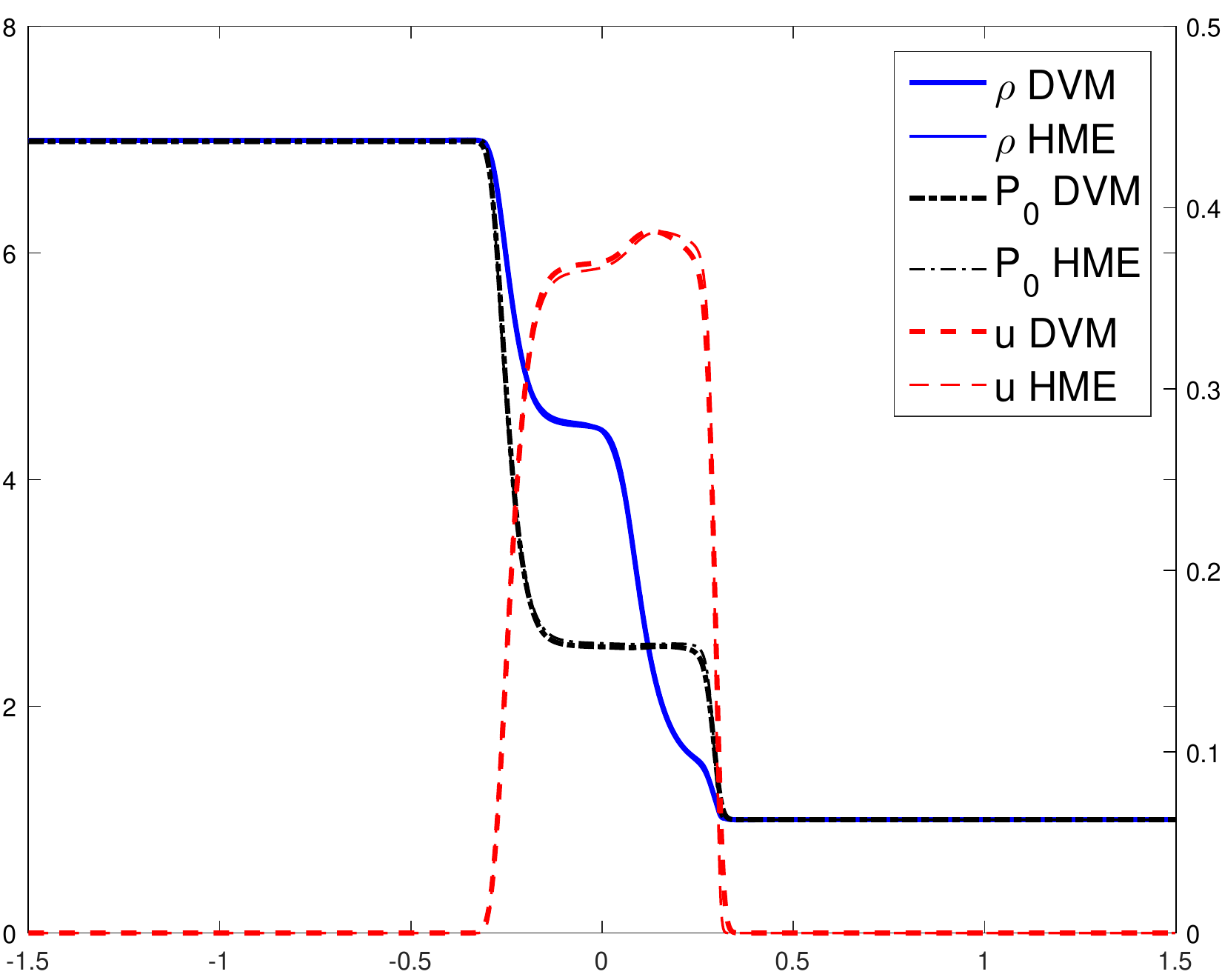}
  }
  \subfigure[M=6]
  {
  \includegraphics[width=1.7in]{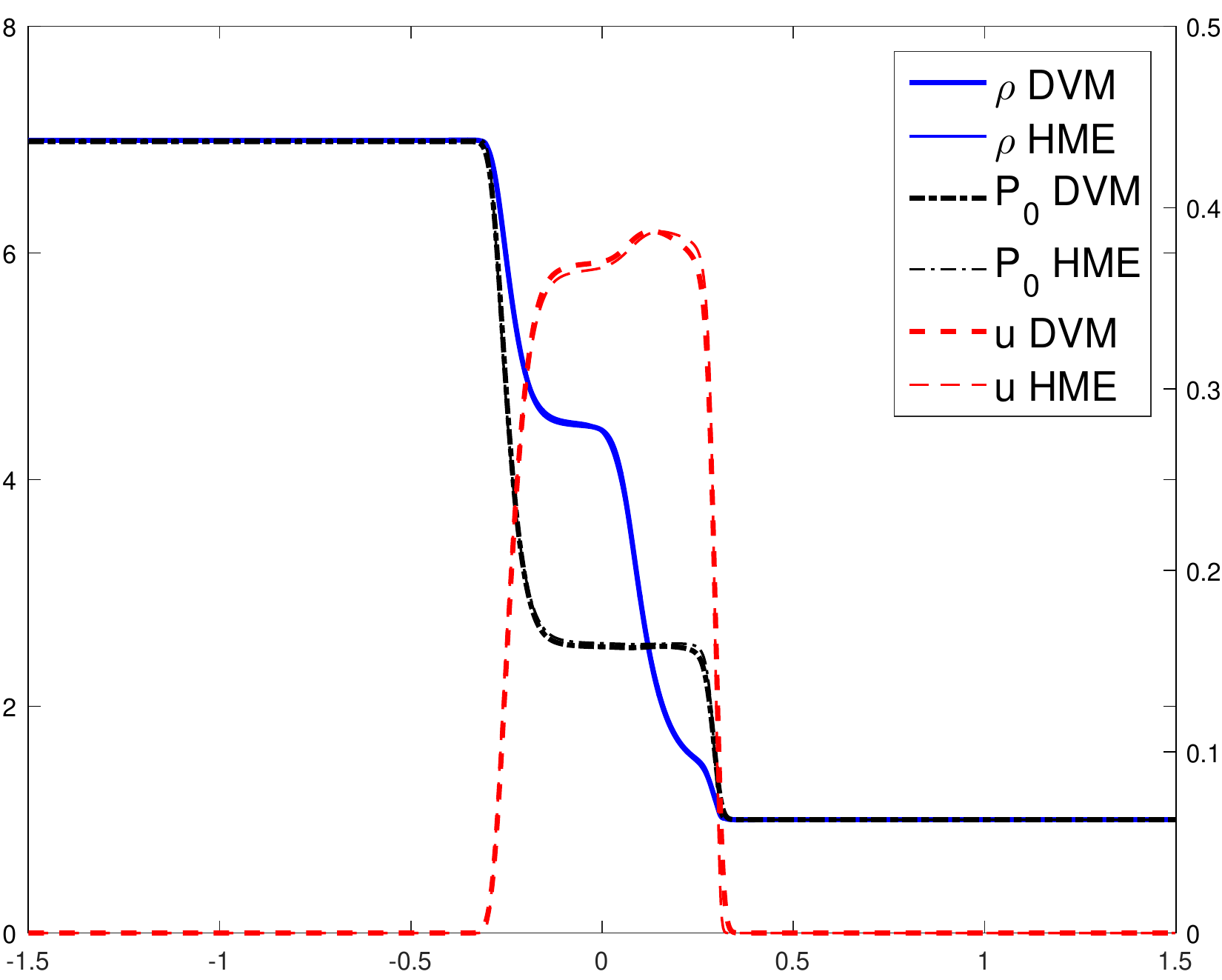}
  }

  \centering
   \subfigure[M=7]
  {
  \includegraphics[width=1.7in]{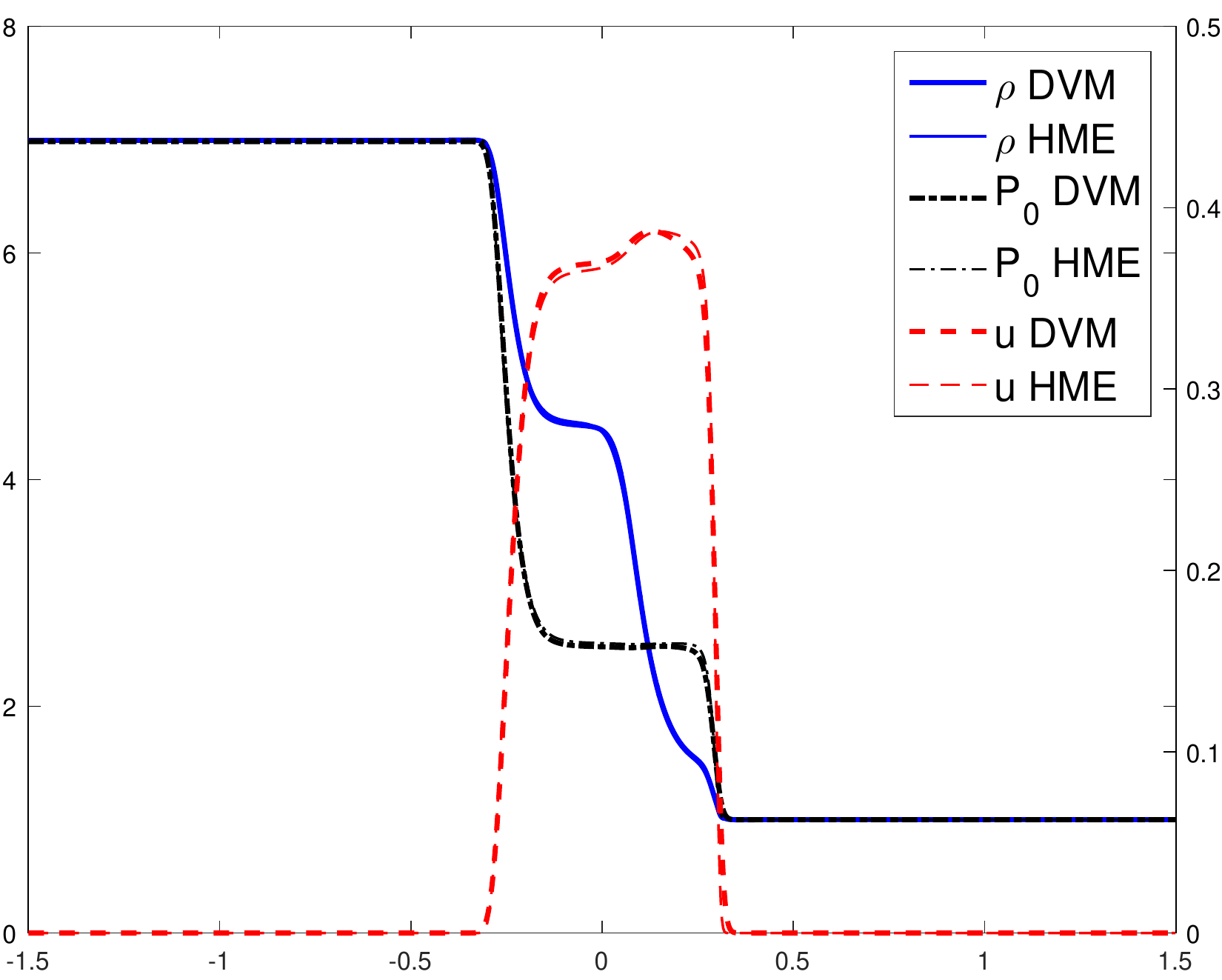}
  }
   \subfigure[M=8]
  {
  \includegraphics[width=1.7in]{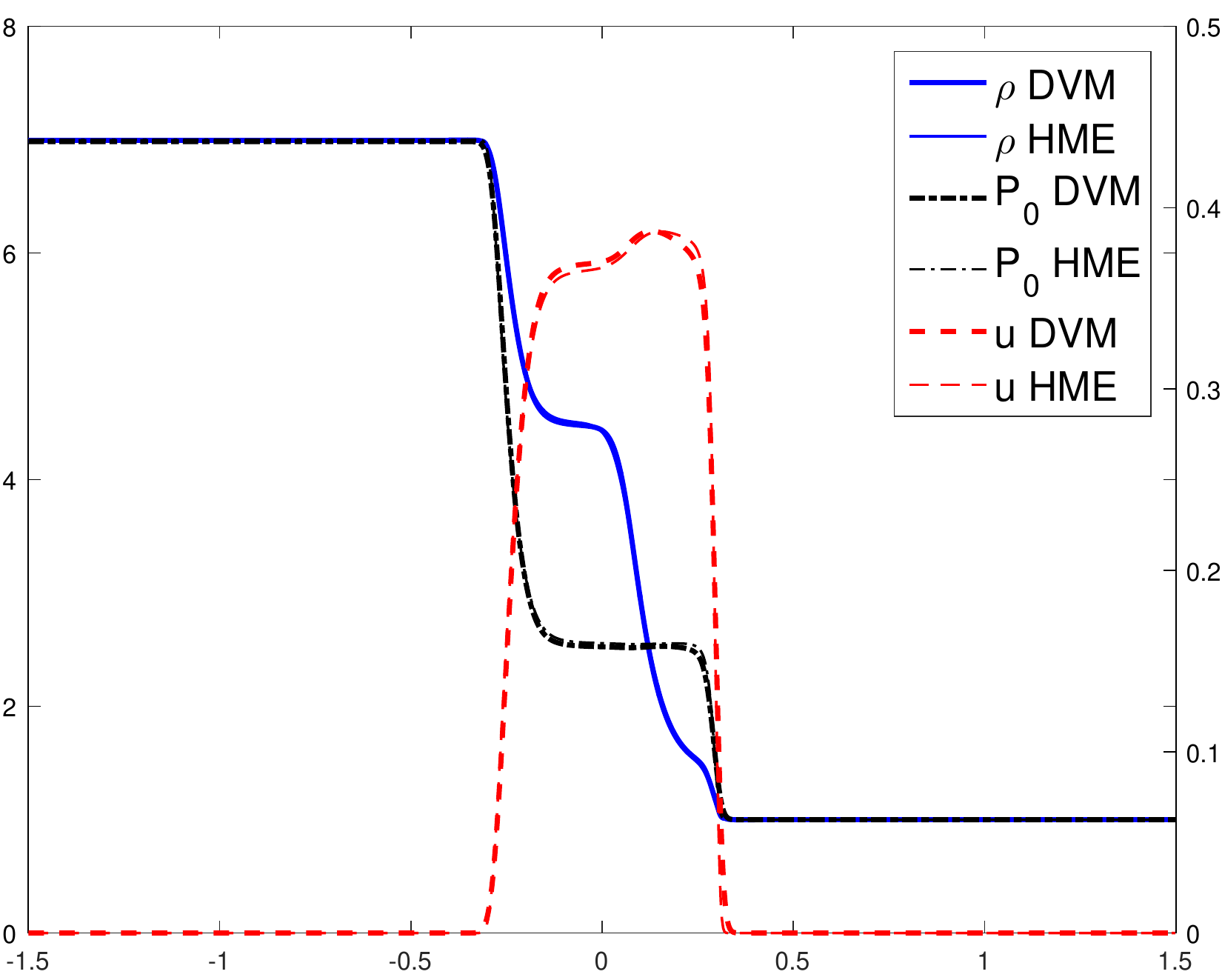}
  }
    \subfigure[M=9]
  {
  \includegraphics[width=1.7in]{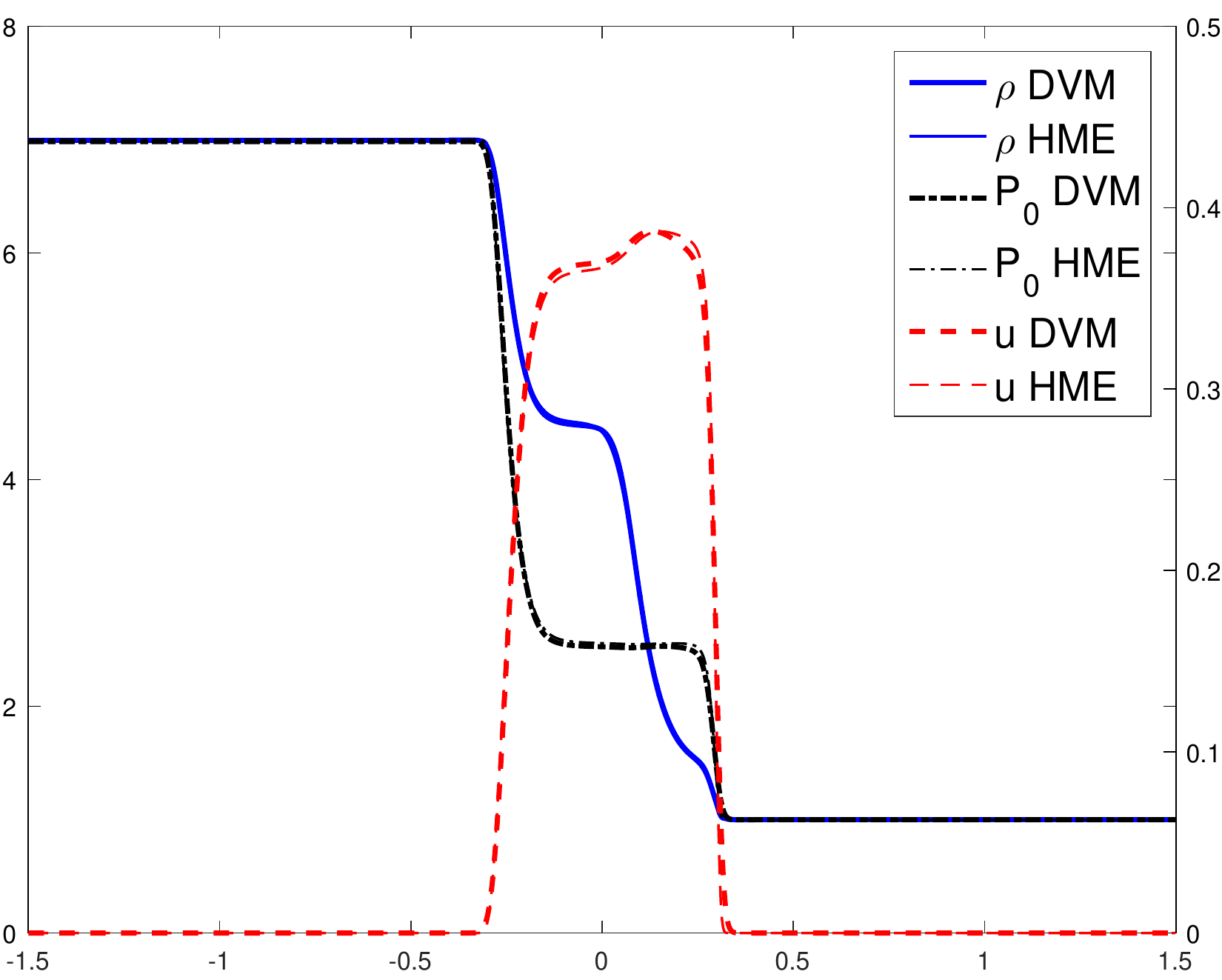}
   }

  \caption{\small Numerical results of the shock tube problem for $Kn=0.05$. The left $y$-axis is for $\rho$ and $P_{0}$, and the right $y$-axis
  is for $u$. The thin lines are the numerical results of the HME \eqref{eq:moment1}, and the thick lines are the results of DVM. The solid blue lines denote $\rho$, dashed red lines denote $u$, and dash-dotted black lines denote $P_{0}$.}
\end{figure}

Fig. \ref{fig:005} shows
the profiles of the density $\rho$, velocity $u$ and thermodynamic pressure $P_{0}$
at  $t=0.3$ obtained by using our scheme  \eqref{eq:advfinal}  and  \eqref{eq:semicol}
 with $M=1,2,\cdots,9$, 
 where  $Kn=0.05$, and the thin lines
are the numerical results of the HME \eqref{eq:moment1}, and the thick lines are the results of DVM, provided
as reference solutions. The solid lines denote $\rho$, dashed lines denote $u$, and dash-dotted lines denote $P_{0}$.

It is clear that the numerical solutions of the HME \eqref{eq:moment1} converge to the reference solution
of the special relativistic Boltzmann equation \eqref{eq:1DBoltz} as $M$ increases.
When \\$M=1$,  the contact discontinuity and shock wave {can} be obviously observed.
It is reasonable because the HME \eqref{eq:moment1} are the same as
the macroscopic RHD equations \eqref{eq:conser}.
When $M=2$, the discontinuities can also observed, but they have been damped.
When $M\geq 3$, the discontinuities are fully damped and the solutions are almost
in agreement with  the reference solutions.
 It is similar to the phenomena in the non-relativistic case \cite{DAMP:2001,1DB:2013}.

\begin{figure}
\label{fig:05}
  \centering
  \subfigure[M=1 RHD]
  {
  \includegraphics[width=1.7in]{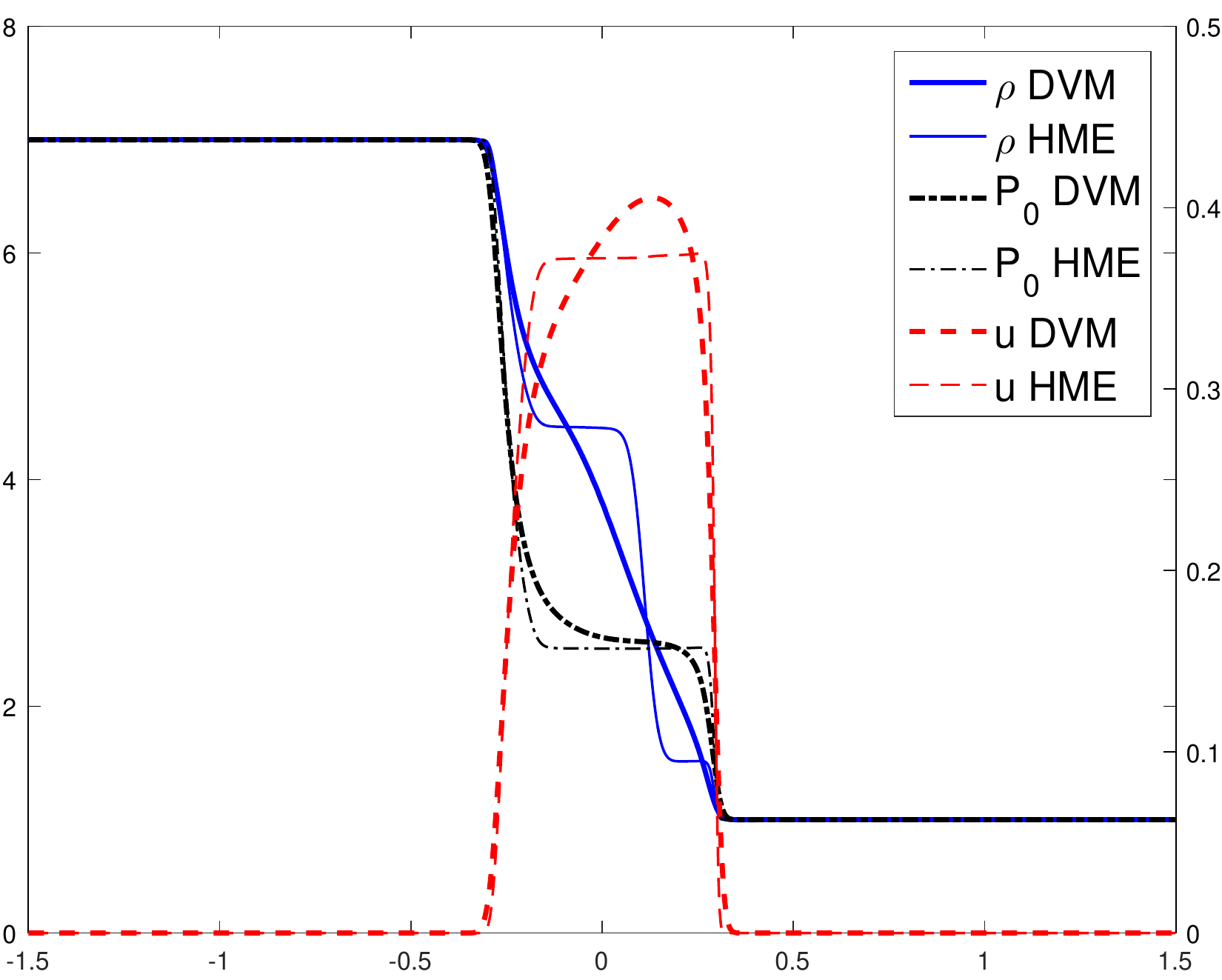}
  }
  \subfigure[M=2]
  {
  \includegraphics[width=1.7in]{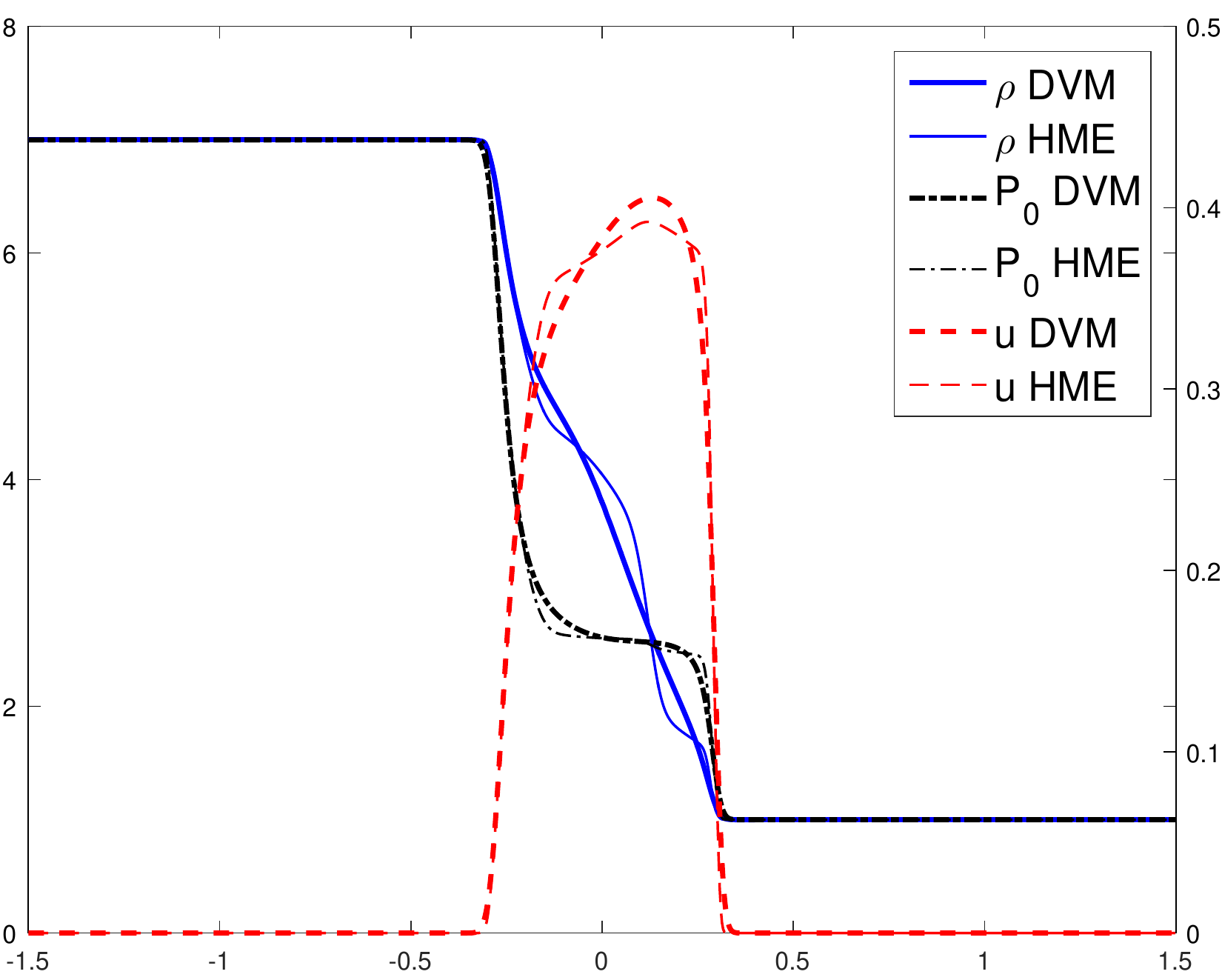}
   }
  \subfigure[M=3]
  {
  \includegraphics[width=1.7in]{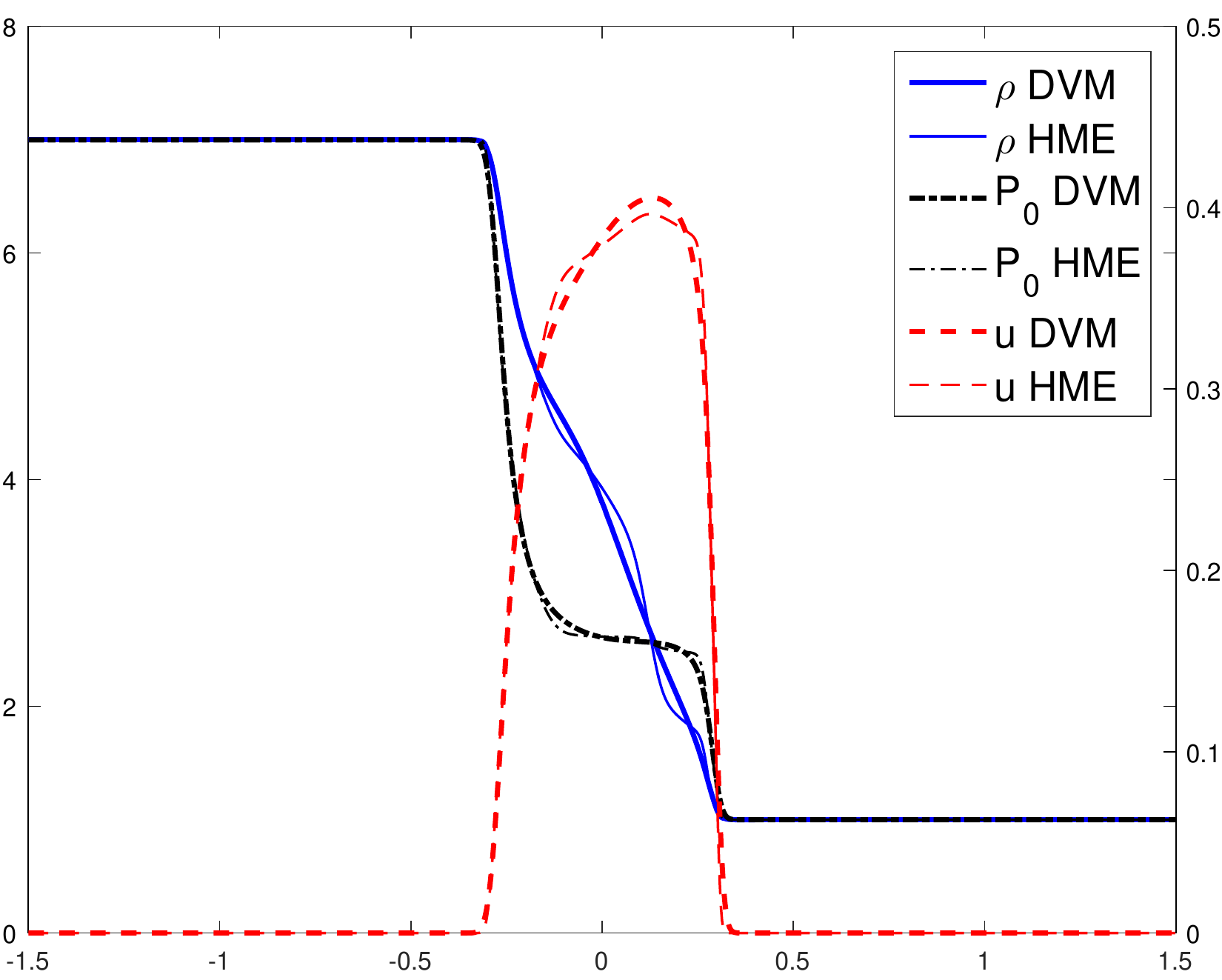}\
  }

  \centering
   \subfigure[M=4]
  {
  \includegraphics[width=1.7in]{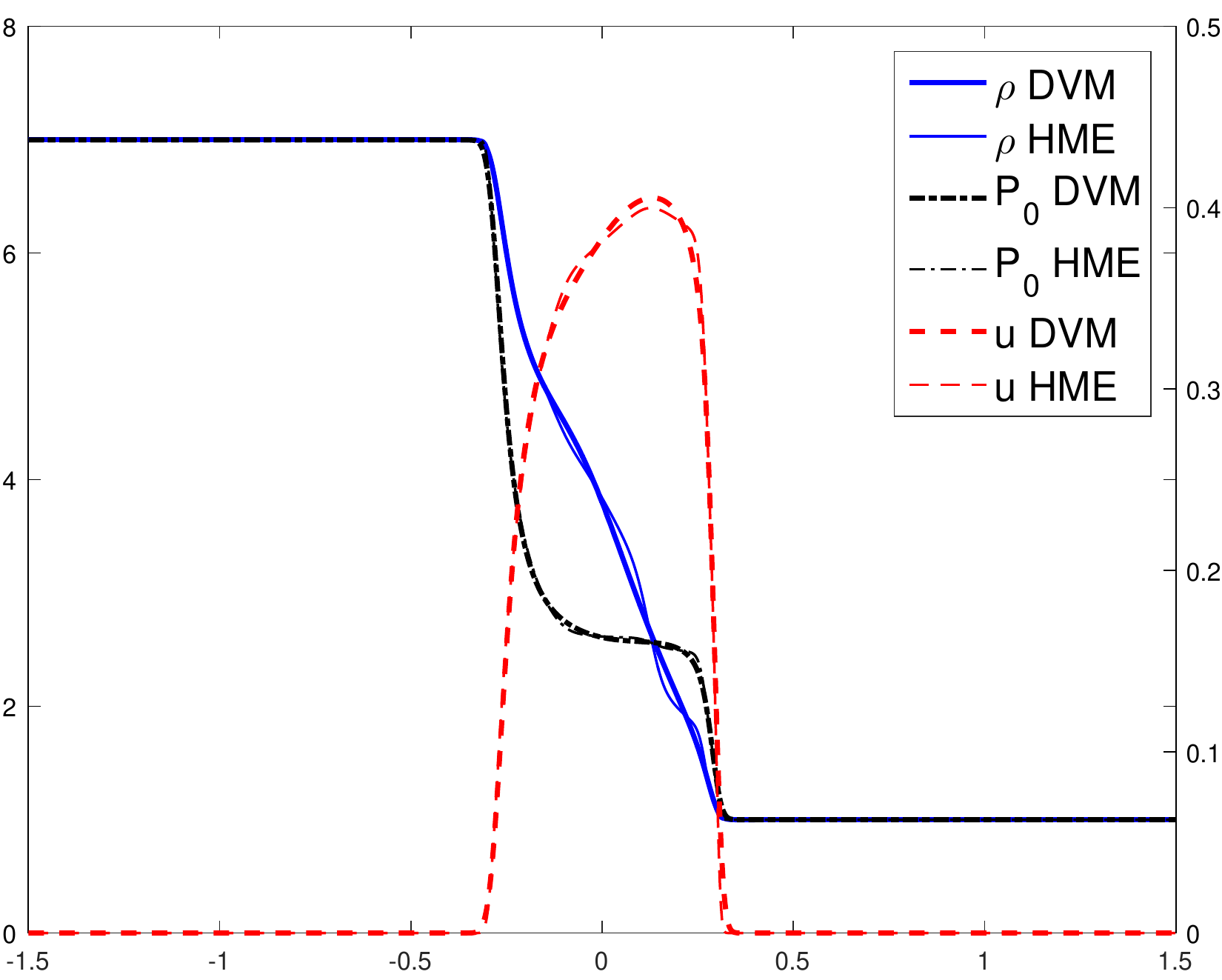}
  }
  \subfigure[M=5]
  {
  \includegraphics[width=1.7in]{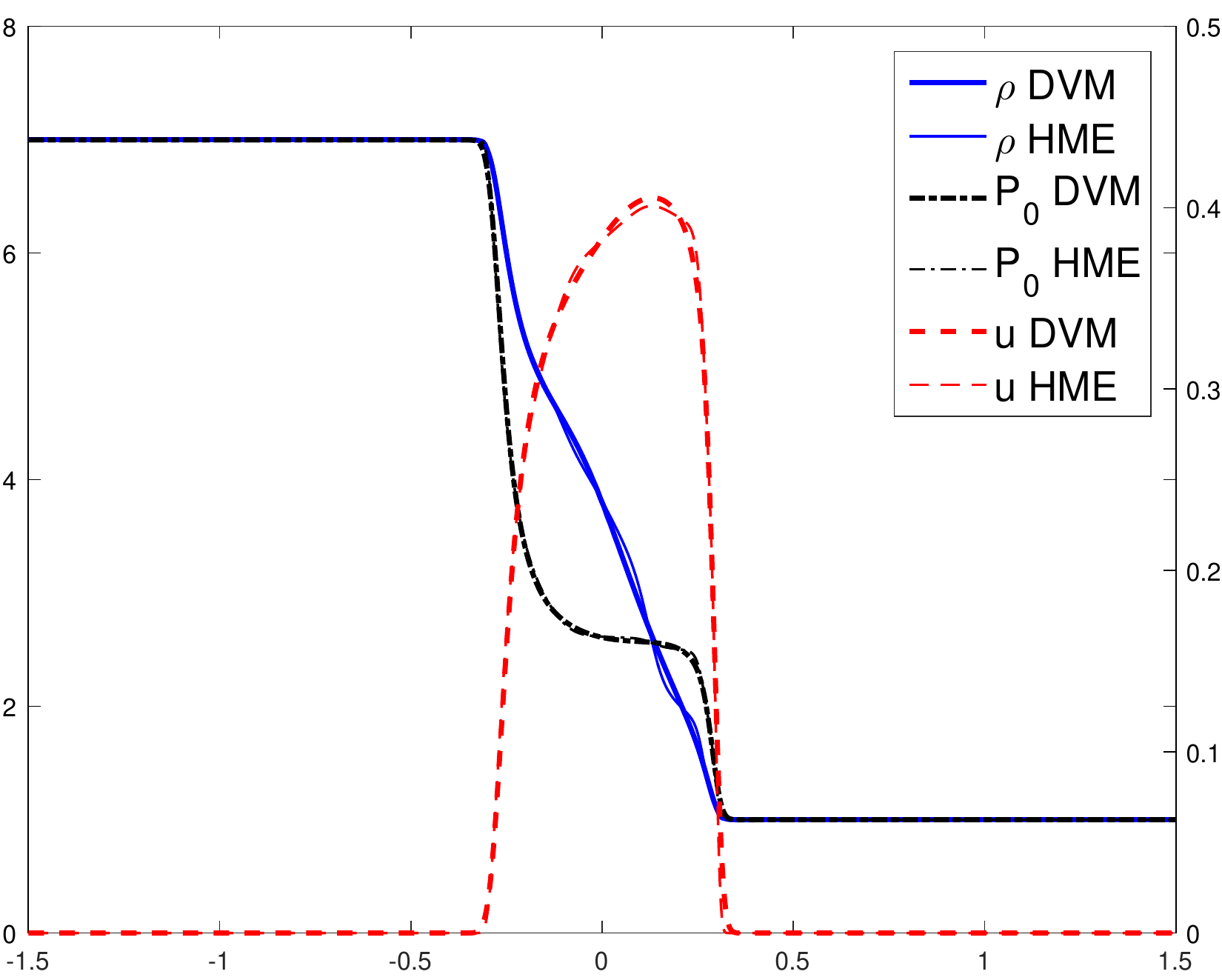}
  }
  \subfigure[M=6]
  {
  \includegraphics[width=1.7in]{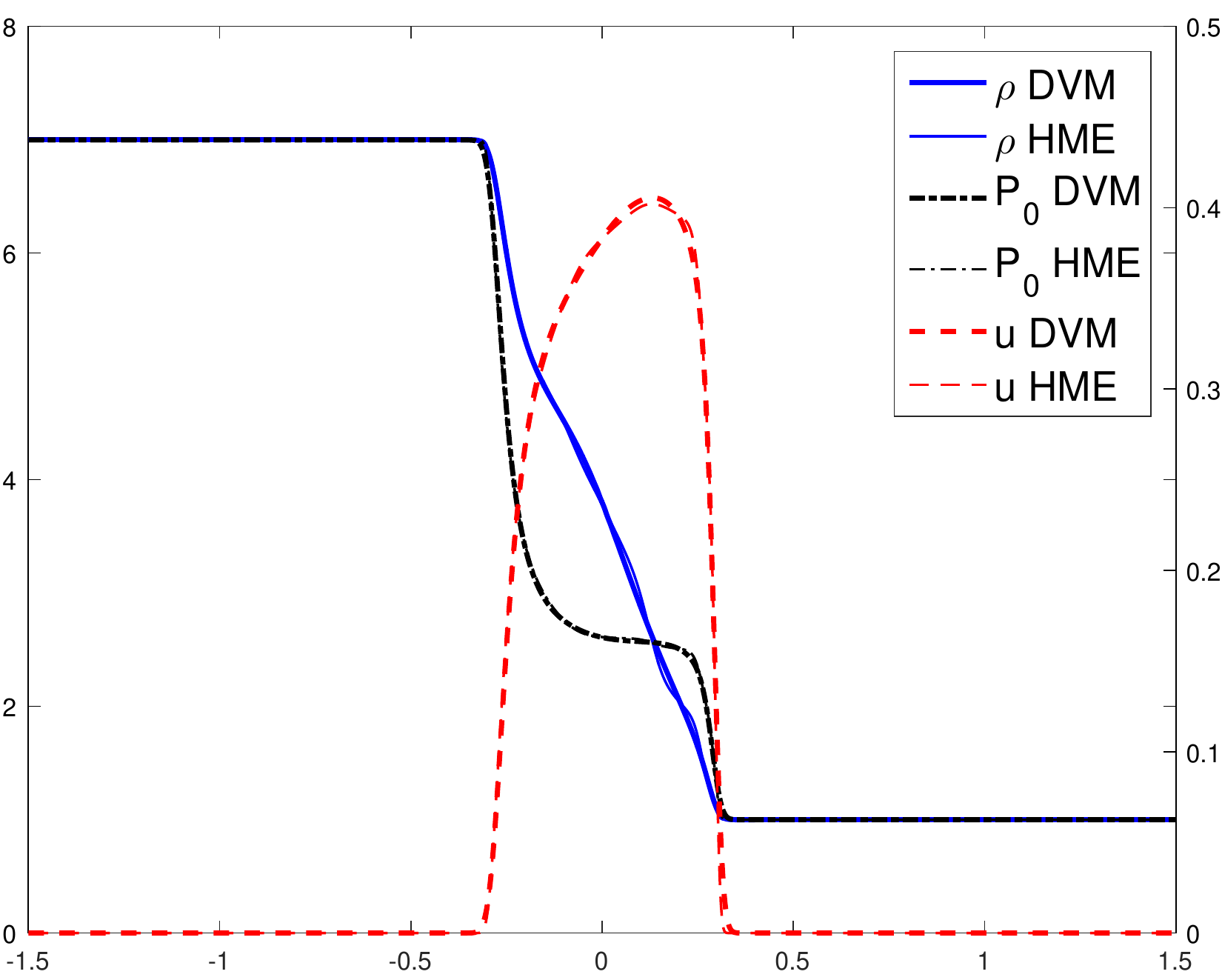}
  }

  \centering
   \subfigure[M=7]
  {
  \includegraphics[width=1.7in]{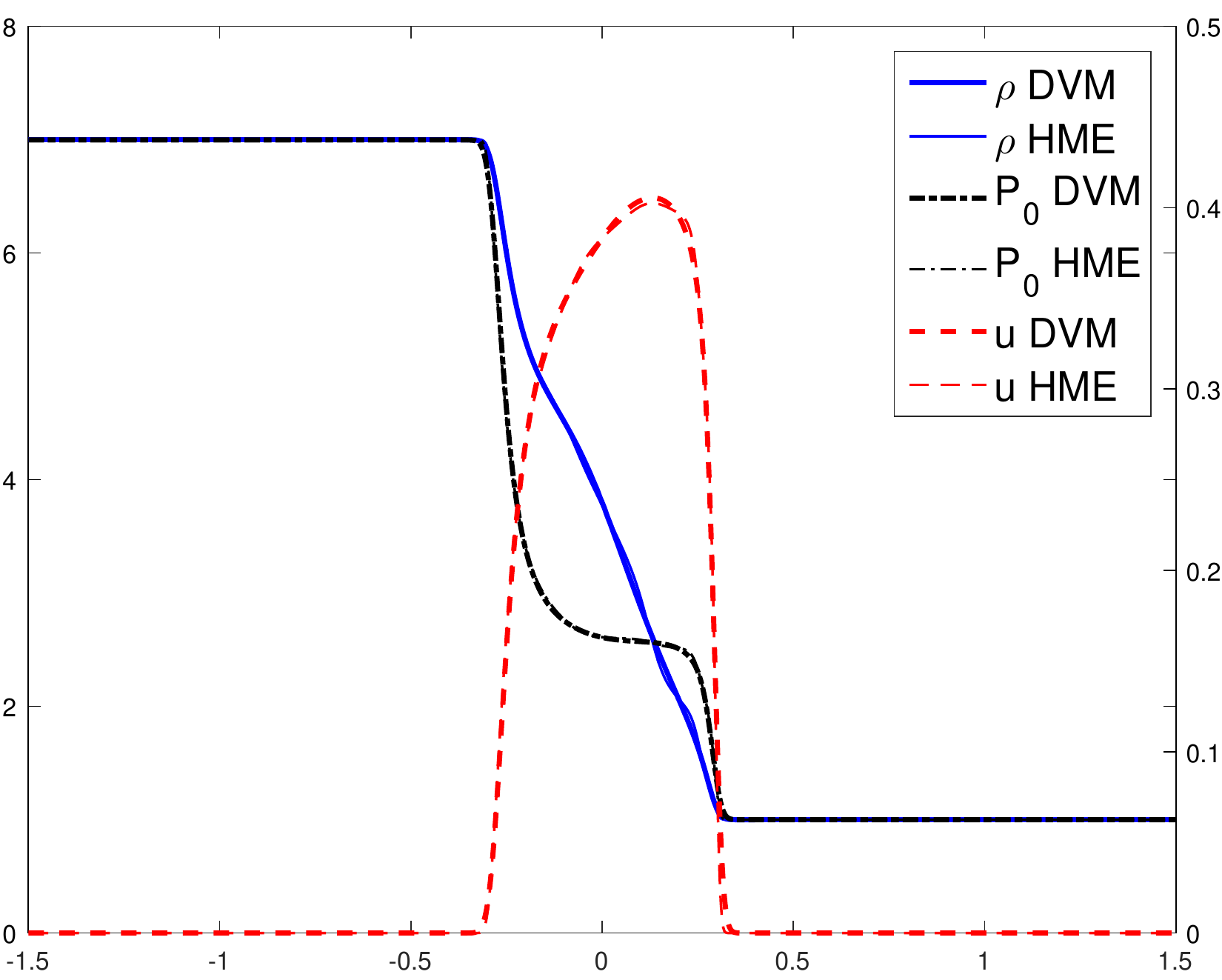}
  }
   \subfigure[M=8]
  {
  \includegraphics[width=1.7in]{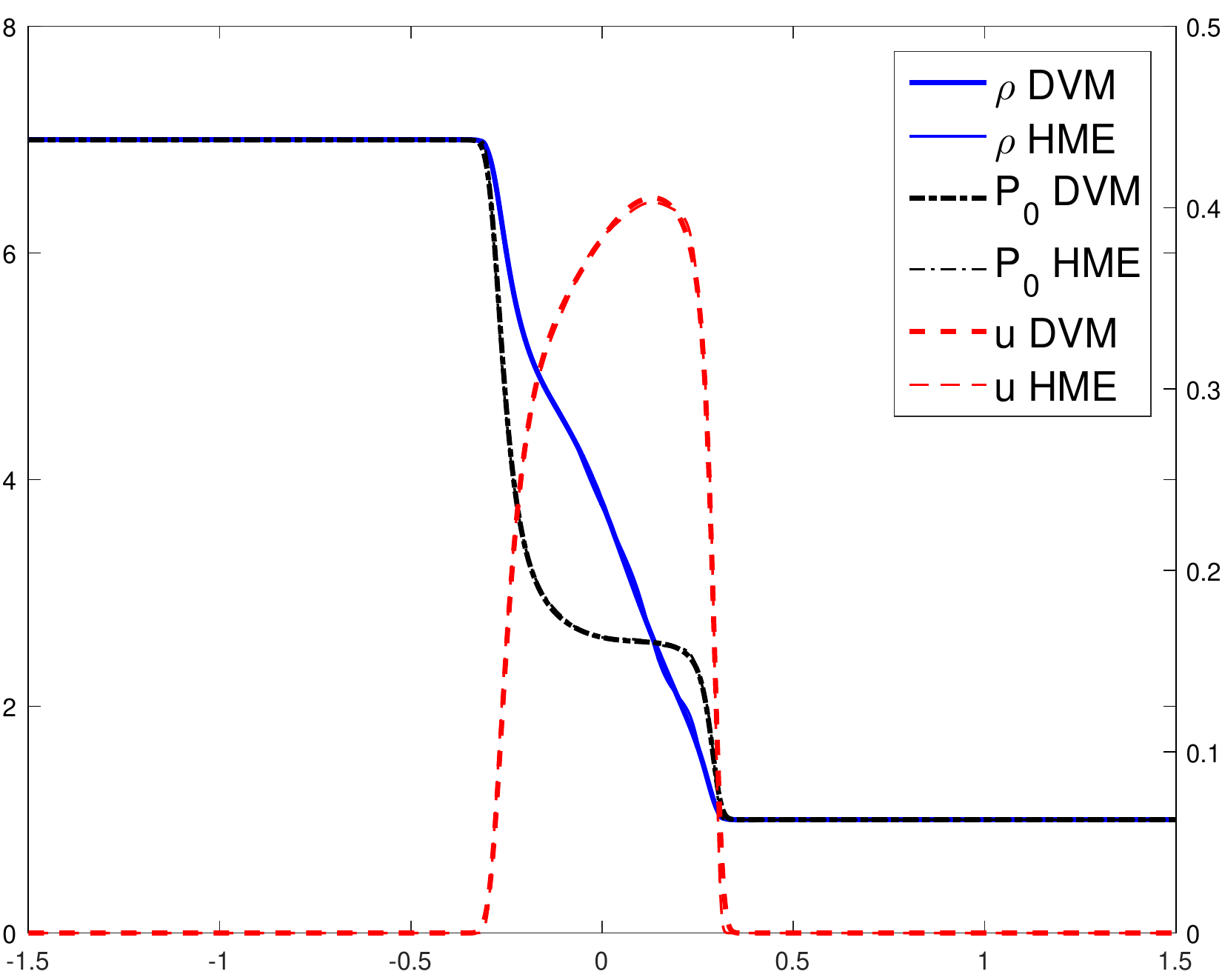}
  }
    \subfigure[M=9]
  {
  \includegraphics[width=1.7in]{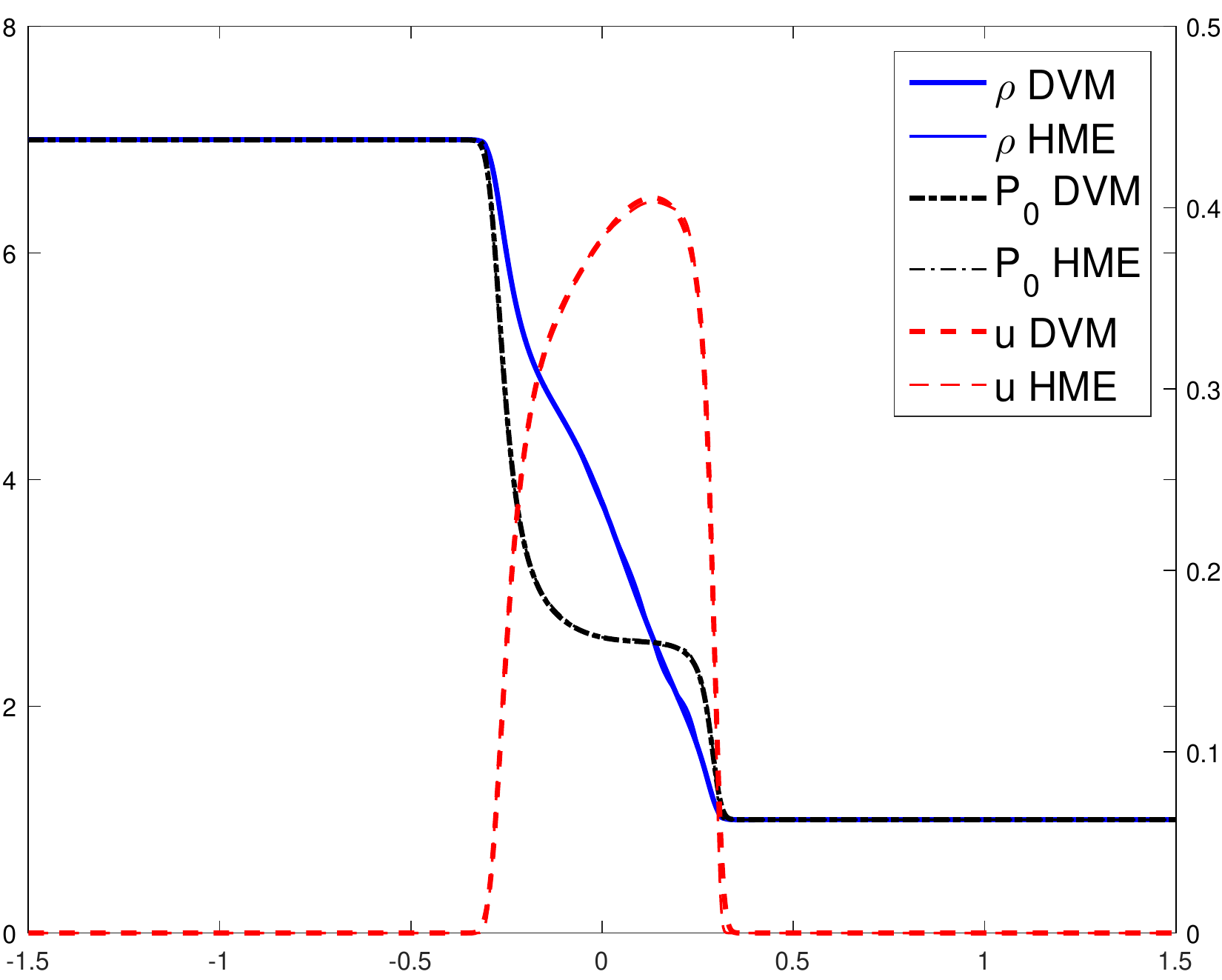}
   }

  \caption{\small Same as Fig. \ref{fig:005} except for $Kn=0.5$.}
\end{figure}

 The results at $t=0.3$ for the case of $Kn=0.5$ are shown in Fig. \ref{fig:05}. 
 The discontinuities are clearer than the case of $Kn=0.05$ when $M=1,2,\cdots 9$,
  and the convergence of the moment method can also be readily observed,
   but it is slower than the case of $Kn=0.05$.
   The contact discontinuities and  shock waves are obvious when $M\leq 2$,
  but when $M>6$, the discontinuities are fully damped and the
solutions are almost the same as the reference solutions.

\section{Conclusions}
\label{sec:conclud}
The paper derived the  arbitrary order globally hyperbolic moment system of
the one-dimensional (1D) special relativistic Boltzmann equation for the first time
 and studied the properties of the moment system:
 the  eigenvalues and their bound as well as  eigenvectors, hyperbolicity,
    characteristic fields, linear stability, and Lorentz covariance.
 The key contribution was the careful study of two families of the
 complicate Grad type orthogonal polynomials depending on a parameter.
 We derived the recurrence relations and derivative relations with respect to the independent variable and the parameter  respectively, and studied their zeros and coefficient matrices
 in the recurrence formulas.
%
Built on the knowledges of two families of the  Grad type orthogonal polynomials with a parameter,
 the model reduction method by the operator projection \cite{MR:2014} might be extended to the
1D special relativistic Boltzmann equation.

 A semi-implicit operator-splitting type numerical scheme
was presented for our hyperbolic moment system and  a Cauchy problem  was solved
to verify the convergence behavior of the moment method in comparison
with the discrete velocity method. The results showed that
the solutions of our hyperbolic moment system could converge to the solution of
the special relativistic Boltzmann equation as the order of the hyperbolic moment
system increases.

 {
  Now we are deriving the globally hyperbolic moment model of arbitrary order for the 3D special relativistic Boltzmann equation.
    Moreover, it is interesting to  develop  robust, high order accurate numerical schemes for the moment system
and   find other basis for the derivation of  moment system with some good property, e.g. non-negativity.}

%




\section*{Acknowledgements}
This work was partially supported by
{ the Special Project on High-performance Computing under the National Key R\&D Program (No. 2016YFB0200603),
Science Challenge Project (No. JCKY2016212A502), and}
the National Natural Science Foundation
of China (Nos.  91330205, 91630310, \& 11421101).

\newpage
\begin{appendices}

\section{Proofs in Section \ref{sec:RB}}
\label{sec:App2}
\subsection{Proof of Theorem \ref{thm:admissible}}
\label{p:thm:admissible}
\begin{proof}
		For the nonnegative distribution  $f(x,p,t)$,  which is not identically zero,
		using \eqref{eq:NTab} gives
	\[
	T^{\alpha\alpha}>0, \alpha=0,1;\quad
T^{00}+T^{11}{\pm}2T^{01}={c}\int_{\mathbb{R}}(p^{0}{\pm} p^{1})^2f\frac{dp}{p^{0}}>0,
	\]
 which implies the first inequality in \eqref{eq:admissible}.

Using the definition of $\Delta^{\alpha\beta}$ in \eqref{eq:sym-tensor}
and  the tensor decomposition of $T^{\alpha\beta}$ in  \eqref{eq:Tabdiv}
{gives} \eqref{eq:admu1},
which is a {quadratic} equation with respect to $u$.
The first   inequality in \eqref{eq:admissible}	tells us that
  \eqref{eq:admu1} has two different solutions whose product is equal to ${c^{2}}$,
 {while one of them with a smaller absolute value is \eqref{eq:admu}}.
	
Using further \eqref{eq:NTab} gives
	\[
	{N^{0}-c^{-1}u N^{1}=c\int_{\mathbb{R}}(p^{0}-c^{-1}up^{1})f\frac{dp}{p^{0}}}>0,
	\]
	i.e.  the second  inequality in \eqref{eq:admissible},
	and then using  the tensor decomposition of $N^\alpha$ in  \eqref{eq:Ndiv} gives
	\[
	  \rho={c^{-1}m\frac{N^{0}-c^{-1}u N^{1}}{\sqrt{1-c^{-2}u^2}}}>0.
	\]

Using the second identity in \eqref{eq:condition},  the expression of
	$\varepsilon_0$ in \eqref{eq:variable0}, and
	 \eqref{eq:Tabdiv} {gives} \eqref{eq:admT}.

And the inequality ${E\geq mc^{2}}$ holds because
\[
E=U_{\alpha}p^{\alpha}={(1-c^{-2}u^2)^{-\frac{1}{2}}\left(c\sqrt{m^2c^2+p^2}-up\right)}>0,
\]
and
\[
{E^2-m^2c^4=(1-c^{-2}u^2)^{-1}(u\sqrt{m^2c^2+p^2}-cp)^2=\left(\frac{c^2}{U^{0}}p_{<1>}\right)^2}\geq0.
\]
Thus
\[
{ T^{00}-c^{-1}uT^{01}-c^{2}\rho=\frac{T^{00}- 2c^{-1}u T^{01}+ c^{-2}u^2 T^{11}}{1-c^{-2}u^2}
-c^{2}\rho=c^{-1}\int_{\mathbb{R}}E(E-mc^2)f\frac{dp^1}{p^{0}}}>0,
\]
the third   inequality in \eqref{eq:admissible}  holds, and thus
implies   that $G(\theta^{-1})-\theta>1$ for $\theta \in (0,+\infty)$.

On the other hand,
one has 
	\[
	\lim_{\theta\rightarrow0}\left(G(\theta^{-1})-\theta\right)=1,\
	\lim_{\theta\rightarrow +\infty}\left(G(\theta^{-1})-\theta\right)=\lim_{\theta\rightarrow +\infty}\theta= +\infty,
	\]
	and
{\begin{align*}
	\frac{\partial (G(\theta^{-1})-\theta)}{\partial \theta}
	=&-\theta^{-2}
	\left(G(\theta^{-1})^2-3G(\theta^{-1})\theta+\theta^{2}-1\right)=:\tilde{\psi}(G(\theta^{-1}),\theta).
	\end{align*}
Because
\begin{align*}
0&<c^{-1}\int_{\mathbb{R}}(E-mc^2)f^{(0)}\frac{dp}{p^{0}}=-m^{-1}\rho (G(\theta^{-1})-2\theta-1),\\
0&<c^{-1}\int_{\mathbb{R}}(E-mc^2)^2f^{(0)}\frac{dp}{p^{0}}=\rho c^{2}(2G(\theta^{-1}) - 3\theta-2),\\
0&<c^{-1}\int_{\mathbb{R}}(E-mc^2)^3f^{(0)}\frac{dp}{p^{0}}=-\rho mc^{4}((4 - \theta)G(\theta^{-1}) - 5\theta-4),
\end{align*}
one {obtain}s
\begin{align*}
\frac{3}{2}\theta+1<  G(\theta^{-1})
 <
 \begin{cases}
 \min\left\{2\theta+1,(4-\theta)^{-1}(5\theta+4)\right\},&  0<\theta<4,\\
 2\theta+1,& \theta\geq4,
  \end{cases}
\end{align*}
which   is equivalent  to the following inequality
\begin{align*}
\frac{3}{2}\theta+1&<G(\theta^{-1})<
\begin{cases}
(4-\theta)^{-1}(5\theta+4),  & 0<\theta<1,\\
2\theta+1, & \theta\geq1.
\end{cases}
\end{align*}
Thus, one has
\begin{align*}
\tilde{\psi}(G(\theta^{-1}),\theta)>
\begin{cases} \tilde{\psi}\left(2\theta+1,\theta\right)>\theta^{3}(\theta-1)>0, &\theta\geq1,\\
\\
\tilde{\psi}\left((4-\theta)^{-1}(5\theta+4),\theta\right)>(4-\theta)^{-2}\theta^4(\theta+8)(1-\theta)>0, &\theta<1,
\end{cases}\end{align*}
i.e.
\[
\frac{\partial (G(\theta^{-1})-\theta)}{\partial \theta}>0,
\]}
which implies that 	$G(\theta^{-1})-\theta$  is a strictly monotonic function  of $\theta$ in the interval $(0,+\infty)$.

%
Thus   \eqref{eq:admT} has a unique solution in the interval $(0,+\infty)$.
	The proof is completed.
	\qed\end{proof}

\subsection{Proof of Theorem \ref{lem:admissible1}}
\label{p:lem:admissible1}
	\begin{proof}
	Under Theorem \ref{thm:admissible},
		for the nonnegative distribution  $f(x,p,t)$,  which is not identically zero,
one {obtain}s $\{\rho,u,\theta\}$ satisfying
\begin{equation}
\label{eq:positivity}
		\rho>0, \quad |u|<{c}, \quad \theta>0.
\end{equation}		
Due to the last equations in \eqref{eq:variable} and \eqref{eq:variable0}, one {obtain}s
		\[
		{\Pi=-\int_{\mathbb{R}}\Delta_{\alpha\beta}p^{\alpha}p^{\beta}f\frac{dp}{p^{0}}-c^{2}\rho \theta=c^{-1}\int_{\mathbb{R}}(E^2-m^2c^4)f\frac{dp}{p^{0}}-\rho c^{2}\theta>-\rho c^{2}\theta},
		\]
which 	 completes 	the proof.
		\qed\end{proof}

\section{Proofs in Section \ref{sec:orth}}
\label{sec:App3}
\subsection{Proof of Theorem \ref{thm:huxiang}}
\label{p:thm:huxiang}
\begin{proof}
(i)  For $k\leq n+2$,  taking the inner product with respect to $\omega^{(0)}$ between the polynomials  $P_{k}^{(0)}(x;\zeta)$
and $(x^2-1)P_{n}^{(1)}(x;\zeta)$ 
gives
\begin{align*}
\left((x^2-1)P_{n}^{(1)},P_{n+2}^{(0)}\right)_{\omega^{(0)}}&=\left(c_{n}^{(1)}x^{n+2},P_{n+2}^{(0)}\right)_{\omega^{(0)}}=
\frac{c_{n}^{(1)}}{c_{n+2}^{(0)}}\left(P_{n+2}^{(0)},P_{n+2}^{(0)}\right)_{\omega^{(0)}}=r_{n+1},\\
\left((x^2-1)P_{n}^{(1)},P_{n+1}^{(0)}\right)_{\omega^{(0)}}&=
\left(c_{n}^{(1)}\left(x^{n+2}-\sum_{i=1}^{n+2}x_{i,n+2}^{(0)}x^{n+1}+
\left(\sum_{i=1}^{n+2}x_{i,n+2}^{(0)}-\sum_{i=1}^{n}x_{i,n}^{(1)}\right)x^{n+1}\right),P_{n+1}^{(0)}\right)_{\omega^{(0)}}\\
 &=r_{n+1}\left(P_{n+2}^{(0)},P_{n+1}^{(0)}\right)_{\omega^{(0)}}+
 q_{n}\left(P_{n+1}^{(0)},P_{n+1}^{(0)}\right)_{\omega^{(0)}}=q_{n},\\
\left((x^2-1)P_{n}^{(1)},P_{n+1}^{(0)}\right)_{\omega^{(0)}}&=\left(P_{n}^{(1)},c_{n+1}^{(0)}
\left(x^{n+1}-\sum_{i=1}^{n+1}x_{i,n+1}^{(1)}x^{n}+
\left(\sum_{i=1}^{n+1}x_{i,n+1}^{(1)}-\sum_{i=1}^{n+1}x_{i,n+1}^{(0)}\right)x^{n}\right)\right)_{\omega^{(1)}}\\
 &=p_{n+1}\left(P_{n}^{(1)},P_{n+1}^{(1)}\right)_{\omega^{(1)}}+\frac{c_{n+1}^{(0)}}{c_{n}^{(1)}}\left(\sum_{i=1}^{n+1}x_{i,n+1}^{(1)}-
 \sum_{i=1}^{n+1}x_{i,n+1}^{(0)}\right)\left(P_{n}^{(1)},P_{n}^{(1)}\right)_{\omega^{(1)}}\\
 &=\frac{c_{n+1}^{(0)}}{c_{n}^{(1)}}
\sum_{i=1}^{n+1} \left(x_{i,n+1}^{(1)}-x_{i,n+1}^{(0)}\right)=q_{n},
 \\
\left((x^2-1)P_{n}^{(1)},P_{n}^{(0)}\right)_{\omega^{(0)}}&
=\left(P_{n}^{(1)},P_{n}^{(0)}\right)_{\omega^{(1)}}=
   \left(P_{n}^{(1)},c_{n}^{(0)}x^{n}\right)_{\omega^{(1)}}
   =p_{n}\left(P_{n}^{(1)},P_{n}^{(1)}\right)_{\omega^{(1)}}=
   p_{n},\\
    \left((x^2-1)P_{n}^{(1)},P_{k}^{(0)}\right)_{\omega^{(0)}}
    &=\left(P_{n}^{(1)},P_{k}^{(0)}\right)_{\omega^{(1)}}=0, \quad k\leq n-1,
\end{align*}
Substituting them into  \eqref{eq:P01orth-3} gives  \eqref{eq:recP01}.

(ii) Taking the inner product with respect to $\omega^{(1)}$ between $P_{n+1}^{(0)} (x;\zeta)$
 and $P_{k}^{(1)} (x;\zeta)$  with $k\leq n+1$ 
\begin{align*}
\left(P_{n+1}^{(0)},P_{n+1}^{(1)}\right)_{\omega^{(1)}}&= \left(c_{n+1}^{(0)}x^{n+1},P_{n+1}^{(1)}\right)_{\omega^{(1)}}
=p_{n+1}\left(P_{n+1}^{(1)},P_{n+1}^{(1)}\right)_{\omega^{(1)}}=p_{n+1},\\
\left(P_{n+1}^{(0)},P_{n}^{(1)}\right)_{\omega^{(1)}}&= \left(P_{n+1}^{(0)},(x^2-1)P_{n}^{(1)}\right)_{\omega^{(0)}}=q_{n},
\\
\left(P_{n+1}^{(0)},P_{n-1}^{(1)}\right)_{\omega^{(1)}}&=\left(P_{n+1}^{(0)},(x^2-1)P_{n-1}^{(1)}\right)_{\omega^{(0)}}=
r_{n}\left(P_{n+1}^{(0)},P_{n+1}^{(0)}\right)_{\omega^{(0)}}=r_{n},\\
    \left(P_{n+1}^{(0)},P_{k}^{(1)}\right)_{\omega^{(1)}}
    &=\left(P_{n+1}^{(0)},(x^2-1)P_{k}^{(1)}\right)_{\omega^{(0)}}=0, \quad k\leq n-2.
\end{align*}
Similarly, substituting them into  \eqref{eq:P01orth-3} gives    \eqref{eq:recP10}.

(iii) If using \eqref{eq:recP0P1} to eliminate $P_{n+2}^{(0)}$ and $P_{n+1}^{(1)}$ in \eqref{eq:recP01} and \eqref{eq:recP10} respectively, then one {obtains}
 \begin{equation*}
    (x^2-1)P_{n}^{(1)}=\tilde{p}_{n}(x+\tilde{q}_{n})P_{n+1}^{(0)}+\tilde{r}_{n}P_{n}^{(0)},\quad
    P_{n+1}^{(0)}=\frac{1}{\tilde{\tilde{p}}_{n}}(x-\tilde{\tilde{q}}_{n})P_{n}^{(1)}-\frac{a_{n-1}^{(1)}}{a_{n}^{(0)}}\tilde{\tilde{r}}_{n}P_{n-1}^{(1)},
\end{equation*}
with 
\[
    \tilde{p}_{n}=\frac{r_{n+1}}{a_{n+1}^{(0)}}=\frac{c_{n}^{(1)}}{c_{n+1}^{(0)}}=\frac{a_{n}^{(1)}}{p_{n+1}}=\tilde{\tilde{p}}_{n},
\]
\begin{equation*}
  \tilde{q}_{n}=\frac{1}{\tilde{p}_{n}}q_{n}-b_{n+1}^{(0)}=\sum_{i=1}^{n+1}x_{i,n+1}^{(0)}-\sum_{i=1}^{n}x_{i,n}^{(1)}=b_{n}^{(1)}-\tilde{p}_{n}q_{n}=\tilde{\tilde{q}}_{n},
\end{equation*}
\begin{equation*}
 \tilde{r}_{n}=p_{n}-\tilde{p}_{n}a_{n}^{(0)}={p_{n}(1-\tilde{p}_{n}^{2})}
 =\frac{a_{n}^{(0)}}{a_{n-1}^{(1)}}\left(-r_{n}+\frac{1}{\tilde{p}_{n}}a_{n-1}^{(1)}\right)=
 \tilde{\tilde{r}}_{n}.
\end{equation*}
The proof is completed.
\qed\end{proof}

\subsection{Proof of Theorem \ref{thm:derivezeta}}
\label{p:thm:derivezeta}
\begin{proof}
With the aid of  definition and   recurrence relation of the second kind modified
Bessel function in \eqref{eq:bessel} and \eqref{eq:besselrec}, one has
\begin{align*}
\frac{\partial }{\partial \zeta}\omega^{(\ell)}(x;\zeta)
&=\frac{K_{2}(\zeta)+K_{0}(\zeta)-2xK_{1}(\zeta)}{2K_{1}(\zeta)}\left(\frac{1}{K_{1}(\zeta)}(x^2-1)^{\ell-\frac{1}{2}}\exp(-\zeta x)\right)\\
&=\left(G(\zeta)-\zeta^{-1}-x\right)\omega^{(\ell)}(x;\zeta).
\end{align*}
Taking the partial derivative of both sides of identities
\[
    \left(P_{n+1}^{(\ell)},P_{k}^{(\ell)}\right)_{\omega^{(\ell)}}=\delta_{n+1,k}, k=0,\cdots,n+1,
\]
 with respect to $\zeta$ and using \eqref{eq:abn} gives
\begin{align*}
    \frac{\partial }{\partial \zeta}
    \left(P_{n+1}^{(\ell)},P_{n+1}^{(\ell)}\right)_{\omega^{(\ell)}}
    =&2\left(\frac{\partial }{\partial \zeta}P_{n+1}^{(\ell)},P_{n+1}^{(\ell)}\right)_{\omega^{(\ell)}}
    +\left(G(\zeta)-\zeta^{-1}\right)\left(P_{n+1}^{(\ell)},P_{n+1}^{(\ell)}\right)_{\omega^{(\ell)}}-\left(xP_{n+1}^{(\ell)},P_{n+1}^{(\ell)}\right)_{\omega^{(\ell)}}\\
    =&2\left(\frac{\partial }{\partial \zeta}P_{n+1}^{(\ell)},P_{n+1}^{(\ell)}\right)_{\omega^{(\ell)}}
    +\left(G(\zeta)-\zeta^{-1}-b_{n+1}^{(\ell)}\right)=0,\\
    \frac{\partial }{\partial \zeta}
    \left(P_{n+1}^{(\ell)},P_{n}^{(\ell)}\right)_{\omega^{(\ell)}}
    =&\left(\frac{\partial }{\partial \zeta}P_{n+1}^{(\ell)},P_{n}^{(\ell)}\right)_{\omega^{(\ell)}}
    +\left(P_{n+1}^{(\ell)},\frac{\partial }{\partial \zeta}P_{n}^{(\ell)}\right)_{\omega^{(\ell)}}\\
    &+\left(G(\zeta)-\zeta^{-1}\right)\left(P_{n+1}^{(\ell)},P_{n}^{(\ell)}\right)_{\omega^{(\ell)}}-\left(xP_{n}^{(\ell)},P_{n+1}^{(\ell)}\right)_{\omega^{(\ell)}}\\
    =&\left(\frac{\partial }{\partial \zeta}P_{n+1}^{(\ell)},P_{n}^{(\ell)}\right)_{\omega^{(\ell)}}
    -a_{n}^{(\ell)}=0,\\
    \frac{\partial }{\partial \zeta}
    \left(P_{n+1}^{(\ell)},P_{k}^{(\ell)}\right)_{\omega^{(\ell)}}
    =&\left(\frac{\partial }{\partial \zeta}P_{n+1}^{(\ell)},P_{k}^{(\ell)}\right)_{\omega^{(\ell)}}
    +\left(P_{n+1}^{(\ell)},\frac{\partial }{\partial \zeta}P_{k}^{(\ell)}\right)_{\omega^{(\ell)}}\\
    &+\left(G(\zeta)-\zeta^{-1}\right)\left(P_{n+1}^{(\ell)},P_{k}^{(\ell)}\right)_{\omega^{(\ell)}}-\left(xP_{k}^{(\ell)},P_{n+1}^{(\ell)}\right)_{\omega^{(\ell)}}\\
    =&\left(\frac{\partial }{\partial \zeta}P_{n+1}^{(\ell)},P_{k}^{(\ell)}\right)_{\omega^{(\ell)}}
    =0,\quad k\leq n-1.
\end{align*}
Thus one has
\begin{align*}
    \left(\frac{\partial }{\partial \zeta}P_{n+1}^{(\ell)},P_{n+1}^{(\ell)}\right)_{\omega^{(\ell)}}
    &=-\frac{1}{2}\left(G(\zeta)-\zeta^{-1}-b_{n+1}^{(\ell)}\right),\\
    \left(\frac{\partial }{\partial \zeta}P_{n+1}^{(\ell)},P_{n}^{(\ell)}\right)_{\omega^{(\ell)}}
    &= a_{n}^{(\ell)},\quad
    \left(\frac{\partial }{\partial \zeta}P_{n+1}^{(\ell)},P_{k}^{(\ell)}\right)_{\omega^{(\ell)}}
    =0,\quad k\leq n-1.
\end{align*}
Because  $\frac{\partial P_{n+1}^{(\ell)}}{\partial \zeta}$ is
a polynomial and its degree is not larger than $ n+1$,
using \eqref{eq:P01orth-3} gives  \eqref{eq:partialPnzeta}.
The proof is completed.
\qed\end{proof}

\subsection{Proof of Theorem \ref{thm:derivex}}
\label{p:thm:derivex}
\begin{proof}
Similar to the proof of Theorem \ref{thm:derivezeta}, one has
\[
{\frac{\partial }{\partial x}\omega^{(1)}(x;\zeta)
=x\omega^{(0)}(x;\zeta)-\zeta\omega^{(1)}(x;\zeta).}
\]
Because the degrees of polynomials
$\frac{\partial P_{n+1}^{(0)}}{\partial x}$ and
$(x^2-1)\frac{\partial P_{n}^{(1)}}{\partial x}+xP_{n}^{(1)}$
are not larger than  $ n$ and  $n+1$, respectively,
and
\[
\lim_{x\rightarrow+\infty}P_{i}^{(0)}(x;\zeta)P_{j}^{(1)}(x;\zeta)\omega^{(1)}(x;\zeta)=0,\quad \lim_{x\rightarrow1}P_{i}^{(0)}(x;\zeta)P_{j}^{(1)}(x;\zeta)\omega^{(1)}(x;\zeta)=0,\quad \forall i,j\in\mathbb{N},
\]
one {can} calculate the expansion coefficients in  \eqref{eq:P01orth-3} as follows
\begin{align*}
\left(\frac{\partial }{\partial x}P_{n+1}^{(0)},P_{n}^{(1)}\right)_{\omega^{(1)}}=&
\left((n+1)c_{n+1}^{(0)}x^{n},P_{n}^{(1)}\right)_{\omega^{(1)}}=
\frac{n+1}{\tilde{p}_{n}}\left(P_{n}^{(1)},P_{n}^{(1)}\right)_{\omega^{(1)}}
=\frac{n+1}{\tilde{p}_{n}},\\
\left(\frac{\partial }{\partial x}P_{n+1}^{(0)},P_{n-1}^{(1)}\right)_{\omega^{(1)}}=&
\int_{1}^{+\infty}\frac{\partial }{\partial x}\left(P_{n+1}^{(0)}P_{n-1}^{(1)}\omega^{(1)}\right)dx
-\left(P_{n+1}^{(0)},(x^2-1)\frac{\partial }{\partial x}P_{n-1}^{(1)}\right)_{\omega^{(0)}}\\
&-\left(P_{n+1}^{(0)},xP_{n-1}^{(1)}\right)_{\omega^{(0)}}+\zeta\left(P_{n+1}^{(0)},(x^2-1)P_{n-1}^{(1)}\right)_{\omega^{(0)}}
=\zeta r_{n},\\
\left(\frac{\partial }{\partial x}P_{n+1}^{(0)},P_{k}^{(1)}\right)_{\omega^{(1)}}=&
\int_{1}^{+\infty}\frac{\partial }{\partial x}\left(P_{n+1}^{(0)}P_{k}^{(1)}\omega^{(1)}\right)dx
-\left(P_{n+1}^{(0)},(x^2-1)\frac{\partial }{\partial x}P_{k}^{(1)}\right)_{\omega^{(0)}}\\
&-\left(P_{n+1}^{(0)},xP_{k}^{(1)}\right)_{\omega^{(0)}}
+\zeta\left(P_{n+1}^{(0)},(x^2-1)P_{k}^{(1)}\right)_{\omega^{(0)}}
=0,k\leq n-2.
\end{align*}
and
%
\begin{align*}
\left((x^2-1)\frac{\partial }{\partial x}P_{n}^{(1)}+xP_{n}^{(1)},P_{n+1}^{(0)}\right)_{\omega^{(0)}}=&
\left((n+1)c_{n}^{(1)}x^{n+1},P_{n+1}^{(0)}\right)_{\omega^{(0)}}\\
=&(n+1)\tilde{p}_{n}\left(P_{n+1}^{(0)},P_{n+1}^{(0)}\right)_{\omega^{(0)}}
=(n+1)\tilde{p}_{n},\\
\left((x^2-1)\frac{\partial }{\partial x}P_{n}^{(1)}+xP_{n}^{(1)},P_{n}^{(0)}\right)_{\omega^{(0)}}=&
\int_{1}^{+\infty}\frac{\partial }{\partial x}\left(P_{n}^{(1)}P_{n}^{(0)}\omega^{(1)}\right)dx-\left(P_{n}^{(1)},\frac{\partial }{\partial x}P_{n}^{(0)}\right)_{\omega^{(1)}}\\
&+\zeta\left(P_{n}^{(1)},P_{n}^{(0)}\right)_{\omega^{(1)}}=\zeta\left(P_{n}^{(1)},P_{n}^{(0)}\right)_{\omega^{(1)}}=\zeta p_{n},\\
\left((x^2-1)\frac{\partial }{\partial x}P_{n}^{(1)}+xP_{n}^{(1)},P_{k}^{(0)}\right)_{\omega^{(0)}}=&
\int_{1}^{+\infty}\frac{\partial }{\partial x}\left(P_{n}^{(1)}P_{k}^{(0)}\omega^{(1)}\right)dx-\left(P_{n}^{(1)},\frac{\partial }{\partial x}P_{k}^{(0)}\right)_{\omega^{(1)}}\\
&+\zeta\left(P_{n}^{(1)},P_{k}^{(0)}\right)_{\omega^{(1)}}=0,k\leq n-1.\\
\end{align*}
The proof is completed.
\qed\end{proof}
\subsection{Proof of Theorem \ref{thm:zerointer01}}
\label{p:thm:zerointer01}
\begin{proof}
Substituting $\{x_{i,n+1}^{(0)}\}_{i=1}^{n+1}$ into \eqref{eq:recP01x} gives
\[
   \left((x_{i,n+1}^{(0)})^2-1\right)P_{n}^{(1)}(x_{i,n+1}^{(0)};\zeta)=\tilde{r}_{n}P_{n}^{(0)}(x_{i,n+1}^{(0)};\zeta).
\]
which implies that $\tilde{r}_{n}\neq0$. In fact, if assuming   $\tilde{r}_{n}=0$,
then the above identity and  the fact that $(x_{i,n+1}^{(0)})^2-1>0$
imply $P_{n}^{(1)}(x_{i,n+1}^{(0)};\zeta)=0$,
which contradicts with $P_{n}^{(1)}$ being a polynomial of degree $n$.

Using Theorem \ref{thm:zerointerself} gives
\[
{\rm sign}\left(P_{n}^{(1)}(x_{i,n+1}^{(0)};\zeta)P_{n}^{(1)}(x_{i+1,n+1}^{(0)};\zeta)\right)={\rm sign}\left(P_{n}^{(0)}(x_{i,n+1}^{(0)};\zeta)P_{n}^{(0)}(x_{i+1,n+1}^{(0)};\zeta)\right)<0.
\]
Thus there exist no less than one {zero} of  the polynomial $P_{n}^{(1)}$ in each subinterval
$\left(x_{i,n+1}^{(0)},x_{i+1,n+1}^{(0)}\right)$.
The proof is completed.
\qed\end{proof}

\subsection{Proof of Corollary \ref{col:sign}}
\label{p:col:sign}
\begin{proof}
It is obvious that
\[
p_{n}=\frac{c_{n}^{(0)}}{c_{n}^{(1)}}>0,\quad r_{n}=\frac{c_{n-1}^{(1)}}{c_{n+1}^{(0)}}>0,\quad
\tilde{p}_{n}=\frac{c_{n}^{(1)}}{c_{n+1}^{(0)}}>0.
\]
Using  Theorems \ref{thm:rec} and \ref{thm:zerointer01} gives
{\begin{align*}
\tilde{q}_{n}&=\sum_{i=1}^{n+1}x_{i,n+1}^{(0)}-\sum_{i=1}^{n}x_{i,n}^{(1)}=\sum_{i=1}^{n}\left(x_{i+1,n+1}^{(0)}-x_{i,n}^{(1)}\right)
+x_{1,n+1}^{(0)}>0,\\
q_{n}&=\tilde{p}_{n}\left(\sum_{i=1}^{n+2}x_{i,n+2}^{(0)}-\sum_{i=1}^{n}x_{i,n}^{(1)}\right)=\tilde{p}_{n}\left(b_{n+1}^{(0)}+\tilde{q}_{n}\right)>0,
\end{align*}}
which imply $q_{n}>0$ and $\tilde{q}_{n}>0$.

Comparing the coefficients of the term of order $n$
at two sides of \eqref{eq:recP01x} gives
\begin{align*}
    &\tilde{r}_{n}=p_{n}^{-1}
\left(\sum_{i=1}^{n}\sum_{j=i+1}^{n}x_{i,n}^{(1)}x_{j,n}^{(1)}-1-
    \sum_{i=1}^{n+1}\sum_{j=i+1}^{n+1}x_{i,n+1}^{(0)}x_{j,n+1}^{(0)}+
    \left(\sum_{i=1}^{n+1}x_{i,n+1}^{(0)}-\sum_{i=1}^{n}x_{i,n}^{(1)}\right)\sum_{i=1}^{n+1}x_{i,n+1}^{(0)}\right)\\
    &= p_{n}^{-1}
\left(\sum_{i=1}^{n}\sum_{j=i+1}^{n}x_{i,n}^{(1)}x_{j,n}^{(1)}+\sum_{i=1}^{n+1}\sum_{j=i}^{n+1}x_{i,n+1}^{(0)}x_{j,n+1}^{(0)}
    -\sum_{i=1}^{n}x_{i,n}^{(1)}\sum_{i=1}^{n+1}x_{i,n+1}^{(0)}-1\right)\\
    &=p_{n}^{-1}
\left(\sum_{i=1}^{n}x_{i+1,n+1}^{(0)}(x_{i+1,n+1}^{(0)}-x_{i,n}^{(1)})+(x_{1,n+1}^{(0)})^{2}-1   +\sum_{i=0}^{n}\sum_{j=i+1}^{n}(x_{i+1,n+1}^{(0)}-x_{i,n}^{(1)})(x_{j+1,n+1}^{(0)}-x_{j,n}^{(1)})\right),
\end{align*}
where $x_{0,n}^{(1)}=0$.

{Combining Theorem \ref{thm:zerointer01} gives $\tilde{r}_{n}>0$.}
The proof is completed.
\qed\end{proof}
\subsection{Proof of Corollary \ref{col:varx}}
\label{p:col:varx}
\begin{proof}
Taking partial derivative of $P_{n}^{(\ell)}(x_{i,n}^{(\ell)};\zeta)$ with respect to $\zeta$
and using  Theorem \ref{thm:derivezeta}
gives
\[
    \frac{\partial x_{i,n}^{(\ell)}}{\partial \zeta}=-\left(\frac{\partial P_{n}^{(\ell)}}{\partial x}(x_{i,n}^{(\ell)};\zeta)\right)^{-1}
    \left(\frac{\partial P_{n}^{(\ell)}}{\partial \zeta}(x_{i,n}^{(\ell)};\zeta)\right)=-a_{n-1}^{(\ell)}\left(\frac{\partial P_{n}^{(\ell)}}{\partial x}(x_{i,n}^{(\ell)};\zeta)\right)^{-1}P_{n-1}^{(\ell)}(x_{i,n}^{(\ell)};\zeta).
\]
Due to Theorem \ref{thm:zerointerself}, one has
\[
{\rm sign}(P_{n-1}^{(\ell)}(x_{i,n}^{(\ell)};\zeta))=(-1)^{n+i}={\rm sign}\left(\frac{\partial P_{n}^{(\ell)}}{\partial x}(x_{i,n}^{(\ell)};\zeta)\right).
\]
Combining them  completes the proof.
\qed\end{proof}

\subsection{Proof of Lemma \ref{lem:zerointer}}
\label{p:lem:zerointer}
\begin{proof}
According to the definition of $Q_{2n}(x;\zeta)$ in \eqref{eq:Q2n},
it is not difficult to know that $Q_{2n}(x;\zeta)$ is an even function and
a polynomial of degree $2n$.

If taking $x$ in \eqref{eq:Q2n} as the zero of $P_{n+1}^{(0)}(x;\zeta)$, i.e.
 $x=x_{i,n+1}^{(0)}$,  $i=1,\cdots,{n+1}$,  then one has
\[
Q_{2n}(x_{i,n+1}^{(0)};\zeta)=P_{n+1}^{(0)}(-x_{i,n+1}^{(0)};\zeta)P_{n}^{(1)}(x_{i,n+1}^{(0)};\zeta).
\]
Since
\[
{\rm sign}\left(P_{n+1}^{(0)}(-x_{i,n+1}^{(0)};\zeta)\right)=(-1)^{n+1}, \quad i=1,\cdots, n+1,
\]
using Theorem \ref{thm:zerointer01} gives
\[
{\rm sign}\left(Q_{2n}(x_{i,n+1}^{(0)};\zeta)Q_{2n}(x_{i+1,n+1}^{(0)};\zeta)\right)
={\rm sign}\left(P_{n}^{(1)}(x_{i,n+1}^{(0)};\zeta)P_{n}^{(1)}(x_{i+1,n+1}^{(0)};\zeta)\right)<0,
\]
for $i=1,\cdots,n$, which {implies} that there exist no less than one {zero} of $Q_{2n}(x;\zeta)$
in each subinterval
 $ (x_{i,n+1}^{(0)},x_{i+1,n+1}^{(0)})$, $i=1,\cdots,n$.
Because $Q_{2n}(x;\zeta)$ is an even polynomial of degree $2n$, there exists exactly one zero of $Q_{2n}(x;\zeta)$ in each subinterval
 $(x_{i,n+1}^{(0)},x_{i+1,n+1}^{(0)})$, $i=1,\cdots,n$.
 The proof is completed.
\qed\end{proof}

\subsection{Proof of Lemma \ref{lem:derivezetax}}
\label{p:lem:derivezetax}
\begin{proof}According to the definition of $Q_{2n}(x;\zeta)$ in \eqref{eq:Q2n}, one has
\begin{align*}
    \frac{\partial Q_{2n}}{\partial \zeta}(z_{i,n};\zeta)=&
    \frac{\partial P_{n+1}^{(0)}}{\partial \zeta}\Big|_{x=z_{i,n}}
     P_{n}^{(1)}(-z_{i,n};\zeta)
    +\frac{\partial P_{n}^{(1)}}{\partial \zeta}\Big|_{x=z_{i,n}}
    P_{n+1}^{(0)}(-z_{i,n};\zeta)\\
    &+\frac{\partial P_{n+1}^{(0)}}{\partial \zeta} \Big|_{x=-z_{i,n}} P_{n}^{(1)}(z_{i,n};\zeta)
    +\frac{\partial P_{n}^{(1)}}{\partial \zeta}\Big|_{x=-z_{i,n}} P_{n+1}^{(0)}(z_{i,n};\zeta).
\end{align*}
Using Theorem \ref{thm:derivezeta} gives
    \begin{align*}
    \frac{\partial Q_{2n}}{\partial \zeta}(z_{i,n};\zeta)=&
    a_{n}^{(0)}\left(P_{n}^{(0)}(z_{i,n};\zeta)P_{n}^{(1)}(-z_{i,n};\zeta)+P_{n}^{(0)}(-z_{i,n};\zeta)P_{n}^{(1)}(z_{i,n};\zeta)\right)\\
    &+a_{n-1}^{(1)}\left(P_{n-1}^{(1)}(z_{i,n};\zeta)P_{n+1}^{(0)}(-z_{i,n};\zeta)+P_{n-1}^{(1)}(-z_{i,n};\zeta)P_{n+1}^{(0)}(z_{i,n};\zeta)\right)\\
    &+\left(\frac{1}{2}(b_{n+1}^{(0)}+b_{n}^{(1)})-(G(\zeta)-\zeta^{-1})\right)Q_{2n}(z_{i,n};\zeta)\\
    =&a_{n}^{(0)}\left(P_{n}^{(0)}(z_{i,n};\zeta)P_{n}^{(1)}(-z_{i,n};\zeta)+P_{n}^{(0)}(-z_{i,n};\zeta)P_{n}^{(1)}(z_{i,n};\zeta)\right)\\
    &+a_{n-1}^{(1)}\left(P_{n-1}^{(1)}(z_{i,n};\zeta)P_{n+1}^{(0)}(-z_{i,n};\zeta)+P_{n-1}^{(1)}(-z_{i,n};\zeta)P_{n+1}^{(0)}(z_{i,n};\zeta)\right).
    \end{align*}
Substituting \eqref{eq:recP01x} and \eqref{eq:recP10x} into it
gives
 \begin{align*}
    \frac{\partial Q_{2n}}{\partial \zeta}(z_{i,n};\zeta)
    =&\frac{2a_{n}^{(0)}}{\tilde{r}_{n}}\left(((z_{i,n})^2-1)P_{n}^{(1)}(z_{i,n};\zeta)P_{n}^{(1)}(-z_{i,n};\zeta)+
    \tilde{p}_{n}z_{i,n}P_{n+1}^{(0)}(-z_{i,n};\zeta)P_{n}^{(1)}(z_{i,n};\zeta)\right)\\
    &-\frac{2a_{n}^{(0)}}{\tilde{r}_{n}}\left(P_{n+1}^{(0)}(z_{i,n};\zeta)P_{n+1}^{(0)}(-z_{i,n};\zeta)-
    \tilde{p}_{n}^{-1}z_{i,n}P_{n+1}^{(0)}(-z_{i,n};\zeta)P_{n}^{(1)}(z_{i,n};\zeta)\right)\\
    =&2\frac{P_{n}^{(1)}(z_{i,n};\zeta)}{P_{n}^{(1)}(-z_{i,n};\zeta)}
    \frac{a_{n}^{(0)}}{\tilde{r}_{n}}\left((\tilde{p}_{n}+
    \tilde{p}_{n}^{-1})z_{i,n}P_{n}^{(1)}(-z_{i,n};\zeta)P_{n+1}^{(0)}(-z_{i,n};\zeta)\right.\\
    &\left.+(z_{i,n}^2-1)P_{n}^{(1)}(-z_{i,n};\zeta)^2
    +P_{n+1}^{(0)}
    (-z_{i,n};\zeta)^2\right)\\
    =&2\frac{P_{n}^{(1)}(z_{i,n};\zeta)}{P_{n}^{(1)}(-z_{i,n};\zeta)}
    \frac{a_{n}^{(0)}}{\tilde{r}_{n}}\left((c_{n+1}^{(0)})^2\prod_{j=1}^{n+1}(z_{i,n}+x_{j,n+1}^{(0)})
    \tilde{I}_{1}-(c_{n}^{(1)})^2\prod_{j=1}^{n}(z_{i,n}+x_{j,n}^{(1)})
    \tilde{I}_{2}\right),
 \end{align*}
 where
  \[
    \tilde{I}_{1}:=\prod_{j=1}^{n+1}(z_{i,n}+x_{j,n+1}^{(0)})-z_{i,n}\prod_{j=1}^{n}(z_{i,n}+x_{j,n}^{(1)}),
    \]
    \[
    \tilde{I}_{2}:=z_{i,n}\prod_{j=1}^{n+1}(z_{i,n}+x_{j,n+1}^{(0)})-(z_{i,n}^2-1)\prod_{j=1}^{n}(z_{i,n}+x_{j,n}^{(1)}).
    \]

Similarly, using Theorem \ref{thm:derivex} and \eqref{eq:recP01x}-\eqref{eq:recP10x} gives
      \begin{align*}
    \frac{\partial Q_{2n}}{\partial x}(z_{i,n};\zeta)
     =&\zeta\frac{c_{n-1}^{(1)}}{c_{n+1}^{(0)}}\left(P_{n-1}^{(1)}(z_{i,n};\zeta)P_{n}^{(1)}(-z_{i,n};\zeta)-P_{n-1}^{(1)}(-z_{i,n};\zeta)P_{n}^{(1)}(z_{i,n};\zeta)\right)\\
     &+(z_{i,n}^2-1)^{-1}\zeta\frac{c_{n}^{(0)}}{c_{n}^{(1)}}\left(P_{n}^{(0)}(z_{i,n};\zeta)P_{n+1}^{(0)}(-z_{i,n};\zeta)-P_{n}^{(0)}(-z_{i,n};\zeta)P_{n+1}^{(0)}(z_{i,n};\zeta)\right)\\
    =&2\zeta\frac{P_{n}^{(1)}(z_{i,n};\zeta)}{P_{n}^{(1)}(-z_{i,n};\zeta)}
    \frac{a_{n}^{(0)}}{\tilde{r}_{n}}\left(z_{i,n}P_{n}^{(1)}(-z_{i,n};\zeta)^{2}\right.\\
    &\left.+(\tilde{p}_{n}+\tilde{p}_{n}^{-1})P_{n+1}^{(0)}(-z_{i,n};\zeta)P_{n}^{(1)}(-z_{i,n};\zeta)
    +(z_{i,n}^2-1)^{-1}z_{i,n}P_{n+1}^{(0)}(-z_{i,n};\zeta)^2 \right)\\
    =&2\zeta\frac{P_{n}^{(1)}(z_{i,n};\zeta)}{P_{n}^{(1)}(-z_{i,n};\zeta)}
    \frac{a_{n}^{(0)}}{\tilde{r}_{n}}\left((c_{n+1}^{(0)})^2\frac{\prod_{j=1}^{n+1}(z_{i,n}+x_{j,n+1}^{(0)})}{z_{i,n}^2-1}
    \tilde{I}_{2}-(c_{n}^{(1)})^2\prod_{j=1}^{n}(z_{i,n}+x_{j,n}^{(1)})\tilde{I}_{1}\right).
   \end{align*}

Using  
Theorem  \ref{thm:zerointer01} gives
\begin{align*}
    z_{i,n}+x_{j+1,n+1}^{(0)}>z_{i,n}+x_{j,n}^{(1)}
    ,\ j=1,\cdots,n,
\\
z_{i,n}+x_{1,n+1}^{(0)}>z_{i,n}+1>z_{i,n}>z_{i,n}-1, 
\end{align*}
for $i=1,\cdots,n$, which imply
    \[
    \tilde{I}_{1}>0,\quad \tilde{I}_{2}>0.
    \]
Thus one has
     \begin{align*}
    \prod_{j=1}^{n+1}&(z_{i,n}+x_{j,n+1}^{(0)})\tilde{I}_{1}-\prod_{j=1}^{n}(z_{i,n}+x_{j,n}^{(1)})
    \tilde{I}_{2}\\
    =&\left(\prod_{j=1}^{n+1}(z_{i,n}+x_{j,n+1}^{(0)})-(z_{i,n}-1)\prod_{j=1}^{n}(z_{i,n}+x_{j,n}^{(1)}
    )\right)\\
    &\cdot\left(\prod_{j=1}^{n+1}(z_{i,n}+x_{j,n+1}^{(0)})-(z_{i,n}+1)\prod_{j=1}^{n}(z_{i,n}+x_{j,n}^{(1)})\right)>0,
    \end{align*}
    \begin{align*}
    \prod_{j=1}^{n+1}&(z_{i,n}+x_{j,n+1}^{(0)})
    \tilde{I}_{2}-((z_{i,n})^2-1)\prod_{j=1}^{n}(z_{i,n}
    +x_{j,n}^{(1)})\tilde{I}_{1}\\
    =&z_{i,n}{\left(\prod_{j=1}^{n+1}(z_{i,n}+x_{j,n+1}^{(0)})\tilde{I}_{1}-\prod_{j=1}^{n}(z_{i,n}+x_{j,n}^{(1)})
    \tilde{I}_{2}\right)}\\
    &+2\prod_{j=1}^{n+1}(z_{i,n}+x_{j,n+1}^{(0)})\prod_{j=1}^{n}(z_{i,n}+x_{j,n}^{(1)})>0.
    \end{align*}
Using Corollaries \ref{col:sign} and \ref{corollary3.2},
and the above results gives
\eqref{EQ3.29zzzz}.
The proof is completed.
\qed\end{proof}
\subsection{Proof of Lemma \ref{lem:varz}}
\label{p:lem:varz}
\begin{proof}
Taking partial derivative of $Q_{2n}(z_{i,n};\zeta)$ with respect to $\zeta$ gives
\[
    \frac{\partial z_{i,n}}{\partial \zeta}=-\left(\frac{\partial Q_{2n}}{\partial x}\Big|_{x=z_{i,n}}\right)^{-1}
    \left(\frac{\partial Q_{2n}}{\partial \zeta}\Big|_{x=z_{i,n}}\right), i=1,\cdots,n.
\]
Using Lemma \ref{lem:derivezetax} {completes} the proof.
\qed\end{proof}

\subsection{Proof of Theorem \ref{thm:eig}}
\label{p:thm:eig}
\begin{proof}
Obviously, both vectors $\vec{u}_{i,n}$ and $\vec{v}_{i,n}$ defined in
\eqref{EQ:3.29bbbb} are not zero  at the same time, $i=\pm 1,\cdots,\pm n$.
The nonzero eigenvalues and eigenvectors of the matrix pair $\vec{A}^{0}_{n}$ and $\vec{A}^{1}_{n}$
in \eqref{eq:eigvalue} and \eqref{eq:eigvec}
{can} be obtained with the aid of  \eqref{eq:eqm1final}-\eqref{eq:eqm2final} and
 Lemma \ref{lem:zerointer}.
Using  Lemma \ref{lem:varz} further gives \eqref{eq:eigparzeta}.

In the following, let us discuss the eigenvector $\vec{y}_{0,n}$.
Multiplying \eqref{eq:recQP} by $P_{n+1}^{(0)}(-x;\zeta)$ gives
\begin{equation}
  (x^2-1)\vec{P}_{n-1}^{(1)}(x;\zeta)P_{n+1}^{(0)}(-x;\zeta)=\vec{J}_{n-1}\vec{P}_{n}^{(0)}(x;\zeta)P_{n+1}^{(0)}(-x;\zeta)
  +r_{n}P_{n+1}^{(0)}(x;\zeta)P_{n+1}^{(0)}(-x;\zeta)\vec{e}_{n}.\label{eq:eig0}
\end{equation}
Transforming \eqref{eq:eig0} by $x$ to $-x$
and then subtracting it from \eqref{eq:eig0} and
 letting $x=1$  gives
as follows
\begin{equation*}
    0=\vec{J}_{n-1}\left(\vec{P}_{n}^{(0)}(1;\zeta)P_{n+1}^{0}(-1;\zeta)-\vec{P}_{n}^{(0)}(-1;\zeta)P_{n+1}^{(0)}(1;\zeta)\right)=\vec{J}_{n-1}\vec{u}_{0,n},
\end{equation*}
which is a special case of \eqref{eq:eigeq1} with $ \hat{\lambda} =0$.
 The proof is completed.
\qed\end{proof}

\section{Proofs in Section \ref{sec:moment}}
\label{sec:App4}

\subsection{Proof of Lemma \ref{lem:basis}}
\label{p:lem:basis}
\begin{proof}
(i)  Due to the definition of $E$ and $p_{<1>}$, each component of $\vec{\mathcal{P}}_{\infty}[u,\theta]$
(resp.  $\vec{\mathcal{P}}_{M}[u,\theta]$)
belongs to $\mathbb{H}^{g^{(0)}_{[u,\theta]}}$  (resp. $\mathbb{H}^{g^{(0)}_{[u,\theta]}}_{M}$) obviously.

(ii) The mathematical induction is  used to prove
that any element in the space $\mathbb{H}^{g^{(0)}_{[u,\theta]}}$  (resp. $\mathbb{H}^{g^{(0)}_{[u,\theta]}}_{M}$)
{can} be written into a linear combination  of
vectors in $\vec{\mathcal{P}}_{\infty}[u,\theta]$
(resp.  $\vec{\mathcal{P}}_{M}[u,\theta]$) .
For $M=1$, it is clear to have the linear combination
\begin{align*}
p^{\alpha}g^{(0)}_{[u,\theta]}\overset{\eqref{eq:divp}}{=}&\left(p^{<\alpha>}+U^{\alpha}E\right)g^{(0)}_{[u,\theta]}
\overset{\eqref{eq:divpp}}{=}
\left(-(U^{0})^{-1}U^{1-\alpha}p_{<1>}+U^{\alpha}E\right)g^{(0)}_{[u,\theta]}\\
\overset{\eqref{eq:Preprezero}}{=}&-(c_{0}^{(1)})^{-1}U^{1-\alpha}\tilde{P}_{0}^{(1)}[u,\theta]+
(c_{1}^{(0)})^{-1}U^{\alpha}\tilde{P}_{1}^{(0)}[u,\theta]+(c_{0}^{(0)})^{-1}U^{\alpha}x_{1,1}^{(0)}\tilde{P}_{0}^{(0)}[u,\theta],
\end{align*}
where  the decomposition of the particle velocity vector  \eqref{eq:divp} has been used.

Assume that the  linear combination
\begin{align*}
p^{\mu_{1}}p^{\mu_{2}}\cdots p^{\mu_{M}}g^{(0)}_{[u,\theta]}=&\sum_{i=0}^{M}c_{i,0}^{\mu_{1},\cdots,\mu_{M}}\tilde{P}_{i}^{(0)}[u,\theta]
+\sum_{i=0}^{M-1}c_{i,1}^{\mu_{1},\cdots,\mu_{M}}\tilde{P}_{i}^{(1)}[u,\theta],\\
&\mu_{i}=0,1, i\in\mathbb{N}, i\leq M, \quad c_{i,0}^{\mu_{1},\cdots,\mu_{M}},c_{i,1}^{\mu_{1},\cdots,\mu_{M}}\in \mathbb{R},
\end{align*}
 holds.  One has to show that $p^{\mu_{1}}p^{\mu_{2}}\cdots p^{\mu_{M+1}}g^{(0)}_{[u,\theta]}$
 {can} be written into a linear combination of components of $\vec{\mathcal{P}}_{M+1}[u,\theta]$.
 Because
\begin{align*}
&p^{\mu_{1}}p^{\mu_{2}}\cdots p^{\mu_{M+1}}g^{(0)}_{[u,\theta]}\\
&=\left(\sum_{i=0}^{M}c_{i,0}^{\mu_{1},\cdots,\mu_{M}}\tilde{P}_{i}^{(0)}[u,\theta]
+\sum_{i=0}^{M-1}c_{i,1}^{\mu_{1},\cdots,\mu_{M}}\tilde{P}_{i}^{(1)}[u,\theta]\right)
\left(-(U^{0})^{-1}U^{1-\mu_{M+1}}p_{<1>}+U^{\mu_{M+1}}E\right)\\
&=\sum_{i=0}^{M}c_{i,0}^{\mu_{1},\cdots,\mu_{M}}U^{\mu_{M+1}}E\tilde{P}_{i}^{(0)}[u,\theta]
-\sum_{i=0}^{M-1}c_{i,1}^{\mu_{1},\cdots,\mu_{M}}U^{1-\mu_{M+1}}(E^2-1)P_{i}^{(1)}(E;\zeta)\\
&-\sum_{i=0}^{M}c_{i,0}^{\mu_{1},\cdots,\mu_{M}}(U^{0})^{-1}U^{1-\mu_{M+1}}P_{i}^{(0)}(E;\zeta)p_{<1>}+
\sum_{i=0}^{M-1}c_{i,1}^{\mu_{1},\cdots,\mu_{M}}U^{\mu_{M+1}}\tilde{P}_{i}^{(1)}[u,\theta],
\end{align*}
one has
\begin{align*}
p^{\mu_{1}}p^{\mu_{2}}\cdots p^{\mu_{M+1}}g^{(0)}_{[u,\theta]}=&\sum_{i=0}^{M}c_{i,0}^{\mu_{1},\cdots,\mu_{M}}U^{\mu_{M+1}}
\left(a_{i-1}^{(0)}\tilde{P}_{i-1}^{(0)}[u,\theta]+b_{i}^{(0)}\tilde{P}_{i}^{(0)}[u,\theta]+a_{i}^{(0)}\tilde{P}_{i+1}^{(0)}[u,\theta]\right)\\
&-\sum_{i=0}^{M-1}c_{i,1}^{\mu_{1},\cdots,\mu_{M}}U^{1-\mu_{M+1}}
\left(p_{i}\tilde{P}_{i}^{(0)}[u,\theta]+q_{i}\tilde{P}_{i+1}^{(0)}[u,\theta]+r_{i+1}\tilde{P}_{i+2}^{(0)}[u,\theta]\right)\\
&-\sum_{i=0}^{M}c_{i,0}^{\mu_{1},\cdots,\mu_{M}}U^{1-\mu_{M+1}}
\left(r_{i-1}\tilde{P}_{i-2}^{(1)}[u,\theta]+q_{i-1}\tilde{P}_{i-1}^{(1)}[u,\theta]+p_{i}\tilde{P}_{i}^{(1)}[u,\theta]\right)\\
&+\sum_{i=0}^{M-1}c_{i,1}^{\mu_{1},\cdots,\mu_{M}}U^{\mu_{M+1}}
\left(a_{i-1}^{(1)}\tilde{P}_{i-1}^{(1)}[u,\theta]+b_{i}^{(1)}\tilde{P}_{i+1}^{(1)}[u,\theta]+a_{i}^{(1)}\tilde{P}_{i+1}^{(1)}[u,\theta]\right)\\
=:&\sum_{i=0}^{M+1}c_{i,0}^{\mu_{1},\cdots,\mu_{M+1}}\tilde{P}_{i}^{(0)}[u,\theta]
+\sum_{i=0}^{M}c_{i,1}^{\mu_{1},\cdots,\mu_{M+1}}\tilde{P}_{i}^{(1)}[u,\theta].
\end{align*}
 by using    the three-term recurrence relations \eqref{eq:recP0P1},  \eqref{eq:recP01}, and \eqref{eq:recP10} for
 the orthogonal polynomials $\{P_{n}^{(\ell)}(x;\zeta), \ell=0,1\}$.

(iii)  Using  \eqref{eq:P01orth} gives
\begin{equation}
\label{eq:orth}
<\tilde{P}_{i}^{(\ell)}[u,\theta],\tilde{P}_{j}^{(\ell)}[u,\theta]>_{g^{(0)}_{[u,\theta]}}=
\left(P_{i}^{(\ell)},P_{j}^{(\ell)}\right)_{ \omega^{(\ell)}}=\delta_{i,j},\ \ell=0,1.
\end{equation}
Because of  \eqref{eq:divp}, one has
\[
\frac{dp}{p^{0}}=dp_{<1>}\frac{-1+u(U^{0}E)^{-1}p_{<1>}}{-up_{<1>}+U^{0}E}=-\frac{dp_{<1>}}{U^{0}E}, \quad E=\sqrt{\left((U^{0})^{-1}p_{<1>}\right)^2+1}.
\]
Thus one obtains
\begin{align}
\nonumber <\tilde{P}_{i}^{(0)}[u,\theta],\tilde{P}_{j}^{(1)}[u,\theta]>_{g^{(0)}_{[u,\theta]}}=&
\int_{\mathbb{R}}g^{(0)}_{[u,\theta]}P_{i}^{(0)}(E;\zeta)P_{j}^{(1)}(E;\zeta)(U^{0})^{-1}p_{<1>}\frac{dp}{p^{0}}\\
=&-\int_{\mathbb{R}}g^{(0)}_{[u,\theta]}P_{i}^{(0)}(E;\zeta)P_{j}^{(1)}(E;\zeta)(U^{0})^{-1}p_{<1>}\frac{dp_{<1>}}{U^{0}E}=0.\label{eq:orthhuxiang}
\end{align}

Combining (i) and (ii) with (iii) completes
the proof.
\qed\end{proof}

\subsection{Proof of Lemma \ref{lem:derive}}
\label{p:lem:derive}
\begin{proof}
	For $s=t$ and $x$,
it is clear to have
\begin{equation*}
\frac{\partial E}{\partial s}=\frac{\partial u}{\partial s}\frac{1}{(1-u^2)}(U^{0})^{-1}p_{<1>},\quad \frac{\partial \left((U^{0})^{-1}p_{<1>}\right)}{\partial s}=\frac{\partial u}{\partial s}\frac{1}{1-u^2}E.
\end{equation*}
Using  those above identities and  \eqref{eq:equm1-zzzzz}  gives
\[
\frac{\partial g^{(0)}_{[u,\theta]}}{\partial s}=-\left(\frac{\partial \theta}{\partial s}\zeta^2\left(G(\zeta)-\zeta^{-1}-E\right)+
\frac{\partial u}{\partial s}\frac{1}{\theta(1-u^2)}(U^{0})^{-1}p_{<1>}\right) g^{(0)}_{[u,\theta]}.
\]
The derivation rule of compound function gives
\begin{align*}
  \frac{\partial \tilde{P}_{n}^{(0)}[u,\theta]}{\partial s}=&\frac{\partial P_{n}^{(0)}}{\partial E}\frac{\partial E}{\partial s}g^{(0)}_{[u,\theta]}
  -\zeta^2\frac{\partial P_{n}^{(0)}}{\partial \zeta}\frac{\partial \theta}{\partial s}g^{(0)}_{[u,\theta]}
  +P_{n}^{(0)}\frac{\partial g^{(0)}_{[u,\theta]}}{\partial s},\\
  \frac{\partial \tilde{P}_{n-1}^{(1)}[u,\theta]}{\partial s}=&\frac{\partial P_{n-1}^{(1)}}{\partial E}\frac{\partial E}{\partial s}(U^{0})^{-1}p_{<1>}g^{(0)}_{[u,\theta]}
  -\zeta^2\frac{\partial P_{n-1}^{(1)}}{\partial \zeta}\frac{\partial \theta}{\partial s}(U^{0})^{-1}p_{<1>}g^{(0)}_{[u,\theta]}\\
  &+P_{n-1}^{(1)}\frac{\partial \left((U^{0})^{-1}p_{<1>}\right)}{\partial s}g^{(0)}_{[u,\theta]}
  +P_{n-1}^{(1)}(U^{0})^{-1}p_{<1>}\frac{\partial g^{(0)}_{[u,\theta]}}{\partial s},
\end{align*}
Combining them and using Theorems \ref{thm:rec}-\ref{thm:derivex} 
{complete} the proof.
\qed\end{proof}

\subsection{Proof of Lemma \ref{lem:rec}}
\label{p:lem:rec}
\begin{proof}
Using the three-term recurrence relations
 \eqref{eq:recP0P1mat},
 \eqref{eq:recPQ},  and \eqref{eq:recQP} gives
\begin{align*}
    E\vec{\mathcal{\tilde{P}}}_{M}[u,\theta] &=\vec{A}^{0}_{M} \vec{\mathcal{\tilde{P}}}_{M}[u,\theta]+a_{M}^{(0)}\tilde{P}_{M+1}^{(0)}[u,\theta]\vec{e}_{2M+1}^{3}
    +a_{M-1}^{(1)}\tilde{P}_{M}^{(1)}[u,\theta]\vec{e}_{2M+1}^{2}, \\
    (U^{0})^{-1}p_{<1>}\vec{\mathcal{\tilde{P}}}_{M}&=\vec{A}^{1}_{M}\vec{\mathcal{\tilde{P}}}_{M}[u,\theta] +p_{M}\tilde{P}_{M}^{(1)}[u,\theta]\vec{e}_{2M+1}^{3}
    +r_{M}\tilde{P}_{M+1}^{(0)}[u,\theta]\vec{e}_{2M+1}^{2},
 \end{align*}
 where $\vec{e}_{2M+1}^{3}$ is the $(M+1)th$ column of the identity matrix of order $(2M+1)$.
Thus one has
\begin{align*}
E\vec{\mathcal{P}}_{M}[u,\theta] &= \vec{P}_{M}^{p}\vec{A}^{0}_{M}(\vec{P}_{M}^{p})^{T} \vec{\mathcal{P}}_{M}[u,\theta]+a_{M}^{(0)}\tilde{P}_{M+1}^{(0)}[u,\theta]\vec{e}_{2M+1}^{1}
    +a_{M-1}^{(1)}\tilde{P}_{M}^{(1)}[u,\theta]\vec{e}_{2M+1}^{2},\\
(U^{0})^{-1}p_{<1>}\vec{\mathcal{P}}_{M}&= \vec{P}_{M}^{p}\vec{A}^{1}_{M}(\vec{P}_{M}^{p})^{T}\vec{\mathcal{P}}_{M}[u,\theta] +p_{M}\tilde{P}_{M}^{(1)}[u,\theta]\vec{e}_{2M+1}^{1}
    +r_{M}\tilde{P}_{M+1}^{(0)}[u,\theta]\vec{e}_{2M+1}^{2}.
 \end{align*}
{Combining} them with \eqref{eq:divp} completes the proof.
\qed\end{proof}

\subsection{Proof of Lemma \ref{lem:project}}
\label{p:lem:project}
\begin{proof}
It is obvious that $\Pi_{M}[u,\theta]$ is a linear bounded operator and
 $\Pi_{M}[u,\theta]f\in\mathbb{H}^{g^{(0)}_{[u,\theta]}}_{M}$ for all $f\in\mathbb{H}^{g^{(0)}_{[u,\theta]}}$.

For all $f\in\mathbb{H}^{g^{(0)}_{[u,\theta]}}_{M}$, 
besides \eqref{EQ-projection-aaaa},  one has by using Lemma \ref{lem:basis}
\begin{align*}
f=&\sum_{i=0}^{M}\tilde{f}_{i}^{0}\tilde{P}_{i}^{(0)}[u,\theta]+
  \sum_{j=0}^{M-1}\tilde{f}_{j}^{1}\tilde{P}_{j}^{(1)}[u,\theta].
\end{align*}
Taking respectively the inner product with $\tilde{P}_{i}^{(0)}[u,\theta]$ and $\tilde{P}_{j}^{(1)}[u,\theta]$ from both sides
of the last equation gives
\begin{align*}
    f_{i}^{0}=<f,\tilde{P}_{i}^{(0)}[u,\theta]>_{g^{(0)}_{[u,\theta]}},i\leq M,\quad
    f_{j}^{1}=<f,\tilde{P}_{j}^{(1)}[u,\theta]>_{g^{(0)}_{[u,\theta]}},j\leq M-1.
\end{align*}
Comparing them with the coefficients in \eqref{eq:deff1} shows that
  $\tilde{f}_{i}^{0}=f_{i}^{0}$, $\tilde{f}_{j}^{1}=f_{j}^{1}$, \\$i=0,\cdots, M$, $j=1,\cdots,M-1$.
The proof is completed.
\qed\end{proof}

\section{Proofs in Section \ref{sec:prop}}
\label{sec:App5}
\subsection{Proof of Lemma \ref{lem:Dinv}}
\label{p:lem:Dinv}
 \begin{proof}
 It is obvious that for $M=1$, the matrix $D_{M}$ is invertible
 because \\$\det(\vec{D}_M)=\rho\zeta^2c_0^{(1)}(c_0^{(0)}c_1^{(0)}(1-u^2))^{-1}>0$.
{For $M\geq2$,} according to the form of $\vec{D}_M$ in Section \ref{subsec:deduction},
 one has
\[
\det(\vec{D}_{M})=\det(\vec{D}_{2})=\zeta^3c_{2}^{(0)}c_{1}^{(1)}(x_{1,2}^{(0)}+x_{2,2}^{(0)})(\rho G(\zeta)+\Pi)\rho(c_{1}^{(0)}c_{0}^{(0)}(1-u^2))^{-1}.
\]
Using {$\Pi>-\rho\theta$} gives
\[
\det(\vec{D}_{M})>\zeta^3c_{2}^{(0)}c_{1}^{(1)}(x_{1,2}^{(0)}+x_{2,2}^{(0)})\rho^2( G(\zeta)-\zeta^{-1})(c_{1}^{(0)}c_{0}^{(0)}(1-u^2))^{-1}>0.
\]
The proof is completed.
\qed\end{proof}

\subsection{Proof of Theorem \ref{thm:eigmoment}}
\label{p:thm:eigmoment}
 \begin{proof}
 	Consider the following generalized eigenvalue problem (2nd sense): Find a vector $\vec r$
 	that obeys $\lambda \vec{B}_{M}^{0}\vec{r}=\vec{B}_{M}^{1}\vec{r}$ or
 	$ \lambda \vec{M}_{M}^{t}\vec{D}_{M}\vec{r}=\vec{M}_{M}^{x}D_{M}\vec{r}$. Thanks to \eqref{eq:MtMx},
 	this eigenvalue problem is equivalent to
 \[
  (\lambda-u)    \vec{A}_{M}^{0}(\vec{P}_{M}^{p})^{T}\vec{D}_{M}\vec{r}= (\lambda u-1) \vec{A}_{M}^{1}(\vec{P}_{M}^{p})^{T}\vec{D}_{M}\vec{r}.
 \]
Because Theorem \ref{thm:eig}  tells us that
 $ \hat{\lambda} _{i,M}$ and $\vec{y}_{i,M}$ satisfy
\begin{equation*}
     \hat{\lambda} _{i,M}\vec{A}_{M}^{0}\vec{y}_{i,M}=\vec{A}_{M}^{1}\vec{y}_{i,M},\ | \hat{\lambda} _{i,M}|<1,
\end{equation*}
the scalar  $\lambda_{i,M}$  in \eqref{eq:lambda}  and vector $\vec{r}_{i,M}$ in \eqref{eq:r}
solve  the above generalized eigenvalue problem,
%
and satisfy
\[
|\lambda_{i,M}|<\frac{1-u}{1-u}=1.
\]
The proof is completed.
\qed\end{proof}

\subsection{Proof of Lemma \ref{lem:Mpositive}}
\label{p:lem:Mpositive}
 \begin{proof}
Because $U^{0}\vec{M}_{M}^{t}-U^{1}\vec{M}_{M}^{x}=\vec{P}_{M}^{p}A_{M}^{0}(\vec{P}_{M}^{p})^{T}$ and
 the permutation matrix $\vec{P}_{M}^{p}$ in  \eqref{eq:MtMx-22222} satisfies $\vec{P}_{M}^{p}(\vec{P}_{M}^{p})^{T}=(\vec{P}_{M}^{p})^{T}\vec{P}_{M}^{p}=I$,
 two matrices  $U^{0}\vec{M}_{M}^{t}-U^{1}\vec{M}_{M}^{x}$ and $\vec{A}_{M}^{0}$ are  similar
 and thus have the same eigenvalues.
The definition of  $\vec{A}_{M}^{0}$  in \eqref{EQ:3.20-000000}
tells us that the eigenvalues of  $\vec{A}_{M}^{0}$ are the zeros of $P_{M+1}^{(0)}(x;\zeta)$ and $P_{M}^{(1)}(x;\zeta)$
which are larger than one \cite[Theorem 3.4]{SP:2011},
so the matrix $U^{0}\vec{M}_{M}^{t}-U^{1}\vec{M}_{M}^{x}$ is positive definite.

 Theorem \ref{thm:eig} implies
\[
    \rho\left((\vec{A}_{M}^{0})^{-\frac{1}{2}}\vec{A}_{M}^{1}(\vec{A}_{M}^{0})^{-\frac{1}{2}}\right)=\rho\left((\vec{A}_{M}^{0})^{-1}\vec{A}_{M}^{1}\right)<1,
\]
where $\rho(\cdot)$ is the spectral radius of the matrix.
Then
$I-\left((U^{0}\vec{A}_{M}^{0})^{-\frac{1}{2}}U^{1}\vec{A}_{M}^{1}(U^{0}\vec{A}_{M}^{0})^{-\frac{1}{2}}\right)$ is positive definite, so the matrix $M_{M}^{t}$ is positive definite.
\qed\end{proof}

\subsection{Proof of Theorem \ref{thm:hyper}}
\label{p:thm:hyper}
\begin{proof}
 Lemmas \ref{lem:Dinv} and \ref{lem:Mpositive} show that
the matrix $\vec{B}_{M}^{0}=\vec{M}^{t}_{M}\vec{D}_{M}$ is invertible,
 and Theorem \ref{thm:eigmoment} implies that $\vec{B}_{M}$ is  diagonalizable with real eigenvalues and the spectral radius of $\vec{B}_{M}$ is less than one. The proof is completed.
\qed\end{proof}

\subsection{Proof of Theorem \ref{thm:Deg}}
\label{p:thm:Deg}
\begin{proof}
Because
\begin{align*}
    \nabla_{\vec{W}_{M}}\lambda_{i,M}=\frac{1}{\left(1-u   \hat{\lambda} _{i,M}\right)^{2}}\left(0,1- \hat{\lambda} _{i,M}^2,
    -(1-u^2)\frac{\partial  \hat{\lambda} _{i,M}}{\partial \theta},0\cdots,0\right)^{T},\
\end{align*}
and
$ \vec{r}_{i,M}=\vec{D}_{M}^{-1}\vec{P}_{M}^{p}\left((\vec{u}_{i,M})^{T},(\vec{v}_{i,M})^{T}\right)^{T}$,
$i=-M,\cdots M$,
one has
\begin{align}\label{EQ-lineargeneration}
\nabla_{\vec{W}_{M}}\lambda_{i,M}\cdot \vec{r}_{i,M}=\frac{1}
{\left(1-u  \hat{\lambda} _{i,M}\right)^{2}}\left(\left(1- \hat{\lambda} _{i,M}^2\right)
\vec{d}_{2}\vec{P}_{2}^{p}\vec{\tilde{r}}_{i}^{M}
-(1-u^2)\frac{\partial  \hat{\lambda} _{i,M}}{\partial \theta}\vec{d}_{3}\vec{P}_{2}^{p}\vec{\tilde{r}}_{i}^{M}\right),
\end{align}
where $\vec{\tilde{r}}_{i}^{M}=\left((\vec{u}_{i,M}^{(3)})^{T},(\vec{v}_{i,M}^{(2)})^{T}\right)^{T}$, $\vec{u}_{i,M}^{(3)}$ and $\vec{v}_{i,M}^{(2)}$ denote two vectors formed by  first three and two
components of  $\vec{u}_{i,M}$ and $\vec{v}_{i,M}$ respectively,
and $\vec{d}_{2}$ and $\vec{d}_{3}$ are the second and third row of $\vec{D}_{2}^{-1}$, specifically
    \begin{align*}
    \vec{d}_{2}=&-\frac{G(\zeta)(1-u^2)}{(\rho G(\zeta)+\Pi)\sqrt{\zeta}}\left(0,0,G(\zeta),0,\sqrt{-G(\zeta)^{2}+3\zeta^{-1}G(\zeta)+1}\right),
    \\
    \vec{d}_{3}=&\frac{1}{\rho\zeta(G(\zeta)^2\zeta^2-3G(\zeta)\zeta-\zeta^2+1)}\left(0,
    \frac{\sqrt{G(\zeta)-2\zeta^{-1}}(G(\zeta)^2\zeta^2-3G(\zeta)\zeta-\zeta^2+1)}{\sqrt{G(\zeta)^2\zeta^2-3G(\zeta)\zeta-\zeta^2+2}},
    \right.\\
    &\left.\frac{(G(\zeta)\zeta-1)\tilde{n}^{1}G(\zeta)}{\sqrt{\zeta}(\rho G(\zeta)+\Pi)},
    -\frac{\sqrt{2G(\zeta)^3\zeta^2-7G(\zeta)^2\zeta-2G(\zeta)\zeta^2+6G(\zeta)+\zeta}}{\sqrt{G(\zeta)^2\zeta^2-3G(\zeta)\zeta-\zeta^2+2}},\right.\\
    &\left.\frac{(G(\zeta)-\zeta^{-1})\sqrt{-G(\zeta)+3\zeta^{-1}G(\zeta)+1}}{\rho G(\zeta)+\Pi}\right).
    \end{align*}

 The identity \eqref{EQ-lineargeneration} always holds,
  because   $ \hat{\lambda} _{0,M}=0$ and $\vec{u}_{0,M}$ and $\vec{v}_{0,M}$ are given  in \eqref{EQ:3.29bbbb--zzz}.
The proof is completed.
\qed\end{proof}

\subsection{Explaination of Remark \ref{remark5.2}}
\label{p:remark5.2}
\begin{figure}
  \centering
  \begin{minipage}{6.5cm}
  \includegraphics[width=1\textwidth]{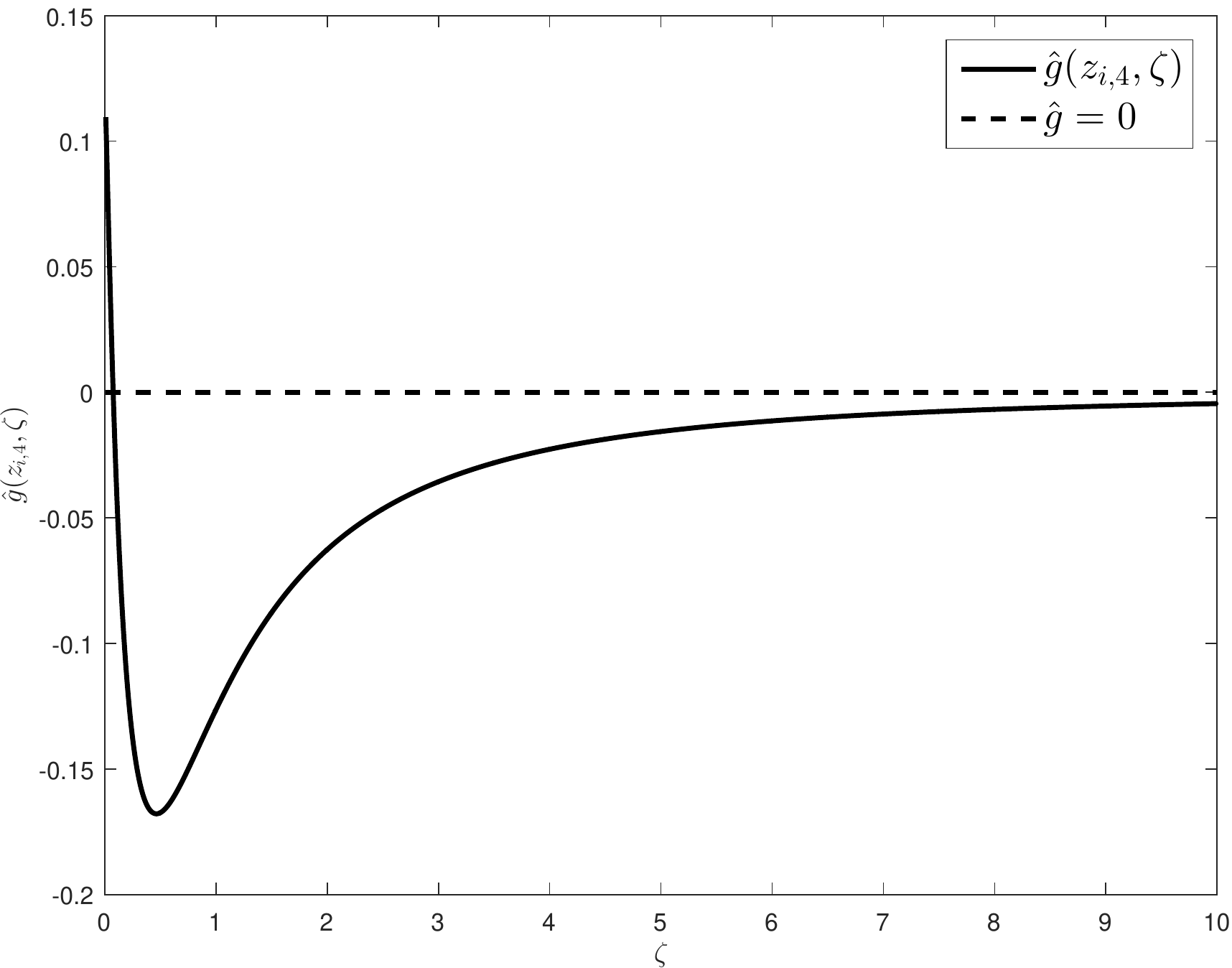}
  \end{minipage}
  \quad
  \begin{minipage}{6.5cm}
  \includegraphics[width=1\textwidth]{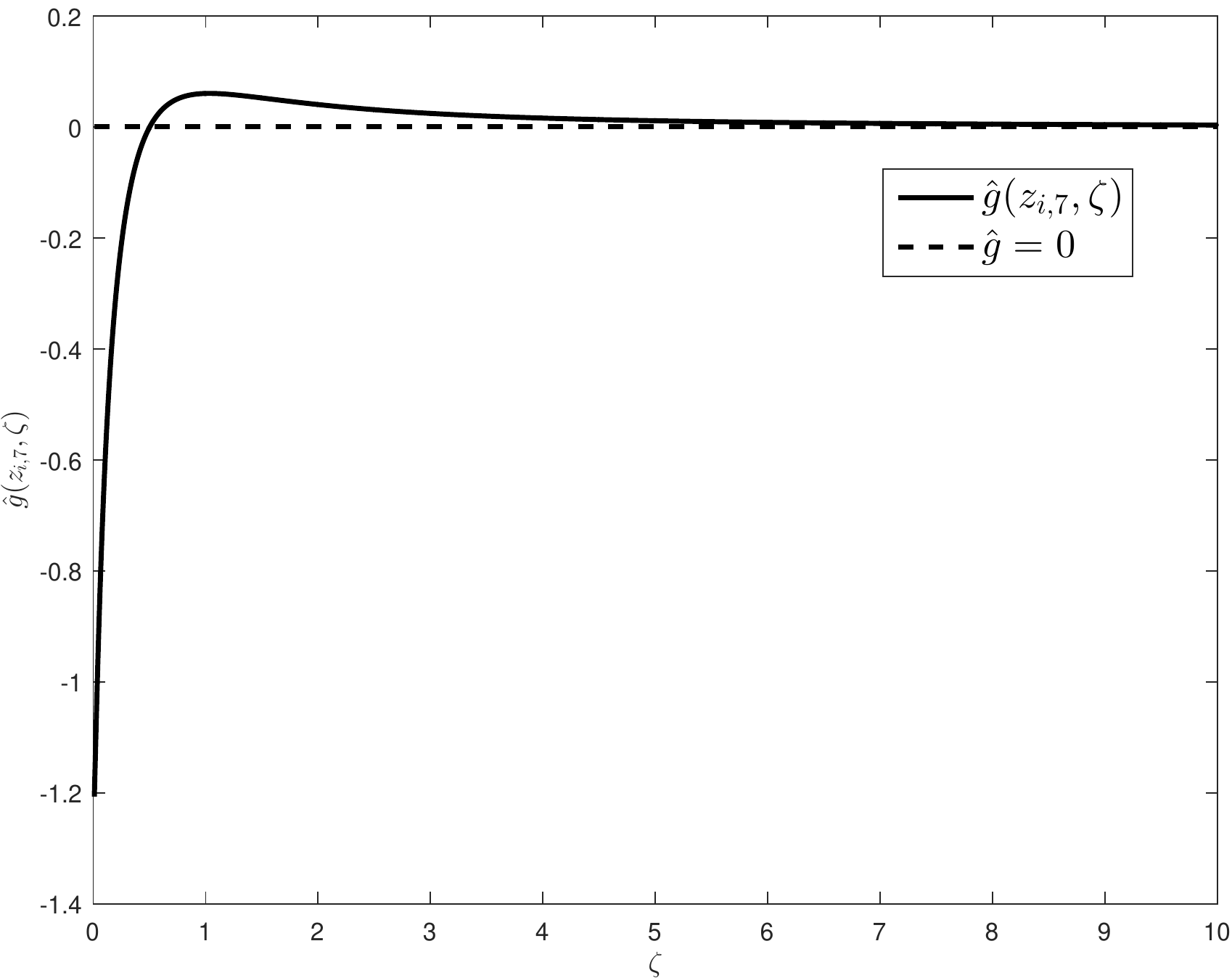}
  \end{minipage}
  \caption{\small Plots of   $\hat g(z_{1,M};\zeta)$  in terms of $\zeta$ for   $M=4$ and $7$ from the left to right. }
  \label{fig:notsimple}
\end{figure}

In fact, in order to judge by numerical experiments
whether  the sign of  $ \nabla_{\vec{W}_{M}}\lambda_{i,M}\cdot\vec{r}_{i,M}$  is constant or not,
 \eqref{EQ-lineargeneration} should be reformed.
For  $i=\pm 1,\pm 2,\cdots,\pm M$,  Theorem \ref{thm:eig} and \eqref{EQ-lineargeneration} give
\begin{align*}
\frac{    \nabla_{\vec{W}_{M}}\lambda_{i,M}\cdot\vec{r}_{i,M}} { -(1-u^2) }
=&\left(\frac{z_{i,M}^2-1}{\rho G(\zeta)+\Pi}
    -\zeta\frac{\partial z_{i,M}}{\partial \zeta}\frac{1}{\rho(G(\zeta)^2\zeta^2-3G(\zeta)\zeta-\zeta^2+1)}\right.\\
    &\cdot\left.\left(G(\zeta)\zeta-1-\zeta z_{i,M}-\frac{(G(\zeta)\zeta-1)\tilde{n}^{1}\sqrt{z_{i,M}^2-1}}{\rho G(\zeta)+\Pi}\right)\right)P_{M}^{(1)}(-z_{i,M};\zeta)
    \\
    &+\left(\frac{z_{i,M}^2-1}{\rho G(\zeta)+\Pi}
    +    \zeta\frac{\partial z_{i,M}}{\partial \zeta}\frac{1}{\rho(G(\zeta)^2\zeta^2-3G(\zeta)\zeta-\zeta^2+1)}\right.\\
    &\cdot \left.\left( G(\zeta)\zeta-1+\zeta z_{i,M} + \frac{(G(\zeta)\zeta-1)\tilde{n}^{1}\sqrt{z_{i,M}^2-1}}{\rho G(\zeta)+\Pi}\right)\right)P_{M}^{(1)}(z_{i,M};\zeta).
\end{align*}
Only a   simple case is  discussed in the following.
As shown in Remark \ref{remark2.2}, at  the local thermodynamic equilibrium, $\Pi=0$ and $n^\alpha=0$,
 thus one has
\begin{align*}
   \frac{ \nabla_{\vec{W}_{M}}\lambda_{i,M}\cdot\vec{r}_{i,M}  }{-(1-u^2)}
   =&\left(\frac{z_{i,M}^2-1}{\rho G(\zeta)}
    -\zeta\frac{\partial z_{i,M}}{\partial \zeta}
    \frac{G(\zeta)\zeta-1-\zeta z_{i,M} }{\rho(G(\zeta)^2\zeta^2-3G(\zeta)\zeta-\zeta^2+1)}\right)
    P_{M}^{(1)}(-z_{i,M};\zeta)
    \\
    &+\left(\frac{z_{i,M}^2-1}{\rho G(\zeta)}
    +\zeta\frac{\partial z_{i,M}}{\partial \zeta}
    \frac{G(\zeta)\zeta-1+\zeta z_{i,M} }
    {\rho(G(\zeta)^2\zeta^2-3G(\zeta)\zeta-\zeta^2+1)}\right)
     P_{M}^{(1)}(z_{i,M};\zeta).
\end{align*}
Using the term
$$\frac{\rho (G(\zeta)^2\zeta^2-3G(\zeta)\zeta-\zeta^2+1)}{(z_{i,M}^2-1)(G(\zeta)\zeta-1)  P_{M}^{(1)}(-z_{i,M};\zeta)},$$
to normalize the above identity
and noting  that
\[
{\rm sign}(G(\zeta)^2\zeta^2-3G(\zeta)\zeta-\zeta^2+1)=-{\rm sign}(x_{1,2}^{(0)}x_{2,2}^{(0)})<0
\]
gives
\[{\rm sign}\left(\nabla_{\vec{W}_{M}}\lambda_{i,M}\cdot\vec{r}_{i,M}\right)=(-1)^{M}{\rm sign}\left(\hat g(z_{i,M};\zeta)\right),
\]
where $\hat g(z_{i,M};\zeta)$ is defined by
\begin{align*}
\hat g(z_{i,M};\zeta)= &\frac{G(\zeta)^2\zeta^2-3G(\zeta)\zeta-\zeta^2+1}{ G(\zeta)\left(G(\zeta)\zeta-1\right)}\left(1+\frac{P_{M}^{(1)}(z_{i,M};\zeta)}{P_{M}^{(1)}(-z_{i,M};\zeta)}\right)\\
     &-\frac{\zeta}{(z_{i,M})^2-1} \left(1-\frac{P_{M}^{(1)}(z_{i,M};\zeta)}{P_{M}^{(1)}(-z_{i,M};\zeta)}\right)\frac{\partial z_{i,M}}{\partial \zeta}\\
     &+\frac{\zeta^2 z_{i,M}}{ \big( (z_{i,M})^2-1\big)  \left(G(\zeta)\zeta-1\right)}\left(1+\frac{P_{M}^{(1)}(z_{i,M};\zeta)}{P_{M}^{(1)}(-z_{i,M};\zeta)}\right)
     \frac{\partial z_{i,M}}{\partial \zeta}, \ i\geq 1,
\end{align*}
and   $\hat g(z_{i,M};\zeta):=\hat g(z_{-i,M};\zeta)$ for $i\leq -1$.
It is relatively easy to judge by numerical experiments
whether the sign of  $\hat g(z_{i,M};\zeta)$ is constant or not.
Fig. \ref{fig:notsimple} shows plots of  $\hat g(z_{1,4};\zeta)$  and  $\hat g(z_{1,7};\zeta)$ in terms of $\zeta$.
Similar to the special case of $M=4$ and 7,
our observation in numerical experiments is that  the sign of  $\hat g(z_{1,M};\zeta)$  is not constant  when $M\geq 4$ so that
 both $\lambda_{1,M}$ and  $\lambda_{-1,M}$  characteristic fields
 are neither linearly degenerate
 nor genuinely nonlinear when $M\geq 4$. 
 Such phenomenon is still not found in the case of $M\leq 3$.

\subsection{Proof of Theorem \ref{thm:LS}}
\label{p:thm:LS}
\begin{proof}
      Because
  the matrix $\vec{D}_M$ in \eqref{eq:Dmatrix} at  $\vec{W}_{M}=\vec{W}_{M}^{(0)}$
{can} be reformed as follows
\[
\vec{D}_{M}=\begin{pmatrix}
\vec{D}_{3\times3}^{11} &\vec{D}_{3\times2}^{12} & \vec{O}\\
\vec{O} & \vec{D}_{2\times2}^{22} & \vec{O}\\
\vec{O} & \vec{O} & \vec{I}_{2M-4}
\end{pmatrix},
\]
and its inverse  is given by
\[
\vec{D}_{M}^{-1}=\begin{pmatrix}
\left(\vec{D}_{3\times3}^{11}\right)^{-1} & -\left(\vec{D}_{3\times3}^{11}\right)^{-1}\vec{D}_{3\times2}^{12}\left(\vec{D}_{2\times2}^{22}\right)^{-1} & \vec{O}\\
\vec{O} & \left(\vec{D}_{2\times2}^{22}\right)^{-1} & \vec{O}\\
\vec{O} & \vec{O} & \vec{I}_{2M-4}
\end{pmatrix},
\]
as well as
\[
   \tilde{\vec{D}}_{M}^{W}=\begin{pmatrix}
   \vec{O} &\vec{D}_{3\times2}^{12} & \vec{O}\\
   \vec{O} & \vec{D}_{2\times2}^{22} & \vec{O}\\
   \vec{O} & \vec{O} & \vec{I}_{2M-4}
   \end{pmatrix},
   \]
 the product of $\tilde{\vec{D}}_{M}^{W}$ and $\vec{D}_{M}^{-1}$
 is of the following form
\[
\tilde{\vec{D}}_{M}^{W}\vec{D}_{M}^{-1}=\begin{pmatrix}
\vec{O}_{3\times3}& \vec{D}_{3\times2}^{12}\left(\vec{D}_{2\times2}^{22}\right)^{-1} & \vec{O}_{3\times(2M-4)}\\
\vec{O}_{2\times3}& \vec{I}_{2}& \vec{O}_{2\times(2M-4)}\\
\vec{O}_{(2M-4)\times3}&\vec{O}_{(2M-4)\times2}&\vec{I}_{2M-4}
\end{pmatrix},
\]
where $\vec{D}_{3\times3}^{11}$ is the $3\times 3$ subblock of $\vec{D}_{2}$ in the upper left corner, $\vec{D}_{3\times2}^{12}$ denotes the
$3\times 2$ subblock of $\vec{D}_{2}$ in the upper right corner,
and $\vec{D}_{2\times2}^{22}$ is $2\times 2$ subblock of $\vec{D}_{2}$ in the bottom right corner.
{It is obvious that each eigenvalue of $-\tilde{\vec{D}}_{M}^{W}\vec{D}_{M}^{-1}$ is non-positive, so does the matrix
\[
\bar{\vec{Q}}_{M}:=-\frac{1}{\tau}\left(U^{0}\vec{M}_{M}^{t}-U^{1}\vec{M}_{M}^{x}\right)^{\frac{1}{2}}\tilde{\vec{D}}_{M}^{W}\vec{D}_{M}^{-1}\left(U^{0}\vec{M}_{M}^{t}-U^{1}\vec{M}_{M}^{x}\right)^{-\frac{1}{2}}.
\]

The matrix $ U^{0}\vec{M}_{M}^{t}-U^{1}\vec{M}_{M}^{x} $ can be written as follows
\[
\begin{pmatrix}
\vec{M}_{3\times3}^{11}&\vec{M}_{3\times2}^{12}&\vec{O}_{3,2M-4}\\
(\vec{M}_{3\times2}^{12})^{T}&\vec{M}_{2\times2}^{22}&\vec{M}_{2\times(2M-4)}^{23}\\
\vec{O}_{2M-4,3}&(\vec{M}_{2\times(2M-4)}^{23})^{T}&\vec{M}_{(2M-4)\times(2M-4)}^{33}
\end{pmatrix},
\]
where $\vec{M}_{3\times3}^{11}$ is the $3\times 3$ subblock of $\vec{P}_{2}^{p}\vec{A}_{2}^{0}(\vec{P}_{2}^{p})^{T}$ in the upper left corner, $\vec{M}_{3\times2}^{12}$ denotes the
$3\times 2$ subblock of $\vec{P}_{2}^{p}\vec{A}_{2}^{0}(\vec{P}_{2}^{p})^{T}$ in the upper right corner,
and $\vec{M}_{2\times2}^{22}$ is $2\times 2$ subblock of $\vec{P}_{2}^{p}\vec{A}_{2}^{0}(\vec{P}_{2}^{p})^{T}$ in the bottom right corner,
the rest subblocks form the $(2M-2)\times (2M-2)$ bottom right corner of $\vec{P}_{2}^{p}\vec{A}_{M}^{0}(\vec{P}_{M}^{p})^{T}$.
Thus one has
\[
\vec{M}_{D}:=(\vec{M}_{3\times2}^{12})^{T}\vec{D}_{3\times2}^{12}\left(\vec{D}_{2\times2}^{22}\right)^{-1}
=-\left(\vec{D}_{3\times2}^{12}\left(\vec{D}_{2\times2}^{22}\right)^{-1}\right)^{T}\vec{M}_{3\times3}^{11}
\left(\vec{D}_{3\times2}^{12}\left(\vec{D}_{2\times2}^{22}\right)^{-1}\right),
\]
which is symmetric because
$
\vec{M}_{3\times3}^{11}\vec{D}_{3\times2}^{12}\left(\vec{D}_{2\times2}^{22}\right)^{-1}+
\vec{M}_{3\times2}^{12}=\vec{O}_{3\times2}.
$

On  the other hands, because the first three components of $\vec{S}(\vec{W}_M)$ are zero,
all elements in  the first three rows and  the first three columns of the matrix
\[
\vec{Q}_{M}=-\frac{1}{\tau}\left(U^{0}\vec{M}_{M}^{t}-U^{1}\vec{M}_{M}^{x}\right)\tilde{\vec{D}}_{M}^{W}\vec{D}_{M}^{-1},
\]
are zero, and the matrix $\vec{Q}_{M}$ is of form
\[
\vec{Q}_{M}=-\frac{1}{\tau}
\begin{pmatrix}
\vec{O}_{3,3}&\vec{O}_{3,2}&\vec{O}_{3,2M-4}\\
\vec{O}_{2,3}&\vec{M}_{2\times2}^{22}+\vec{M}_{D}&\vec{M}_{2\times(2M-4)}^{23}\\
\vec{O}_{2M-4,3}&(\vec{M}_{2\times(2M-4)}^{23})^{T}&\vec{M}_{(2M-4)\times(2M-4)}^{33}
\end{pmatrix}.
\]
Hence the matrix $\vec{Q}_{M}$  is symmetric.
It is obvious that $\vec{Q}_{M}$ is congruent with $\bar{\vec{Q}}_{M}$,
so it is negative semi-definite.

Because both matrices $\vec{D}_{M}$ and $\vec{M}_{M}^{t}$ are invertible, and
$\vec{M}_{M}^{t}$ is positive definite,
\eqref{eq:LS} is equivalent to
\begin{align}
 \det\left(i\omega \vec{I}-ik\vec{M}_{M}-\hat{\vec{Q}}_{M}\right)=0,
\end{align}
where
\begin{align*}
\hat{\vec{Q}}_{M}:=& \left(\vec{M}_{M}^{t}\right)^{-\frac{1}{2}}\vec{Q}_{M}\left(\vec{M}_{M}^{t}\right)^{-\frac{1}{2}},
\end{align*}
and
$$\vec{M}_{M}:=\left(\vec{M}_{M}^{t}\right)^{-\frac{1}{2}}\vec{M}_{M}^{x}\left(\vec{M}_{M}^{t}\right)^{-\frac{1}{2}}.$$
It is obvious that the matrix $\hat{\vec{Q}}_{M}$ is congruent with $\vec{Q}_{M}$ and negative semi-definite, and $\vec{M}_{M}$ is symmetric. Using Lemmas 1 and  2 in \cite{LS:2016} completes the
proof.}
\qed\end{proof}

\subsection{Proof of Lemma \ref{lem:lorentz}}
\label{p:lem:lorentz}
\begin{proof}
\begin{itemize}
\item[(i)] Under the given Lorentz boost ($x$ direction)
\[
t'=\gamma(v) (t-vx),\ x'=\gamma(v)(x-vt),\ \gamma(v)=(1-v^2)^{-\frac{1}{2}},
\]
where $v$ is the relative velocity between frames in the $x$-direction,
one has
	\begin{align*}
	&(p^{0})'=\gamma(v)(p^{0}-p^{1}v),\quad (p^{1})'=\gamma(v)(p^{1}-p^{0}v), \\
	&(U^{0})'=\gamma(v)(U^{0}-U^{1}v),\quad (U^{1})'=\gamma(v)(U^{1}-U^{0}v).
	\end{align*}
Thus one further {obtain}s
	\[
	E'= (U^{0})'(p^{0})'-(U^{1})'(p^{1})'=U^{0}p^{0}-U^{1}p^{1}=E,
	\]
	and
	\begin{align*}
	\left(\frac{p_{<1>}}{U^{0}}\right)'=& \frac{-(p^{<1>})'}{(U^{0})'}=-\frac{p^{<0>}-p^{<1>}v}{U^{1}-U^{0}v}=-\frac{(U^{0})^{-1}U^{1}p^{<1>}-p^{<1>}v}{U^{1}-U^{0}v}=\frac{p_{<1>}}{U^{0}},
	\\
	\left(\frac{dp}{p^{0}}\right)'=& \frac{d(p^{1})'}{(p^{0})'}=\frac{dp^{0}-dp^{1}v}{p^{1}-p^{0}v}=\frac{(p^{0})^{-1}p^{1}dp^{1}-dp^{1}v}{p^{1}-p^{0}v}=\frac{dp}{p^{0}}.
		\end{align*}
Combining them with  \eqref{eq:deff1} gives that each component of $\vec{f}_{M}$ is Lorentz invariant, such that
 the last $(2M-4)$ components of $\vec{W}_{M}$ are also Lorentz invariant.

	From \eqref{eq:variable0},
it is not difficult to prove that $\rho$ and $\theta$ are Lorentz invariant.
On the other hand, because
\[
\tilde{n}^{1}=\int_{\mathbb{R}}\frac{p^{<1>}}{U^{0}}f\frac{dp}{p^{0}},
\]
the quantity	 $\tilde{n}^{1}$ is  Lorentz invariant.
Moreover, on has
\[
\left(\frac{du}{1-u^2}\right)'=\frac{d(U^{1})'}{(U^{0})'}=\frac{dU^{0}-dU^{1}v}{U^{1}-U^{0}v}=\frac{(U^{0})^{-1}U^{1}dU^{1}-dU^{1}v}{U^{1}-U^{0}v}=\frac{dU^{1}}{U^{0}}=\frac{du}{1-u^2}.
\]
Using the above results completes the proof of the first part.

\item[(ii)]  Because   $\vec{A}_{M}^{0}$ and  $\vec{A}_{M}^{1}$  only depend on $\theta$,
they are Lorentz invariant.

The source term $\vec{S}(\vec{W}_{M})$ in \eqref{eq:colAWfinal}	{can} be rewritten into
	\[
	\vec{S}(\vec{W}_{M})=-\frac{1}{\tau}\vec{P}_{M}^{p}\vec{A}_{M}^{0}(\vec{P}_{M}^{p})^{T}\left(\vec{f}_{M}-\vec{f}^{(0)}_{M}\right),
	\]
	which has been expressed  in terms of 	the Lorentz covariant quantities.
	In fact,  the general source term $\vec{S}(\vec{W}_{M})$ in
	the  moment system   \eqref{eq:moment1} is also Lorentz invariant.
	The proof is completed.
\end{itemize}
	\qed\end{proof}

\subsection{Proof of Theorem \ref{thm:lorentz}}
\label{p:thm:lorentz}
\begin{proof}
	From  the 3rd step in Sec. \ref{subsec:deduction} and Lemma  \ref{lem:lorentz},  one  knows that
	$\hat{\vec{D}}_{M}=\vec{D}_M(\vec{D}_{M}^{u})^{-1}$ {can} be expressed in terms of the Lorentz covariant quantities, so it is Lorentz invariant.
	Because
	\begin{align*}
	(\vec{M}_{M}^{t})' &=-\gamma(v)(U^{1}-U^{0}v)\vec{P}_{M}^{p}\vec{A}_{M}^{1}(\vec{P}_{M}^{p})^{T}+\gamma(v)(U^{0}-U^{1}v)\vec{P}_{M}^{p}\vec{A}_{M}^{0}(\vec{P}_{M}^{p})^{T},\\
	 (\vec{M}_{M}^{x})'&=-\gamma(v)(U^{0}-U^{1}v)\vec{P}_{M}^{p}\vec{A}_{M}^{1}(\vec{P}_{M}^{p})^{T}+\gamma(v)(U^{1}-U^{0}v)\vec{P}_{M}^{p}\vec{A}_{M}^{0}(\vec{P}_{M}^{p})^{T},
	\end{align*}
	and
	\begin{align*}
	\left(\frac{\partial }{\partial t}\right)'=\gamma(v)\left(\frac{\partial }{\partial t}+v\frac{\partial }{\partial x}\right),\quad
	\left(\frac{\partial }{\partial x}\right)'=\gamma(v)\left(\frac{\partial }{\partial x}+v\frac{\partial }{\partial t}\right),
	\end{align*}
	%
	%
	%
	one has
	\begin{align*}
	\left(\vec{D}_{M}^{u}\frac{\partial \vec{W}_{M}}{\partial t}\right)'=&{\rm diag}\left\{1,\big((U^{0})^2\big)',1,\cdots,1\right\}\gamma(v)\left(\frac{\partial (\rho,u',\theta,\Pi,\tilde{n}^{1},f_{3}^{0},\cdots,f_{M-1}^{1})^{T}}{\partial t}\right.\\
&	+\left.v\frac{\partial (\rho,u',\theta,\Pi,\tilde{n}^{1},f_{3}^{0},\cdots,f_{M-1}^{1})^{T}}{\partial x}\right)\\
	=&{\rm diag}\left\{1,\big((U^{0}\big)^{-1})',1,\cdots,1\right\}\gamma(v)\left(\frac{\partial (\rho,(U^{1})',\theta,\Pi,\tilde{n}^{1},f_{3}^{0},\cdots,f_{M-1}^{1})^{T}}{\partial t}\right.\\
&	+\left.v\frac{\partial (\rho,(U^{1})',\theta,\Pi,\tilde{n}^{1},f_{3}^{0},\cdots,f_{M-1}^{1})^{T}}{\partial x}\right)\\
	=&\vec{D}_{M}^{u}\gamma(v)\left(\frac{\partial \vec{W}_{M}}{\partial t}+v\frac{\partial \vec{W}_{M}}{\partial x}\right),
	\end{align*}
	where  the last equal sign is derived by  following the proof of Lemma \ref{lem:lorentz}.
	Similarly, one has
	\begin{equation*}
	\left(\vec{D}_{M}^{u}\frac{\partial \vec{W}_{M}}{\partial x}\right)'=\vec{D}_{M}^{u}\gamma(v)\left(\frac{\partial \vec{W}_{M}}{\partial x}+v\frac{\partial \vec{W}_{M}}{\partial t}\right).
	\end{equation*}
	Thus one {obtains}
	\begin{align*}
	&\left(\vec{B}_{M}^{0}\frac{\partial \vec{W}_{M}}{\partial t} +  \vec{B}_{M}^{1}\frac{\partial \vec{W}_{M}}{\partial x}\right)'\\
	=&(\vec{M}_{M}^{t})'\left(\vec{D}_{M}\frac{\partial \vec{W}_{M}}{\partial t}\right)'+(\vec{M}_{M}^{x})'\left(\vec{D}_{M}\frac{\partial \vec{W}_{M}}{\partial x}\right)'\\
	=&\left(-(U^{1})'\vec{P}_{M}^{p}\vec{A}_{M}^{1}(\vec{P}_{M}^{p})^{T}+ (U^{0})'P_{M}^{p}\vec{A}_{M}^{0}(P_{M}^{p})^{T}\right)\vec{D}_{M}\left(\gamma(v)\left(\frac{\partial \vec{W}_{M}}{\partial t}+v\frac{\partial \vec{W}_{M}}{\partial x}\right)\right)\\
	&+\left(-(U^{0})'\vec{P}_{M}^{p}\vec{A}_{M}^{1}(\vec{P}_{M}^{p})^{T}+ (U^{1})'\vec{P}_{M}^{p}\vec{A}_{M}^{0}(\vec{P}_{M}^{p})^{T}\right)\vec{D}_{M}\left(\gamma(v)\left(\frac{\partial \vec{W}_{M}}{\partial x}+v\frac{\partial \vec{W}_{M}}{\partial t}\right)\right)\\
	=&\left(-U^{1}v\vec{P}_{M}^{p}\vec{A}_{M}^{1}(\vec{P}_{M}^{p})^{T}+U^{0}\vec{P}_{M}^{p}\vec{A}_{M}^{0}(\vec{P}_{M}^{p})^{T}\right)
	\vec{D}_{M}\frac{\partial \vec{W}_{M}}{\partial t}\\
	&+\left(-U^{0}\vec{P}_{M}^{p}\vec{A}_{M}^{1}(\vec{P}_{M}^{p})^{T}+U^{1}\vec{P}_{M}^{p}\vec{A}_{M}^{0}(\vec{P}_{M}^{p})^{T}\right)
	\vec{D}_{M}\frac{\partial \vec{W}_{M}}{\partial x}\\
	=&\vec{B}_{M}^{0}\frac{\partial \vec{W}_{M}}{\partial t} +  \vec{B}_{M}^{1}\frac{\partial \vec{W}_{M}}{\partial x}.
	\end{align*}
	Combining it with Lemma \ref{lem:lorentz}  completes the proof.
	\qed\end{proof}

\section{Proofs in Section \ref{sec:NE}}
\label{sec:App6}
\subsection{Proof of Lemma \ref{lem:projectn}}
\label{p:lem:projectn}
\begin{proof}
Using Lemmas \ref{lem:basis} and \ref{lem:project} gives
	\begin{align*}
	 \tilde{f}g^{(0)}_{[u,\theta]}=&\sum_{i=0}^{M}\tilde{f}_{i}^{0}\tilde{P}_{i}^{(0)}[u,\theta]+\sum_{j=0}^{M-1}\tilde{f}_{j}^{1}\tilde{P}_{j}^{(1)}[u,\theta],\\
\Pi_{M}[u,\theta]	f=&\sum_{i=0}^{M}f_{i}^{0}\tilde{P}_{i}^{(0)}[u,\theta]+\sum_{j=0}^{M-1}f_{j}^{(1)}\tilde{P}_{j}^{(1)}[u,\theta],
	\end{align*}
	where
	\begin{align*}
	f_{i}^{0}&=<f,\tilde{P}_{i}^{(0)}[u,\theta]>_{g^{(0)}_{[u,\theta]}},\quad \tilde{f}_{i}^{0}=<\tilde{f}g^{(0)}_{[u,\theta]},\tilde{P}_{i}^{(0)}[u,\theta]>_{g^{(0)}_{[u,\theta]}},\quad i\leq M,\\
	f_{j}^{1}&=<f,\tilde{P}_{j}^{(1)}[u,\theta]>_{g^{(0)}_{[u,\theta]}},\quad
	\tilde{f}_{j}^{1}=<\tilde{f}g^{(0)}_{[u,\theta]},\tilde{P}_{j}^{(1)}[u,\theta]>_{g^{(0)}_{[u,\theta]}},\quad j\leq M-1.
	\end{align*}
Therefore one has
	\begin{align*}
	&<\tilde{f}g^{(0)}_{[u,\theta]},\Pi_{M}[u,\theta]f>_{g^{(0)}_{[u,\theta]}}=<\tilde{f}f,\Pi_{M}[u,\theta]f>_{f}\\
	=&\sum_{i=0}^{M}f_{i}^{0}<\tilde{f}f,\tilde{P}_{i}^{(0)}[u,\theta]>_{f}+
	\sum_{j=0}^{M-1}f_{j}^{1}<\tilde{f}f,\tilde{P}_{j}^{(1)}[u,\theta]>_{f}\\
	=&\sum_{i=0}^{M}f_{i}^{0}<\tilde{f}g^{(0)}_{[u,\theta]},\tilde{P}_{i}^{(0)}[u,\theta]>_{g^{(0)}_{[u,\theta]}}+
	\sum_{j=0}^{M-1}f_{j}^{1}<\tilde{f}g^{(0)}_{[u,\theta]},\tilde{P}_{j}^{(1)}[u,\theta]>_{g^{(0)}_{[u,\theta]}}\\
	=&\sum_{i=0}^{M}f_{i}^{0}\tilde{f}_{i}^{0}+
	\sum_{j=0}^{M-1}f_{j}^{1}\tilde{f}_{j}^{1}\\
	=&\sum_{i=0}^{M}<f,\tilde{f}_{i}^{0}\tilde{P}_{i}^{(0)}[u,\theta]>_{g^{(0)}_{[u,\theta]}}+
	\sum_{j=0}^{M-1}<f,\tilde{f}_{j}^{1}\tilde{P}_{j}^{(1)}[u,\theta]>_{g^{(0)}_{[u,\theta]}}\\
	=&<f,\tilde{f}g^{(0)}_{[u,\theta]}>_{g^{(0)}_{[u,\theta]}}=<f,\tilde{f}f>_{f}.
	\end{align*}
	The proof is completed.
	\qed\end{proof}

\subsection{Proof of Lemma \ref{lem:ut}}
\label{p:lem:ut}
\begin{proof}
Using Lemma \ref{lem:project} gives
	\begin{align*} &\Pi_{M}[u_{1},\theta_{1}]f=\sum_{i=0}^{M}f_{i}^{0}\tilde{P}_{i}^{(0)}[u_{1},\theta_{1}]+
\sum_{j=0}^{M-1}f_{j}^{1}\tilde{P}_{i}^{(1)}[u_{1},\theta_{1}],\\ &\Pi_{M}[u_{1},\theta_{1}]\Pi_{M}[u_{2},\theta_{2}]f=\sum_{i=0}^{M}\tilde{f}_{i}^{0}\tilde{P}_{i}^{(0)}[u_{1},\theta_{1}]+
\sum_{j=0}^{M-1}\tilde{f}_{j}^{1}\tilde{P}_{j}^{(1)}[u_{1},\theta_{1}],
	\end{align*}
	\begin{align*}
	f_{i}^{0}&=<f,\tilde{P}_{i}^{(0)}[u_{1},\theta_{1}]>_{g^{(0)}_{[u_{1},\theta_{1}]}}=
	{<f,{P}_{i}^{(0)}(u_{1},\zeta_{1})f>_{f}},i\leq M,\\
	f_{j}^{1}&=<f,\tilde{P}_{j}^{(1)}[u_{1},\theta_{1}]>_{g^{(0)}_{[u_{1},\theta_{1}]}}=
	{<f,{P}_{j}^{(1)}(u_{1},\zeta_{1})(U_{0})_{1}^{-1}p_{<1>}f>_{f}},j\leq M-1,
	\end{align*}
	\begin{align*}
	\tilde{f}_{i}^{0}&=<\Pi_{M}[u_{2},\theta_{2}]f,\tilde{P}_{i}^{(0)}[u_{1},\theta_{1}]>_{g^{(0)}_{[u_{1},\theta_{1}]}}=
	{<\Pi_{M}[u_{2},\theta_{2}]f,{P}_{i}^{(0)}(u_{1},\zeta_{1})f>_{f}},i\leq M,\\
	\tilde{f}_{j}^{1}&=<\Pi_{M}[u_{2},\theta_{2}]f,\tilde{P}_{j}^{(1)}[u_{1},\theta_{1}]>_{g^{(0)}_{[u_{1},\theta_{1}]}}=
	{<\Pi_{M}[u_{2},\theta_{2}]f,{P}_{j}^{(1)}(u_{1},\zeta_{1})(U_{0})_{1}^{-1}p_{<1>}f>_{f}},j\leq M-1,
	\end{align*}
Because both ${P}_{i}^{(0)}(u_{1},\zeta_{1})f$ and ${P}_{j}^{(1)}(u_{1},\zeta_{1})(U_{0})_{1}^{-1}p_{<1>}f$ belong to
the space $\mathbb{H}_{M}^{f}$,
	using Lemma \ref{lem:projectn} {completes} the proof.
	\qed\end{proof}

\subsection{Proof of Theorem \ref{thm:colS}}
\label{p:thm:colS}
\begin{proof}
Because Eq.	\eqref{eq:semicol} is equivalent to
	\[
	\left(\vec{I}+\frac{\Delta t}{\tau_{i}^{\ast}}\left(\vec{M}_{i,M}^{t\ast}\right)^{-1}\left(U_{i}^{0\ast}\vec{M}_{i,M}^{t\ast}-U_{i}^{1\ast}
	\vec{M}_{i,M}^{x\ast}\right)\left(\vec{I}-\vec{D}_{i,M}^{f_{i}^{(0)\ast}}\right)\right)
 \vec{f}_{i,M}^{n+1}=\vec{f}_{i,M}^{\ast},
	\]
it is unconditionally stable if and only if the modulus  of each eigenvalue of the matrix
	\[
	\vec{I}+\frac{\Delta t}{\tau_{i}^{\ast}}\left(\vec{M}_{i,M}^{t\ast}\right)^{-1}\left(U_{i}^{0\ast}\vec{M}_{i,M}^{t\ast}-U_{i}^{1\ast}
	\vec{M}_{i,M}^{x\ast}\right)\left(\vec{I}-\vec{D}_{i,M}^{f_{i}^{(0)\ast}}\right),
	\]
	is not less than one.
It is true if the real part of each eigenvalue of the matrix
\begin{equation}
\label{EQ:6.700000-THZ} \left(\vec{M}_{i,M}^{t\ast}\right)^{-1}\left(U_{i}^{0\ast}\vec{M}_{i,M}^{t\ast}-U_{i}^{1\ast}\vec{M}_{i,M}^{x\ast}\right)
\left(\vec{I}-\vec{D}_{i,M}^{f_{i}^{(0)\ast}}\right)=:\left(\vec{M}_{i,M}^{t\ast}\right)^{-1}\bar{\vec{M}}_{D}^{*},
\end{equation}
	is non-negative.

{In fact,
thanks to \eqref{EQ.6.3b}, the characteristic polynomial of
the upper triangular matrix $ \vec{I}-\vec{D}_{i,M}^{f_{i}^{(0)\ast}}$
is explicitly given by
\begin{align*}
	0=\det\left(\lambda \vec{I}-\left(\vec{I}-\vec{D}_{i,M}^{f_{i}^{(0)\ast}}\right)\right)=
\det\left((\lambda-1)\vec{I}+\vec{D}_{i,M}^{f^{(0)\ast}}\right)=\lambda(\lambda-1)^{2M},
	\end{align*}
and
 $\bar{\vec{M}}_{D}^{*}$ is a symmetric   matrix
and   congruent with
\[
\left(U_{i}^{0\ast}\vec{M}_{i,M}^{t\ast}-U_{i}^{1\ast}\vec{M}_{i,M}^{x\ast}\right)^{\frac{1}{2}}
\left(\vec{I}-\vec{D}_{i,M}^{f_{i}^{(0)\ast}}\right)\left(U_{i}^{0\ast}
\vec{M}_{i,M}^{t\ast}-U_{i}^{1\ast}\vec{M}_{i,M}^{x\ast}\right)^{-\frac{1}{2}},
\]
which is similar to the  matrix $ \vec{I}-\vec{D}_{i,M}^{f_{i}^{(0)\ast}}$.
Thus the matrix $\bar{\vec{M}}_{D}^{*}$ is positive semi-definite and each eigenvalue of the matrix
$
\left(\vec{M}_{i,M}^{t\ast}\right)^{-1}\bar{\vec{M}}_{D}^{\ast}$ is non-negative
because of the relation
$
\left(\vec{M}_{i,M}^{t\ast}\right)^{-1}\bar{\vec{M}}_{D}^{\ast}
=\left(\vec{M}_{i,M}^{t\ast}\right)^{-\frac12}  \big( \left(\vec{M}_{i,M}^{t\ast}\right)^{-\frac12}  \bar{\vec{M}}_{D}^{\ast}
\left(\vec{M}_{i,M}^{t\ast}\right)^{-\frac12}\big) \left(\vec{M}_{i,M}^{t\ast}\right)^{\frac12}$
. }
	The proof is completed.
	\qed\end{proof}

\end{appendices}

\end{document}